\documentclass[11pt]{article}

\usepackage[margin=1in]{geometry}
\usepackage{authblk}
\usepackage{graphicx}
\usepackage{xcolor}
\usepackage{natbib}
\usepackage{hyperref}

\usepackage{amsmath,amssymb,amsthm,mathtools,bm,amsfonts}
\usepackage{booktabs}

\usepackage[ruled,vlined,algo2e]{algorithm2e}

\theoremstyle{plain}
\newtheorem{theorem}{Theorem}[section]

\newtheorem{proposition}[theorem]{Proposition}
\newtheorem{corollary}[theorem]{Corollary}

\theoremstyle{definition}
\newtheorem{definition}{Definition}[section]
\newtheorem{assumption}{Assumption}[section]

\theoremstyle{remark}

\newtheorem{remark}[theorem]{Remark}

\newcommand{\M}{\mathcal{M}}

\newcommand{\X}{\mathcal{X}}
\newcommand{\popsi}{\tilde{\psi}}
\newcommand{\pssig}{\xi^{\mathrm{S}}}
\newcommand{\SigmaS}{\widetilde\Sigma^{\mathrm{S}}}
\newcommand{\rank}{\operatorname{rank}}
\newcommand{\op}{\operatorname{op}}
\newcommand{\Tr}{\operatorname{Tr}}
\newcommand{\Var}{\operatorname{Var}}

\newcommand{\E}{\mathbb{E}}
\newcommand{\R}{\mathbb{R}}

\newcommand{\samp}{\mathrm{samp}}

\DeclareMathOperator{\Cov}{Cov}
\DeclareMathOperator{\diag}{diag}

\SetKw{KwReturn}{return}

\AtBeginDocument{
  \setlength{\abovedisplayskip}{4pt}  
  \setlength{\belowdisplayskip}{4pt}   
  \setlength{\abovedisplayshortskip}{2pt}
  \setlength{\belowdisplayshortskip}{3pt}
}

\newcommand{\SecTech}{Section~S1}                
\newcommand{\SecCondCLT}{Section~S1.1}           
\newcommand{\SecSecularEq}{Section~S1.2}

\newcommand{\LemBdeloc}{Lemma~S1.1}              
\newcommand{\LemCondCLT}{Lemma~S1.2}             
\newcommand{\LemIsotropy}{Lemma~S1.3}            
\newcommand{\PropDecouple}{Proposition~S1.4}     
\newcommand{\EqSigmaplugin}{Equation~(S20)}

\newcommand{\SecDeltaMethod}{Section~S1.4}       
\newcommand{\EqPiCLT}{Equation~(S18)}            
\newcommand{\EqSigmaPiMult}{Equation~(S19)}      
\newcommand{\MainTableone}{S1}      
\newcommand{\MainTabletwo}{S2}      
\newcommand{\Mainfigp}{S5}      
\newcommand{\Mainpathwaylow}{S8}

\def\spacingset#1{\renewcommand{\baselinestretch}%
{#1}\small\normalsize} 
\title{\sc {Inference for Similarity and Alignability between Noisy High-Dimensional Datasets}}

\author[1]{Hongrui Chen}
\author[1,2,3]{Rong Ma\thanks{Corresponding author: \texttt{rongma@hsph.harvard.edu}}}

\affil[1]{Department of Biostatistics, Harvard University}

\affil[2]{Eric and Wendy Schmidt Center, Broad Institute}

\affil[3]{Department of Data Science, Dana-Farber Cancer Institute}

\begin{document}

\maketitle

\begin{abstract}
The rapid growth of high-dimensional datasets across a wide range of scientific domains has created an urgent need for new statistical methods to compare distributions with underlying low-dimensional structure. Assessing similarity between high-dimensional datasets whose observations concentrate near low-dimensional manifolds is particularly challenging due to the nontrivial effects of noise in high dimensions. We propose a principled framework for statistical inference on the similarity and alignability of high-dimensional datasets with low-dimensional smooth signals under heterogeneous noise. The key idea is to link the spectral properties of the observed data matrices to the geometry of their underlying signal distributions. Under a manifold signal-plus-noise model, we build on the principal variances associated to the underlying signals and develop a scale- and rotation-invariant dissimilarity measure between two datasets that may differ in sample size and noise structure. We further construct an estimator of the dissimilarity and a statistical test of dataset alignability, namely, whether the dissimilarity vanishes. The proposed methodology and its theoretical guarantees under high-dimensional asymptotic settings draw on recent advances in random matrix theory (RMT). The proposed framework accommodates heterogeneous noise across datasets and provides a fast,
theoretically grounded approach to comparing high-dimensional datasets
with low-dimensional  structures. Through extensive simulations
and analyses of multiple single-cell datasets, we demonstrate that the proposed method substantially outperforms existing approaches.
\end{abstract}

\noindent\textbf{Keywords:}
Data integration; High-dimensional statistics; Spectral methods.

\section{Introduction}

\begin{figure}[htbp]
  \centering
  \begin{minipage}{0.48\linewidth}
    \centering
    \includegraphics[width=\linewidth]{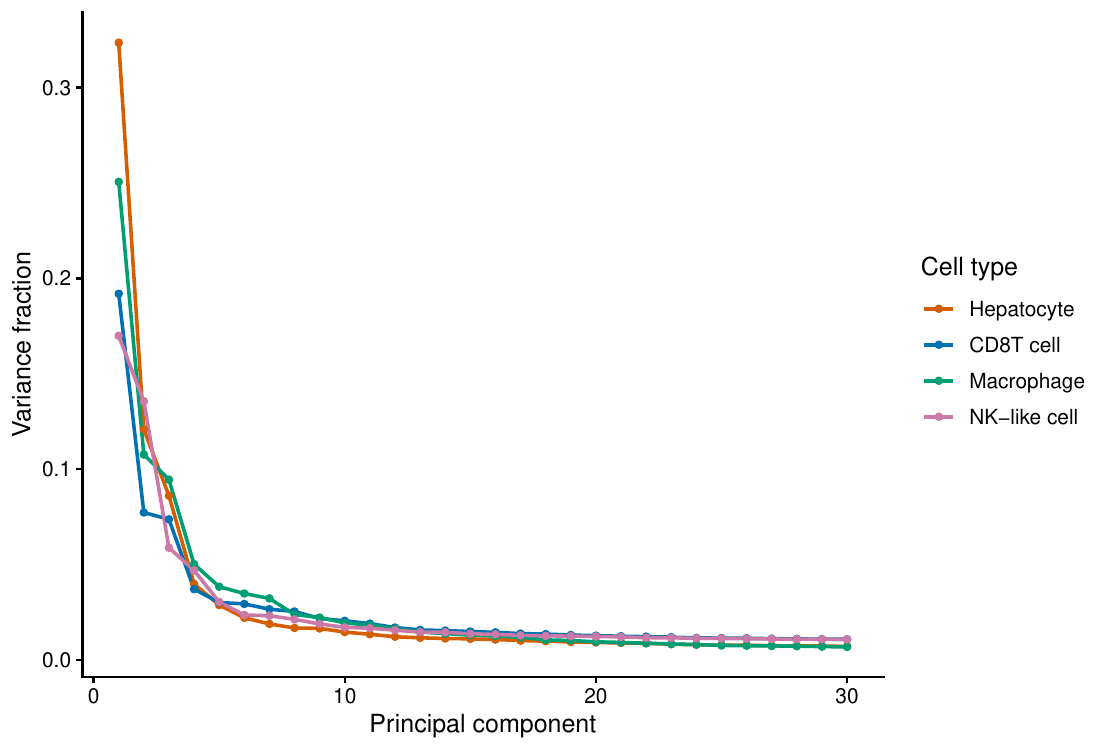}
    \\ \small (a)
  \end{minipage}
  \hfill
  \begin{minipage}{0.48\linewidth}
    \centering
    \includegraphics[width=\linewidth]{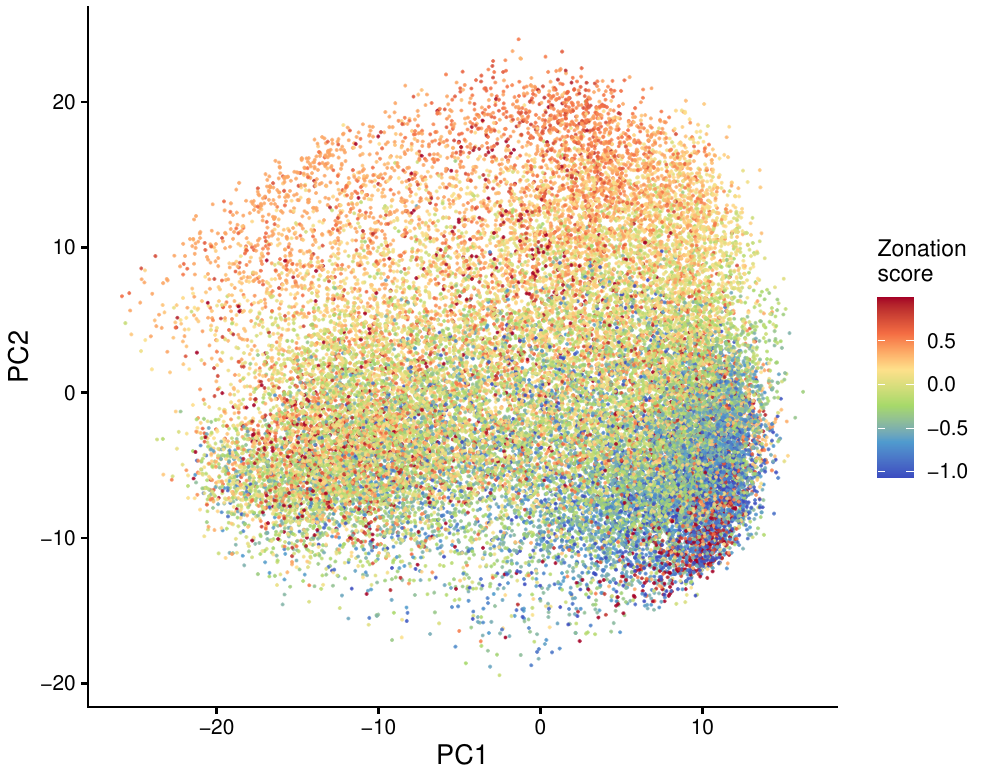}
    \\ \small (b)
  \end{minipage}

  \caption{}
  \label{fig:bioassumption}
\end{figure}

Recent advances in data collection, storage and computation have made large, high-dimensional datasets increasingly available across domains such as astronomy, business analytics, human genetics, and microbiology. Despite their high ambient dimensionality, the underlying signal structure of many datasets can often be effectively modeled and represented by some ideal low-dimensional structure, such as a smooth manifold. This ``manifold hypothesis" \citep{fefferman2016testing} forms the foundation of important fields of data science, such as manifold learning \citep{meilua2024manifold}, domain adaptation \citep{luo2020unsupervised} and representation learning \citep{bengio2013representation}. 
Single-cell data provide a particularly relevant example of this
setting. Although each cell is measured across thousands of genes, coordinated
gene programs often concentrate biologically meaningful variation in a small
number of spectral directions. This motivates low-rank models for single-cell
RNA-seq data \citep{Linderman2022ALRA,Aparicio2020RMT}. Continuous cellular
processes may also form nonlinear trajectories, supporting a low-dimensional
manifold representation \citep{Verma2020Manifold}. For example,
Figure~\ref{fig:bioassumption}(a) shows rapid spectral decay in single-cell RNA-seq datasets across four cell
types, while Figure~\ref{fig:bioassumption}(b) shows a continuous zonation
gradient among HCC hepatocytes. These observations motivate the use of a manifold signal-plus-noise model, previously introduced in \citep{wu2018think,DingMa2023,ding2025kernel} and adopted in this work (Section 2.1).

\color{black}

When analyzing multiple high-dimensional datasets that individually contains some underlying low-dimensional signal structures, it is often crucial to characterize and quantify the (dis)similarity between these datasets in terms of their latent signal distributions. Such analysis enables assessment of the diversity of information contained in these noisy datasets and the extent to which their underlying signals overlap, offering an integrated perspective that reveals high-level relationships among datasets that are not apparent when each dataset is analyzed in isolation.

A closely related problem is testing whether two datasets share the same latent structure, which calls for rigorous and broadly applicable two-sample procedures. Classical parametric tests such as Hotelling’s $T^2$ and Student’s $t$, as well as many nonparametric methods—including energy distance \citep{energydistance1,energydistance2}, maximum mean discrepancy (MMD) \citep{kernelmmd1,kernelmmd3}, graph-based tests \citep{graphtest1,graphtest2}, Wasserstein distance and its projected variants \citep{wasserstein1}, and interpoint distance tests \citep{interpoint2}—were primarily developed in finite-dimensional settings. More fundamentally, most of these approaches compare the entire ambient-space distributions, which may provide an unnecessarily fine resolution for assessing shared latent structure: modest geometric, distributional, or nuisance variation can make two datasets appear substantially different even when their dominant structures are similar. This issue is compounded in high dimensions, where distributional estimation and test calibration become increasingly difficult. 

For comparing  high-dimensional datasets, a more effective alternative approach is to compare \emph{structural properties} of their underlying distributions, such as feature covariance. A prominent line of work focuses on  testing large-scale multiple hypotheses about the entry-wise equality of the sparse population covariance or precision matrices between two datasets \citep{cai2013two,xia2015testing,xia2018multiple,liu2013gaussian}. When the data matrices appear to exhibit low-rank structures, another group of work focus on  the inference of the leading eigenvalues or eigensubspaces of the population covariance matrices \citep{smai,rmt1,mestre2011asymptotic,yan2024inference,bao2021singular}, distance covariance matrices \citep{li2023eigenvalue}, or canonical correlation matrices \citep{bao2019canonical,yang2022limiting}, which can be used to detect structural similarity or differences between two high-dimensional datasets. 
Despite their conceptual simplicity and effectiveness in comparing high-dimensional datasets, existing methods have not exploited the smoothness of the low-dimensional signal structure often present in applications such as single-cell analysis, thereby missing an important aspect to improve inferential power.
Moreover, most existing approaches assume homogeneous noise both within and across datasets, making them vulnerable to inferential bias and  inconsistency when applied to datasets with complex, heteroskedastic noise. 

To address these limitations, in this work we consider a general manifold signal-plus-noise model and propose a statistical framework for inferring similarity and alignability between high-dimensional datasets with  latent low-dimensional signal structures. Our approach differs from existing eigenstructure inference methods in three key aspects: (i) it builds on a distinct statistical model to enable designing statistical procedures with enhanced power for inference on datasets with low-dimensional \emph{smooth} structures;
(ii) it leverages recent advances in random matrix theory to handle and overcome the effects of \emph{high-dimensional heteroskedastic noise} across datasets, thereby enabling consistent inference in high-dimensional asymptotic setting under data heterogeneity;  and (iii) it introduces a continuous \emph{dissimilarity measure} between datasets based on their population-level signal profiles, offering a practical and interpretable characterization of similarity between noisy datasets. 

Suppose we observe two data matrices $Y_1$ and $Y_2$, generated from the manifold signal-plus-noise model:
\begin{equation}\label{eq:signal_plus_noise}
Y_k \;=\; S_{\samp,k} \;+\; \Sigma_k^{1/2} X_k \;\in\; \mathbb{R}^{p\times N_k},\quad k=1,2
\end{equation}
where $S_{\samp,k}=[S_{1,k},\ldots,S_{N_k,k}]$ collects $N_k$ noiseless samples (signals) drawn from a low-dimensional smooth manifold $\mathcal{M}_k$ embedded in $\mathbb{R}^p$ (see Section~\ref{subsec:22} for a more formal description), $\Sigma_k\in\R^{p\times p}$ is a deterministic matrix encoding the noise covariance, and $X_k=[x_{i\mu,k}]\in \mathbb{R}^{p\times N_k}$ has independent entries with mean zero and unit variance. For the ambient space dimensionality $p$, we will consider the high-dimensional regime where $p,N_k\to\infty$ with $p/N_k\to\varphi_k\in(0,\infty)$. 
{Unlike traditional homoskedastic models, we assume a piecewise diagonal noise covariance (see Assumption 2.1 below). This assumption is motivated by considerations in single-cell applications. In many single-cell datasets, the variances of gene expression across cells often exhibit group structures that reflect distinct genetic programs comprising functionally related genes. Genes within the same program tend to share similar variance profiles, whereas those belonging to different programs display distinct patterns;
see Section~S6.2 in the Supplement for empirical evidence  supporting this assumption.} When this structure is approximately valid, pooling information across features yields more stable noise estimates and can improve the statistical power of our method.
\color{black}{}
In addition, we assume the data is already centered and define the sample covariance $Q_k:=Y_kY_k^\top/N_k$, whose population counterpart is $
\widetilde\Sigma_k \;:=\; \E[Q_k] \;=\; \Sigma_k + M_k,$ where $M_k \;:=\; \frac{1}{N_k}\E[S_{\samp,k} S_{\samp,k}^\top].$
Our goal is to conduct statistical inference on the similarity and alignability of the latent signal structures underlying $Y_1$ and $Y_2$. For similarity inference, we define a dissimilarity measure between the corresponding signal distributions and develop an estimator of this quantity, together with uncertainty quantification, based on the observed high-dimensional noisy data. For alignability inference, we construct a rigorous statistical test of whether the dissimilarity between the two latent signal structures vanishes under the proposed metric.

Our key innovation lies in establishing a theoretical connection between the signal distribution on the manifold $\mathcal{M}_k$, and the spectral properties of matrices $\widetilde{\Sigma}_k$, $\Sigma_k$, and $M_k$, which are derived from finite-sample observations. This linkage enables the formulation of a natural dissimilarity measure between noisy, high-dimensional datasets based on their latent population signatures, effectively capturing discrepancies between their underlying low-dimensional signal structures. By bridging population-level geometric properties with finite-sample spectral behavior, this framework provides a foundation for developing consistent inference procedures for manifold learning grounded in RMT.
Our main methodological contributions are:
\begin{itemize}
\item We introduce the normalized Manifold Spectral Dissimilarity (nMSD), a scale-invariant metric that quantifies the difference between the population-level low-dimensional signal structures of two datasets. We further propose a computationally efficient estimator of the nMSD based on noisy observations under heterogeneous noise and unequal sample sizes, and establish its consistency and asymptotic normality in high-dimensional settings, thereby providing theoretical justification for dataset similarity inference.
\item We formulate the data alignability problem as a  two-sample hypothesis test on the equivalence of latent signal structures with respect to the nMSD metric. A corresponding statistical test is proposed, and its asymptotic validity and power consistency are established under high-dimensional asymptotic settings.
\item As a important component of our inferential framework, we develop mean-squared consistent estimators for data-specific noise covariance matrices, allowing us to accurately account for the effects of group-wise heteroskedastic noise. Our approach integrates the low-rank truncation technique with a Potts segmentation method, which can be efficiently solved using dynamic programming.
 \item We extend the concept of nMSD to general nonlinear settings using kernel-based manifold learning. We propose estimators for the kernelized nMSD and establish their consistency in high-dimensional regimes, enabling more flexible comparison of datasets via their spectral embeddings in reproducing kernel Hilbert spaces (RKHS) rather than in the ambient Euclidean space.
\end{itemize}
To support these methodological advances, we develop novel theoretical tools that are broadly applicable beyond the present context. In particular, our derivation of asymptotic normality results relies on key delocalization properties of the signal spectral representations, established by exploiting the latent manifold geometry. Moreover, our inference procedures are constructed through a careful analysis of a system of secular equations under heteroskedastic noise. Details are provided in Section~\ref{subsec:24}, and \SecCondCLT\ and \SecSecularEq \ in the Supplement.

The paper is organized as follows. We first conclude this section by collecting some notation used throughout the paper. In Section~\ref{model} we formally introduce our data-generating model and the nMSD, and characterize the connection between manifold geometry and the matrix spectral properties. Section~\ref{sec:estimator} presents our proposed methods and their asymptotic properties.  In Section~\ref{sec:sim} we conduct extensive simulations and analyze real-world biological data to assess the empirical performance of proposed methods. 
Our derivation of the main theoretical results and further simulation results are deferred to the supplement.

We use $\odot$ for elementwise Hadamard products: for equal-shaped arrays $A,B$, $(A\odot B)_{ij}=A_{ij}B_{ij}$, and $A^{\odot 2}$ denotes the elementwise square.
The symbol $\otimes$ denotes the Kronecker product; for vectors it reduces to the outer (tensor) product, and $u^{\otimes 4}:=u\otimes u\otimes u\otimes u$.
$\mathrm{Cum}_4(S)$ denotes the fourth-order cumulant tensor of a random vector $S$, with entries $\mathrm{Cum}(S_{i_1},S_{i_2},S_{i_3},S_{i_4})$.
For stochastic orders, $X_N=O_{\mathbb P}(a_N)$ means $X_N/a_N$ is bounded in probability, and $X_N=o_{\mathbb P}(a_N)$ means $X_N/a_N\to 0$ in probability.
For a sequence of events $(E_N)$, we say that $E_N$ holds \emph{with high probability} (w.h.p.) if there exist constants $C,c>0$, independent of $N$, such that \(\mathbb{P}(E_N) \ge 1 - C N^{-c}\)
for all sufficiently large $N$.
Matrix norms: $\|A\|_{\mathrm F}=\sqrt{\sum_{i,j}A_{ij}^2}$ (Frobenius) and $\|A\|_{\mathrm{op}}=\sup_{\|x\|_2=1}\|Ax\|_2$.
And $\|\cdot\|_{\psi_2}$ denotes the sub-Gaussian (Orlicz) norm:
$\ \|X\|_{\psi_2}:=\inf\{t>0:\ \E\exp(X^2/t^2)\le 2\}$.
For any bounded linear operator $T:\mathcal H\to\mathcal H$, write $[T]$ for its matrix in a fixed orthonormal basis $(e_1,\ldots,e_p)$.
Likewise, for $A:\mathcal H\to\R^N$, let $[A]\in\R^{N\times p}$ be the matrix of $A$ in the fixed basis of
$\mathcal H$ and the standard basis of $\R^N$.
Throughout the paper, $X_n \xrightarrow{\mathbb{P}} X$ denotes convergence in probability, and 
$X_n \xrightarrow{d} X$ denotes convergence in distribution.
We also denote by $\mathbb O(p)$ the orthogonal group in dimension $p$,
\(
\mathbb O(p) := \{ Q \in \mathbb{R}^{p\times p} : Q^\top Q = I_p \}.
\)

\section{Noisy Manifold Model, Spectral Dissimilarity, and Matrix Spectral Properties}
\label{model}

We formally introduce the manifold signal-plus-noise model and explain how the noiseless signals arise from the underlying manifold structure.
 We then provide a rigorous definition of the proposed dissimilarity between two datasets and describe the key spectral properties of the data matrices induced by their manifold geometry.

\subsection{Signal-plus-noise matrix model}

For the noisy observations, we assume that each data matrix $Y_k$ for $k=1,2$ follows the signal-plus-noise model defined in (\ref{eq:signal_plus_noise}). Without loss of generality, we assume the data is already centered so that \(\E[Y_k] =\E[S_{\samp,k}] = 0\).
We also make the following assumptions about the noise component $\Sigma_k^{1/2}X_k$ and the dimensionality of $Y_k$. For notational simplicity, throughout, we drop the dataset index $k$ whenever it is clear from the context that the statement applies to both datasets, and reinstate the subscript only when it is necessary.

\begin{assumption}[High-dimensional regime and block-heteroskedastic noise]\label{asmp:noise}
Throughout, we assume the following conditions hold.
\begin{enumerate}
\item \textbf{Proportional high-dimensional regime.}
For each dataset, we have $p,N \to \infty$ with $p/N \to \varphi \in (0,\infty)$.

\item \textbf{Entrywise scaling and sub-Gaussian tails.}
The entries $x_{i\mu}$ are independent across $(i,\mu)$, satisfy
$\E x_{i\mu}=0$ and $\Var(x_{i\mu})=1$, and are uniformly sub-Gaussian:
there exists a constant $K<\infty$ such that
$\|x_{i\mu}\|_{\psi_2}\le K$ for all $i,\mu$.

\item \textbf{Zero third cumulant.}
The noise is symmetric with vanishing third cumulant $\kappa_3 = 0$,
that is, $\E x_{i\mu}^3 = 0$ for all $i$ and $\mu$.

\item \textbf{Diagonal block-heteroskedastic noise.}
The noise covariance is diagonal and piecewise constant along a fixed
feature order: there exists a partition of $\{1,\ldots,p\}$ into
consecutive blocks $I_{1},\ldots,I_{K_p}$ such that
$
  \Sigma = \mathrm{diag}(\sigma_1,\ldots,\sigma_p),$ where $
  \sigma_i = \tau_j$ for all $i\in I_j$, with noise levels uniformly bounded: $  0 < \sigma_{\min} \le \tau_j \le \sigma_{\max} < \infty$ for all $j$.
Moreover, the number of blocks $K_p$ may depend on $p$, but satisfies $K_p = o(p),$
and the minimal block length diverges $
  m_p := \min_{1\le j\le K_p} |I_j| \;\to\; \infty$ as 
  $p\to\infty$.
\end{enumerate}
\end{assumption}

Conditions (i) (ii) and (iii) are standard in random matrix theory that facilitate precise asymptotic analysis of the matrix spectrum \citep{erdHos2017dynamical,bai2010spectral}. Condition (iv), whose motivation is discussed above, allows for heteroskedasticity in the noise components typically present in single-cell applications, which is more flexible than the commonly imposed homoskedastic or isotropic noise conditions  \citep{donoho2014minimax,DingMa2023,ding2025kernel} but more informative than the unstructured heteroskedastic setting \citep{Leeb2021MatrixDenoising} in literature. 


\subsection{Latent manifold model for noiseless samples}
\label{subsec:22}

For the latent signal component $S_{\samp,k}=[S_{1,k},\ldots,S_{N_k,k}]$ in the noisy observation model (\ref{eq:signal_plus_noise}), we assume the noiseless samples $\{S_{i,k}\}_{1\le i\le N_k}$ are independently drawn from a distribution supported on a low-dimensional smooth manifold that is embedded into the ambient high-dimensional space $\mathbb{R}^p$. A rigorous formulation of this assumption is provided below. 

\begin{assumption}[Manifold model for signals]\label{asmp:signal}
For $k\in\{1,2\}$, let $(\M_k,g_k)$ be a connected, compact, $m_k$-dimensional
Riemannian manifold and let $\iota_k:\M_k\to\R^{p}$ be an isometric embedding satisfying $\textup{dim}(\iota_k(\mathcal{M}_k))=n_k\le p$.
Let $\mathcal B(\M_k)$ denote the Borel $\sigma$–algebra on $\M_k$, and let
$\mathbb P_k$ be a probability measure on $(\M_k,\mathcal B(\M_k))$.
Assume $\mathbb P_k$ is absolutely continuous with respect to the Riemannian
volume measure $\operatorname{vol}_{g_k}$ on $\M_k$, with density bounded and
strictly positive, that is, $0<c_{k,\min}\ \le\ \frac{d\mathbb P_k}{d\operatorname{vol}_{g_k}}(x)\ \le\ c_{k,\max}<\infty$ for all $x\in \M_k$.
Let $\widetilde{\mathcal F}_k$ be the Borel $\sigma$–algebra on $\iota_k(\M_k)$, and define
the pushforward measure on $(\iota_k(\M_k),\widetilde{\mathcal F}_k)$ by $\rho_k\ :=\ \mathbb P_k\circ \iota_k^{-1}.$
For the noiseless samples, we assume  $S_{1,k},\ldots,S_{N_k,k}\stackrel{\text{i.i.d.}}{\sim} \rho_k$, which are supported on $\iota_k(\mathcal{M}_k)$.
\end{assumption}

Assumption~\ref{asmp:signal} states that the objects of interest are smooth signals; here, the signals are represented by the embedded manifold $\iota(\mathcal{M})$ and the pushforward probability measure $\rho$, which live on the ambient $\R^p
$ space, rather than the abstract manifold $\mathcal{M}$ and the probability measure $\mathbb{P}$.

Under Assumption~\ref{asmp:signal}, since  $n:=\dim(\iota(\M))\le p$,
by an orthogonal change of basis, 
there exists an orthonormal matrix $R\in \mathbb{O}(p)$ such that
$(R\circ\iota)(\M)\subset\R^{n}\times\{0\}^{\,p-n}$. Define the coordinate functions $
  f_j(x)\ :=\ e_j^\top (R\circ\iota)(x)$, for $ j=1,\ldots,n$ and $x\in\M.$
Then for every $x\in\M$, it admits  a sparse Euclidean representation 
\begin{equation}\label{sparse}
     (R\circ\iota)(x)\;=\;\bigl(f_1(x),\ldots,f_n(x),0,\ldots,0\bigr)\in\R^p.
\end{equation}
With a slight abuse of notation, we write
$RS_i=(f_1(S_i),...,f_n(S_i),0,...,0)\in\R^p.$


Under Assumption~\ref{asmp:signal}, the reduction in (\ref{sparse}) implies that the noiseless samples are effectively supported on an $n$-dimensional subspace. This insight enables the inference and comparison of latent signal structures across datasets in a basis-independent manner. 
We next introduce the population and empirical covariance matrices associated with the latent manifold model.
\color{black}

\begin{definition}[Population and empirical signal covariance matrices]\label{def:cov-mat}
Let $S_1,\ldots,S_N\in\R^p$ be centered noiseless samples drawn i.i.d.\ from $\rho$ supported on $\iota(\M)$ as in Assumption~\ref{asmp:signal}.
Denote $S_{\samp}=[\,S_1,\ldots,S_N\,]\in\R^{p\times N}$.
Then the population and empirical signal covariance matrices are defined as $
M := M(\rho) = \mathbb{E}[S_i S_i^\top] = \int_{\iota(\mathcal{M})} S_i S_i^\top d\rho$ and $
\widehat M := \frac{1}{N} S_{\text{samp}} S_{\text{samp}}^\top.$
\end{definition}

The population covariance matrix $M$ is uniquely determined by $\rho:=\mathbb P\circ\iota^{-1}$, the pushforward probability measure supported on $\iota(\mathcal{M})\subset\R^p$ associated with
$(\M,\mathbb P,\iota)$. Because $\iota(\M)$ is compact and the samples are centered,
$M=\E[S_iS_i^\top]$ exists, has $\operatorname{rank}(M)\le n$, and as a linear operator on $\R^p$, vanishes on $V^\perp$ where
$V:=\operatorname{span}(\iota(\M))$. For any unit vector $v\in\R^p$, we have
$
\Var\big(\langle v,S_i\rangle\big):= \int_{\iota(\mathcal{M})} \langle v,S_i\rangle^2 d\rho=v^\top M v,
$
which is the variance of $\rho$ along direction $v$. Hence the top eigenvalue of the population covariance matrix $M$, given by $
\mathcal D_1\;=\;\sup_{\|v\|_2=1} v^\top M v$,
characterizes the principal variance of the distribution $\rho$  on the
embedded manifold $\iota(\mathcal{M})$ along the axis of maximal spread \citep{huckemann2006principal,dai2018principal}. Recursively, the $j$-th eigenvalue
$\mathcal D_j$ arises from the same maximization subject to additional constraint $v\perp\operatorname{span}\{\eta_1,\ldots,\eta_{j-1}\}$; directions orthogonal to $V$ carry zero variance. 
Thus, $\{\mathcal D_j\}_{1\le j\le n}$ capture the overall magnitude of the latent signal structure along its intrinsic principal axes. They are orthogonally invariant to the choice of basis in $\iota(\mathcal{M})$ and remain intrinsic to the ambient Euclidean space $\mathbb{R}^p$.
Following the convention in \cite{DingMa2023,rmt1,LinPanZhaoZhou2024,smai}, hereafter, we refer to these population-level eigenvalues  $\{\mathcal D_j\}_{1\le j\le n}$ as the true signal strength parameters. 


\subsection{Manifold spectral dissimilarity measure}

We next introduce a dissimilarity measure between two datasets based on the eigenvalues $\{\mathcal D_j\}_{1\le j\le n}$ of the population covariance matrix.

\begin{definition}[Normalized population principal variances]\label{def:Pi}
Let $\rho$ be the pushforward measure defined in Assumption~\ref{asmp:signal}, and  $\mathcal{D}_{1}\ge\mathcal{D}_{2}\ge\cdots$ be the ordered eigenvalues of $M$ associated with $\rho$. If we denote $\mathcal{D}^{(r)} \;:=\; \bigl(\mathcal{D}_{1},\dots,\mathcal{D}_{r}\bigr)^{\!\top}$ containing the leading $r$ eigenvalues with $r\le p$, then
the normalized top-$r$ principal variances  of $\rho$ are defined as $
  \Pi_r(\rho)
  := \frac{\mathcal{D}^{(r)}}{\mathbf 1^\top \mathcal{D}^{(r)}}
  \;\in\; \Delta^{\,r-1}
  :=\bigl\{x\in\R^r : x_j\ge0,\ \sum_{j=1}^r x_j=1\bigr\}.$
\end{definition}

\begin{definition}[Manifold spectral dissimilarity]\label{def:dist}
Let $\rho_k$ be a probability measure defined by Assumption~\ref{asmp:signal}, associated with some $(\mathcal{M}_k,\mathbb{P}_k,\iota_k)$ for $k=1,2$.  Let $\Pi_r(\rho_k)$ be the normalized top-$r$ principal variances of $\rho_k$.
The normalized manifold spectral dissimilarity (nMSD) between \(\rho_1\) and \(\rho_2\) is defined as $
  d_r(\rho_1,\rho_2)
  := \bigl\|\Pi_r(\rho_1)-\Pi_r(\rho_2)\bigr\|_2$.
\end{definition}

Our dissimilarity measure is defined on the pushforward distributions  $\rho_1$ and $\rho_2$, defined on their support manifolds  $\iota_1(\mathcal{M}_1)$ and $\iota_2(\mathcal{M}_2)$. Two datasets may reside on the same embedded manifold, i.e., $\iota_1(\mathcal{M}_1) = \iota_2(\mathcal{M}_2)$, yet differ in their sampling densities, thereby inducing a nonzero nMSD. Conversely, two distinct distributions may share the same normalized principal variances and hence have zero dissimilarity under $d_r$. This is an intentional consequence of the resolution at which $d_r$ compares datasets: our goal is not to distinguish arbitrary distributional differences, but to identify datasets with similar dominant structural variation.
In Section~S2 of the Supplement, we provide two examples demonstrating these basic properties of $d_r(\rho_1,\rho_2)$.

Compared with existing approaches that quantify dataset similarity using distance measures on entire distributions, such as KL divergence or MMD, or on covariance matrices, such as the Frobenius distance, the proposed measure is invariant to scale changes and orthogonal transformations and is insensitive to perturbations that do not affect the leading principal variances. This yields a coarser but more stable notion of similarity that is well suited for comparing noisy and high-dimensional datasets. The resulting measure is also computationally simple and readily interpretable in terms of the shape and dimensionality of the latent structure.

\begin{remark}[Choosing the working rank $r$]
If the two embedded manifolds share the same effective Euclidean
dimension $n$, one can naturally set $r=n$. More generally, if their
effective dimensions are $n_1\ne n_2$, it is reasonable to take
$r=\max\{n_1,n_2\}$.
In practice, when $n$ or $(n_1,n_2)$ are unknown, one may consistently estimate the signal rank using existing methods like \cite{jiang2022universal} and \cite{donoho2023screenot}. As such, hereafter we treat $r=n_1=n_2$ as given; in our numerical implementations, we set \(r=\max \{\hat n_1,\hat n_2\}\) where the data-specific rank $n_k$ is estimated by \cite{donoho2023screenot}.
\end{remark}

\subsection{Spectral properties of the signal matrices}
\label{subsec:24}
When the noiseless samples $S_1,\dots,S_N$ are $i.i.d.$ sampled from the above latent manifold model, the corresponding signal matrix $S_\samp$ and the underlying population covariance matrix $M$ exhibit several important spectral properties. These properties, summarized below, play a key role in  our method development.

\begin{proposition}[Signal eigenvalues]\label{prop:pca-align}
Under Assumption~\ref{asmp:signal}, let $S\in\R^p$ be centered with $\E[S]=0$ and let $M:=\E[SS^\top]$.
With the coordinate functions $f_1,\ldots,f_n$ already defined by
$(R\circ\iota)(x)=(f_1(x),\ldots,f_n(x),0,\ldots,0)$, we may, without loss of generality,
apply an additional orthogonal change of coordinates within the first $n$ axes so that
for the same notation $f_1,\ldots,f_n$,  we have $
\E[f_i(S)]=0$,
$\E\big[f_i(S)f_j(S)\big]=0~(i\neq j)$, and $
\mathcal D_i=\E\big[f_i(S)^2\big]$,
with $\mathcal D_1\ge\cdots\ge\mathcal D_n\ge 0$.
In particular, $\{\mathcal D_i\}_{i=1}^n$ are exactly the nonzero eigenvalues of $M$.
\end{proposition}

Let $S_\samp=U_\samp D_\samp V_\samp^\top$ be the SVD of $S_\samp$. For any $R\in\mathbb O(p)$, a change of ambient basis gives
$
S_\samp^{(\mathrm{rot})}=R\,S_\samp=(R\,U_\samp)\,D_\samp\,V_\samp^\top.
$
Hence the singular values $D_\samp$ and the right singular vectors $V_\samp$ are invariant
under ambient orthogonal rotations; only the left singular vectors rotate.

\begin{proposition}[Low-rankness of the sample signal matrix]
\label{pro:lowrank}
Under Assumption~\ref{asmp:signal}, let
$S_\samp=[S_1,\dots,S_N]\in\R^{p\times N}$.
Then $\rank(S_\samp)\le n$, and
$
\Pr\!\left(\rank(S_\samp)=n\right)\to 1
$
as $N\to\infty$.
\end{proposition}

Proposition~\ref{pro:lowrank} is intuitive since the noiseless samples are generated from a  manifold embedded in a low-dimensional subspace. Even though the feature dimension $p$ increases with the sample size $N$, the true signal dimension $n$ and the rank of $S_\samp$ remain small.

\begin{proposition}[Matrix spectral concentration, stability, and delocalization]\label{prop:matrix-spectral}
Under Assumption~\ref{asmp:signal}, let $S_1,\ldots,S_N\in\R^p$ be centered noiseless samples drawn i.i.d.\ from $\rho$ supported on $\iota(\M)$.
Let $(\mathcal D_k,\eta_k)_{k=1}^n$ be the nonzero eigenpairs of $M$ with
$\mathcal D_1\ge\cdots\ge\mathcal D_n>0$ and assume the uniform eigengap $\inf_{1\le k\le n} \min\{\mathcal{D}_{k-1}-\mathcal{D}_k, \mathcal{D}_k-\mathcal{D}_{k+1}\} \ge c > 0, 
\text{where } \mathcal{D}_0 := +\infty, \ \mathcal{D}_{n+1} := 0.$
Since $\iota(\M)$ is compact, set $
  \kappa^2 := \sup_{x\in\iota(\M)} \|x\|_2^2 < \infty.$
Take the thin SVD $S_{\samp}=U_{\samp}D_{\samp}V_{\samp}^\top$ with
$U_{\samp}\in\R^{p\times n}$, $V_{\samp}\in\R^{N\times n}$ containing orthonormal columns and
$D_{\samp}=\mathrm{diag}(d_1,\ldots,d_n)$. Then for any $\delta\in(0,1)$, the following hold.

\begin{enumerate}
\item \textbf{Operator-norm concentration.} For any small constant $\delta>0$, we have $
\Pr\!\Big(\bigl\|\widehat M-M\bigr\|_{\op}\ \ge\ \frac{4\,\kappa^2\,\log(2/\delta)}{\sqrt N}\Big)\ \le\delta.$
\item \textbf{Manifold interpretation of singular values.}
For any fixed $\varepsilon>0$, there exist constants $C,c>0$ such that
$
\Pr\!\left(
\max_{1\le i\le n}
\left|
\frac{d_i}{\sqrt N}-\sqrt{\mathcal D_i}
\right|
>\varepsilon
\right)
\le C N^{-c}.
$
Moreover, the nonzero eigenvalues of $\widehat M$ are
$\{d_i^2/N\}_{i=1}^n$.
\item \textbf{Manifold interpretation of right singular vectors  (sample embeddings).}
For any fixed $i\in[n]$ and $\varepsilon>0$, there exist $C,c>0$ such that $
\Pr\!\left(\max_{1\le j\le N}\Big|\,V_{\samp}(j,i)\;-\;\frac{f_i(S_j)}{\sqrt{N\,\mathcal D_i}}\,\Big|>\varepsilon\right)
\ \le\ C\,N^{-c}.$
\item \textbf{Delocalization of right singular vectors (sample embeddings).}
There exist $C,c>0$ such that, with probability at least $1-N^{-c}$,
 we have $
\forall\,1\le i\le n:\quad
\max_{1\le j\le N}\,|V_{\samp}(j,i)|\ \le\ \frac{C\sqrt{\log N}}{\sqrt N}\,.$
\end{enumerate}
\end{proposition}

Proposition~\ref{prop:matrix-spectral} characterizes the asymptotic relationship between the singular values and vectors of the signal sample matrix $S_\samp$ and the underlying population $\rho$ on the manifold. In particular, part (i) establishes the asymptotic equivalence between the population and the empirical covariance matrices; part (ii) connects the  singular values of the signal matrix $S_\samp$ to the eigenvalues of the population covariance matrix $M$  on the manifold; part (iii) relates the right singular vectors, or the sample spectral embeddings, with the latent Euclidean representation of the noiseless samples on the manifold, providing a geometric interpretation of the singular vectors; and part (iv) describes a key delocalization property of the right singular vectors, which arises naturally from the smoothness of the latent manifold.
Importantly, the delocalization property in part~(iv) does not hold in general; rather, it follows directly from the manifold assumption and plays a pivotal role in deriving the asymptotic distributions of the proposed estimators.
\color{black}

\section{Main results}

\label{sec:estimator}

Next we introduce a unified framework for inferring the nMSD between two noisy high-dimensional datasets. The proposed method proceeds in two stages: first, we construct a consistent estimator of the unknown noise covariance matrix (Section~\ref{subsec:31}); then, after accounting for the noise heteroskedasticity, we develop  an estimator and associated inference procedures for the true nMSD (Sections \ref{subsec:simnMSD} and \ref{subsec:33}). To streamline our presentation, we focus on the statistical procedures and their properties and defer the technical details to \SecTech \ in the Supplement.

\subsection{Noise variance estimation}
\label{subsec:31}

Leveraging the facts that the signal matrix $S_\samp$ is low-rank (Proposition~\ref{pro:lowrank}) with delocalized singular vectors (Proposition~\ref{prop:matrix-spectral}), and that the noise covariance admits a block structure (Assumption \ref{asmp:noise}), we propose computationally efficient noise estimators for each dataset. Specifically, our approach performs a low-rank fitting followed by one-dimensional Potts segmentation, which is solved exactly using dynamic programming \citep{Janssens2007CME,Weinmann2014L1Potts}. The details of the estimation procedure is summarized in Algorithm~\ref{alg:noise}.
 Specifically, a fitted rank-$r$ signal matrix subtracts the signal energy, leaving a residual diagonal that concentrates around the true noise diagonal as the sample size increases. Pooling these residuals along a fixed feature ordering with a one-dimensional Potts penalty then further denoises within blocks and detects change points: the segmentation introduces a small bias but yields substantial variance reduction, and when jumps are well separated, it accurately recovers both block means and boundaries. The following theorem provides a consistency guarantee of the proposed variance estimator and  identifies the appropriate choice of the tuning parameter $\beta$.

\begin{algorithm2e}[t]
\spacingset{1} 
  \small
\caption{Noise variance estimation}
 \label{alg:noise}
\KwIn{data matrix $Y\in\mathbb{R}^{p\times N}$; working rank $r=n\ge 1$; a fixed 1-D feature order on $\{1,\ldots,p\}$; segmentation penalty $\beta>0$.}
\KwOut{$\widehat\Sigma=\mathrm{diag}(\widehat\sigma^{\mathrm{pwc}})$.}

\textbf{1. Low-rank fitting.}
Compute $Q=\frac{1}{N}YY^\top$. Obtain the best rank-$n$ approximation
$\widehat C_n=\widehat U_n\widehat D_n\widehat U_n^\top$ to $Q$ via truncated SVD, where $\widehat U_n\in\R^{p\times n}$ and $\widehat D_n\in\R^{n\times n}$ contain leading $n$ eigenvectors and eigenvalues of $Q$.
Set the raw variance estimators $\widehat\sigma$ as the diagonal elements of the spectral truncation residual 
$Q-\widehat C_n$.

\textbf{2. Piecewise smoothing via Potts segmentation.}
Regularize along the feature order by solving the Potts segmentation problem
\begin{equation}\label{potts}
\widehat\sigma^{\mathrm{pwc}}
\in\arg\min_{\sigma\in\mathbb{R}^p}
\Big\{\sum_{i=1}^p\big(\widehat\sigma_i-\sigma_i\big)^2
\;
+\;\beta\,\#\{1\le i<p:\sigma_{i+1}\neq\sigma_i\}\Big\},
\end{equation}
where $\beta$ is a tuning parameter. This Potts segmentation problem is computed exactly by dynamic programming.

\textbf{3. Return the noise covariance estimator.}
Set $\widehat\Sigma\gets \mathrm{diag}(\widehat\sigma^{\mathrm{pwc}})$.

\end{algorithm2e}


\begin{theorem}[Consistency of noise variance estimators]\label{thm:noise-consistency}
Suppose Assumption~\ref{asmp:noise} and Assumption~\ref{asmp:signal} hold. 
Let $\widehat\sigma^{\mathrm{pwc}}$ be defined as in (\ref{potts}) in Algorithm \ref{alg:noise} with $\beta=c\,\frac{\log p}{N}$. Then with probability $1-o(1)$, we have $\frac{1}{p}\,\big\|\widehat\sigma^{\mathrm{pwc}}-\sigma\big\|_2^2\ \lesssim\ \frac{\log p}{N}.$
\end{theorem}

The rate of convergence ${(\log p)}/{N}$ essentially follows from the concentration of the spectral truncation residual diagonal $\widehat\sigma^{\textup{raw}}$ and a Potts oracle inequality with $\beta=c(\log p)/N$ \citep{Boysen_2009}; proofs and precise constants are deferred to the Supplement. In practice, we found the procedure to be robust with respect to the choice of \(c\) across a broad range; unless otherwise stated, we use the default value \(c=10\) in our own implementation.

Empirically, we find that the Potts segmentation step in  Algorithm \ref{alg:noise} significantly improves the initial variance estimates $\widehat\sigma^{\textup{raw}}$. Detailed numerical results from the simulation studies are provided in the Supplement (Section S4.5). Notably, 
Figure S2 of the Supplement shows,  with a properly chosen penalty parameter $\beta$, the Potts-smoothed estimators achieve a mean squared error nearly an order of magnitude smaller than that of the raw residual estimators.

\subsection{Dataset similarity inference using nMSD}
\label{subsec:simnMSD}
Next we develop an estimator of the nMSD between two noisy datasets. We begin by introducing some necessary notation. 
We define the total population covariance matrix as $\tilde\Sigma=\frac{1}{N}\E[YY^\top]$, where the expectation is with respect to both the noise distribution and the signal distribution $\rho$ from the latent manifold model, and denote $(\xi_k,\popsi_k)$ as its \(k\)-th eigenpair. Conditional on the noiseless samples $S_\samp$, we define  $\SigmaS=\frac{1}{N}\E[YY^\top| S_\samp]$ as the conditional population covariance matrix, and denote 
$(\pssig_k,\psi_k)$ as its \(k\)-th eigenpair. Note that $\tilde\Sigma$ and $\SigmaS$ are different from the population and empirical signal covariance matrices $M$ and $\widehat M$  as the latter ones are defined for the noiseless samples $S_1, ..., S_N$ only.
For a spike location $\pssig\notin\{\sigma_1,\dots,\sigma_p\}$, define functions \begin{align} 
A(\pssig)\;:=\;(\Sigma-\pssig I_p)^{-1}
=\diag\!\Big(\frac{1}{\sigma_1-\pssig},\dots,\frac{1}{\sigma_p-\pssig}\Big),\label{A}\\
g(\pssig)\;:=\;\frac1p\,\Tr A(\pssig),\quad s_2(\pssig)\;:=\;\Tr A(\pssig)^2.\label{g,s2} \end{align}
Fix any matrix \(U\in\mathbb{R}^{p\times r}\) with orthonormal columns, i.e. \(U^\top U = I_r\), and define $
  G_U(\psi) := U^\top A(\psi) U \in \mathbb{R}^{r\times r}.$

Let $\{\lambda_k\}$ be the  eigenvalues of the sample covariance matrix $Q=\frac{1}{N}YY^\top$, and let $\{\theta_k\}$ denote the theoretical limits of $\{\lambda_k\}$ conditional on the noiseless samples $S_\samp$ as $(N,p)\to\infty$; in particular, as will be discussed in \LemBdeloc \ and \LemCondCLT \ in the Supplement, $\theta_k$ is uniquely determined by the conditional population covariance eigenvalues $\{\pssig_k\}$, the noise covariance $\Sigma$ and the aspect ratio limit $\varphi$. Assuming $\Sigma$ and $\varphi$ are known, with a slight abuse of notation we can write $\theta_k:=\theta(\pssig_k)$ so that $\lambda_k\approx \theta(\pssig_k)$ for large $N$ and $p$. When $\theta(\cdot)$ is known or estimable from the data, we can thus construct an estimator of the conditional population  covariance eigenvalue $\pssig$ via the inverse map \(\widehat{\pssig}_k:=\theta^{-1}(\lambda_k)\) for \(k=1,\dots,r\), and subsequently obtain its standard error estimates by deriving a central limit theorem (CLT) of this estimator. 

From this section onward, we also impose a basic signal-strength condition on the
spikes of the conditional population covariance matrix $\SigmaS$:

\begin{assumption}[SNR]\label{asmp:signal-strength}
The leading $r$ eigenvalues $\{\pssig_k\}_{k\le r}$ of $\SigmaS$, with $r\le n$, are
supercritical and mutually separated;
particularly, each $\pssig_k$ lies outside the bulk spectrum of the noise
covariance and is separated from the remaining eigenvalues by a fixed gap.
\end{assumption}

\noindent
A precise definition of the ``supercritical" and  ``separation" conditions, including explicit relationships between the spike eigenvalues and the noise bulk edge, is given in Section~S4 of the Supplementary Material.
Assumption~\ref{asmp:signal-strength} ensures that the leading $r$ spikes of the
conditional population covariance $\SigmaS$ generate statistically detectable outlier
eigenvalues in the sample covariance matrix. In particular, each $\pssig_k$ remains
bounded away from the noise bulk and from other spikes as $(N,p)\to\infty$, so that
both the outlier map $\theta(\cdot)$ and its inverse $\theta^{-1}(\cdot)$ are well
defined and smooth in a neighborhood of these points. This condition is standard in
spiked covariance models \citep{erdHos2017dynamical,bai2010spectral} and is necessary for obtaining consistent estimation and
asymptotic normality of the eigenvalue estimators.

For given datasets $Y_1$ and $Y_2$, to estimate their latent nMSD $d_r(\rho_1,\rho_2)$, inspired by Proposition~\ref{prop:matrix-spectral}, we  first construct estimators of the leading empirical signal covariance eigenvalues $\{d_k^2/N\}$, and use them to approximate the population covariance eigenvalues $\{\mathcal{D}_k\}$; these estimates are then substituted into the definition of the nMSD. The procedure proceeds as follows. First, we  leverage the estimated noise covariance $\widehat\Sigma$ from Algorithm \ref{alg:noise} and the top sample covariance eigenvalues \(\lambda_j\), inverting the outlier eigenvalue map \(\theta\), to obtain estimators \(\widehat\pssig_j\) of the leading conditional population covariance eigenvalues; the concrete expressions of $\theta(\cdot)$ and $\theta'(\cdot)$ are obtained using RMT and are given by Equation (S1) in the Supplement. Second, we convert these conditional population covariance eigenvalues into estimators of the signal covariance eigenvalues by solving for $\{d_1,...,d_r\}$ from the following system of secular equations obtained from RMT that connect the two systems in the asymptotic limit:
\begin{equation}
\label{eq:rank-r-secular}
\det\!\big(I_r+D_r^2\,G_{U_r}(\widehat\pssig_j)\big)=0,  \quad j=1,2,\dots,r,
\end{equation}
where $D_r = \text{diag}(d_1,d_2,...,d_r)$, $U_r\in\R^{p\times r}$ contains the top $r$ left singular vectors of $S_\samp$. Once the solution to (\ref{eq:rank-r-secular}) is obtained, denoted as $\{\widehat d_j\}$, we calculate the normalized signal strength vector
\(
\widehat\Pi_r
={\big(\widehat d_1^{\,2},\dots,\widehat d_r^{\,2}\big)}/
       {\sum_{j=1}^r \widehat d_j^{\,2}}
\in\Delta^{r-1}
\) and compute the dissimilarity estimator. Importantly, the two-step processing by the inversion map $\theta^{-1}(\cdot)$ and solving the secular equations (\ref{eq:rank-r-secular}) accurately adjust for the intrinsic bias in the sample eigenvalues $\{\lambda_k\}$  caused by  noise heteroskedasticity and  high dimensionality reflected by the non-vanishing $\varphi$.  
The derivation of the secular equations (\ref{eq:rank-r-secular})  are provided in \SecCondCLT\ and \SecSecularEq\ in the Supplement. 

To enable efficient computation of the final nMSD estimator, we propose an algorithm that implements the estimation procedure described above (Algorithm~\ref{alg:pi}). 
Specifically, we consider a noninformative prior on $U_r$ and show that under 
this assumption, $G_{U_r}(\xi)$ can be well approximated by $g(\xi) I_r$ (recall $g(\cdot)$ defined in (\ref{g,s2})), 
which is independent of $U_r$ (\LemIsotropy). 
This observation substantially simplifies the secular system (\PropDecouple) and enables a practical and 
fast computation of the estimator.

\begin{algorithm2e}
\caption{nMSD estimation procedure}\label{alg:pi}
\spacingset{1} 
  \small
\KwIn{data matrix $Y_k\in\mathbb{R}^{p\times N_k}, k=1,2$; working rank $r\ge 1$; estimated diagonal noise covariance $\widehat\Sigma_k=\mathrm{diag}(\widehat \sigma_{1,k},\dots,\widehat\sigma_{p,k}), k=1,2$.}
\KwOut{$\widehat\Pi_{r,1},\widehat\Pi_{r,2}\in\Delta^{r-1}$ and the estimated nMSD $\|\widehat\Pi_{r,1}-\widehat\Pi_{r,2}\|_2$}
\BlankLine
\textbf{1. Sample spikes:}
Form $Q_k=\frac{1}{N_k}Y_kY_k^{\top}$ and let $\lambda_{1,k}\ge\cdots\ge\lambda_{r,k}$ be its top $r$ eigenvalues.
\BlankLine
\textbf{2. Invert the outlier map $\theta$:}
For $j=1,\dots,r$ and $k=1,2$, solve the scalar equation $\theta\!\left(s;\widehat\Sigma_k,\varphi_k=p/N_k\right)=\lambda_{j,k}$ for $s>\max_i\widehat\sigma_{i,k}$, and set $\widehat{\xi}^{\mathrm{S}}_{j,k}\gets s$. 
\BlankLine
\textbf{3. Calculate signal strengths:}
 For each $j=1,\dots,r$, set $\widehat d_{j,k}^{\,2}=-\frac{1}{\widehat g\!\big(\widehat{\xi}^{\mathrm{S}}_{j,k}\big)}\,$
where 
$
\widehat g(s)=\frac{1}{p}\sum_{i=1}^p\frac{1}{\widehat\sigma_{i,k}-s},$ for all $s>\max_i \widehat \sigma_{i,k}$.

\textbf{4. Calculate principal variance profiles:}
For $k=1,2,$ set $\displaystyle \widehat\Pi_{r,k}=\frac{(\widehat d_{1,k}^{\,2},\dots,\widehat d_{r,k}^{\,2})}{\sum_{i=1}^r \widehat d_{i,k}^{\,2}}$.

\textbf{5. Calculate dissimilarity estimator:}
\KwReturn{$\|\widehat\Pi_{r,1}-\widehat\Pi_{r,2}\|_2$}
\end{algorithm2e}

The next theorem establishes the consistency of the empirical signal covariance eigenvalue, normalized principal variances, and nMSD estimation in
Algorithm~\ref{alg:pi}.

\begin{theorem}[Consistency of dataset similarity estimation]
\label{thm:signal-consistency-unified-short}
Suppose Assumption~\ref{asmp:noise}, Assumption~\ref{asmp:signal} and Assumption~\ref{asmp:signal-strength} hold, and fix $r\ge1$. Then $\{\widehat{\pssig}_j\}_{j=1}^r$ are consistent estimators of the
conditional population spikes $\{\pssig_j\}_{j=1}^r$. 
Consider the empirical secular system obtained by replacing $(\pssig_j,\Sigma)$ with
$(\widehat{\pssig}_j,\widehat\Sigma)$ in (\ref{eq:rank-r-secular}), where $\widehat\Sigma$ is any consistent estimator of the diagonal noise $\Sigma$.
{If this system admits a solution $\widehat D_r^{\,2}=\diag(\widehat d_1^{\,2},\ldots,\widehat d_r^{\,2})$,}
then
\begin{equation}\label{consis}
\widehat d_j^{\,2} \xrightarrow{\ \mathbb{P}\ } d_j^{\,2}, \quad \text{for all $j=1,\ldots,r,$ and }\quad 
\widehat\Pi_r := \frac{\widehat d^{\,2}}{\mathbf{1}^\top \widehat d^{\,2}} \xrightarrow{\ \mathbb{P}\ } \Pi_r := \frac{d^{\,2}}{\mathbf{1}^\top d^{\,2}}.
\end{equation}
Consequently, for two independent datasets generated under the same framework, we have
$
\bigl\|\widehat\Pi_{r,1}-\widehat\Pi_{r,2}\bigr\|_2\ \xrightarrow{\ \mathbb{P}\ }\ d_r(\rho_1,\rho_2).
$
Furthermore, assume that 
$U_r$ is obtained as the first $r$ columns of a random orthogonal matrix $Q$ drawn from the Haar distribution on  $\mathbb{O}(p)$. Then the secular system admits the following closed-form solution:
\begin{equation}
    \widehat d_j^{\,2}\;=\;-\frac{1}{\widehat g(\widehat{\pssig}_j)}\,,\quad
\widehat g(s)=\frac{1}{p}\sum_{i=1}^p\frac{1}{\widehat\sigma_i-s}\,,\quad
\ j=1,\ldots,r,
\end{equation}
which also satisfies $\widehat d_j^{\,2}\xrightarrow{\mathbb{P}} d_j^{\,2}$ and hence
$\widehat\Pi_r\xrightarrow{\mathbb{P}}\Pi_r$ as above.
\end{theorem}

Our framework is not restricted to the case where $U_r$ follows a noninformative prior. In cases where  $U_r$ has a different prior distribution encoding additional structural knowledge, we may apply a similar argument (as in \LemIsotropy\ in the Supplement) to approximate $G_{U_r}(\xi)$ by some deterministic matrix independent of $U_r$. When such a theoretical reduction is infeasible, one can always  apply  Monte Carlo methods to obtain numerical solutions for $\{d_j\}$. Equation (\ref{consis}) of Theorem~\ref{thm:signal-consistency-unified-short} ensures that, in all such cases, the proposed estimators remain consistent.

In addition to the estimation consistency, the proposed nMSD estimator, as well as the components of $\widehat\Pi_r$, also admit asymptotic normality, which are derived in detail in \SecDeltaMethod \ in the Supplement.
 Building on this theoretical insight, we can construct confidence intervals for the components of $\Pi_r$ and for the underlying nMSD between two datasets, providing uncertainty quantification of these estimators.
 
Let $\widehat\Pi_{r,k}\in\Delta^{r-1}, k=1,2,$ be the principal variance
profile estimators from Algorithm~\ref{alg:pi}. Let
$\widetilde V_\Pi$ denote any consistent estimator of
$\Cov(\widehat\Pi_r)$ satisfying
$\|N\widetilde V_\Pi-\Sigma_\Pi\|_{\mathrm{op}}
\xrightarrow{\mathbb P}0$,
where $\Sigma_\Pi$ is the limiting covariance defined in
\EqPiCLT\ and \EqSigmaPiMult\ in the Supplement. 
For $\alpha\in(0,1)$, let $z_{1-\alpha/2}$ denote the
$(1-\alpha/2)$-quantile of the standard normal distribution. We can define
the componentwise $(1-\alpha)$ confidence intervals
$\mathrm{CI}_{t}^{(\Pi)}
:=
\bigl[
(\widehat\Pi_r)_t
\pm
z_{1-\alpha/2}{[\widetilde V_\Pi]^{1/2}_{tt}}
\bigr]$,
for $t=1,\dots,r$. For each fixed $t$ such that
$[\Sigma_\Pi]_{tt}>0$, these intervals satisfy
$\Pr\!\bigl((\Pi_r)_t\in\mathrm{CI}_{t}^{(\Pi)}\bigr)
\to 1-\alpha$
as $(N,p)\to\infty$. 
Similarly, for two independent datasets producing
$\widehat\Pi_{r,1}$ and $\widehat\Pi_{r,2}$, let
$\widetilde V_\Pi^{(1)}$ and $\widetilde V_\Pi^{(2)}$ denote covariance
estimators satisfying
$\|N_i\widetilde V_\Pi^{(i)}-\Sigma_{\Pi,i}\|_{\mathrm{op}}
\xrightarrow{\mathbb P}0$
for $i=1,2$. Define
$\Delta\widehat\Pi
:=
\widehat\Pi_{r,1}-\widehat\Pi_{r,2}$
and
$\widetilde V_\Delta
:=
\widetilde V_\Pi^{(1)}+\widetilde V_\Pi^{(2)}$.
The component-wise $(1-\alpha)$ confidence intervals for
$\Delta\Pi:=\Pi_r(\rho_1)-\Pi_r(\rho_2)$ are given by
$\mathrm{CI}_t^{(\Delta)}
:=
\bigl[
(\Delta\widehat\Pi)_t
\pm
z_{1-\alpha/2}{[\widetilde V_\Delta]^{1/2}_{tt}}
\bigr]$, and
the
$(1-\alpha)$ confidence interval for $d_r(\rho_1,\rho_2)$ is given by
$\mathrm{CI}_{\textup{nMSD}}
:=
\biggl[
\widehat d_r(\rho_1,\rho_2)
\pm
z_{1-\alpha/2}
\sqrt{
\frac{
\Delta\widehat\Pi^\top
\widetilde V_\Delta
\Delta\widehat\Pi
}{
\widehat d_r(\rho_1,\rho_2)^2
}
}
\biggr]$.

\begin{proposition}
\label{pro:CInMSD}
Under Assumptions~\ref{asmp:noise}, \ref{asmp:signal}, and
\ref{asmp:signal-strength}, for each $t=1,...,r,$ such that $[\Sigma_\Pi^{(1)}+\Sigma_\Pi^{(2)}]_{tt}>0$, and
$\widetilde V_\Pi^{(1)}$ and $\widetilde V_\Pi^{(2)}$ satisfy the consistency conditions stated above, we have   $\Pr\!\bigl(
(\Delta \Pi)_t
\in
\mathrm{CI}_{t}^{(\Delta)}
\bigr)
\to 1-\alpha$
as $(N,p)\to\infty$. Suppose that
$d_r(\rho_1,\rho_2)$ is bounded away from zero, and the limiting variance
of the nMSD estimator is positive. Then
$\Pr\!\bigl(
d_r(\rho_1,\rho_2)
\in
\mathrm{CI}_{\textup{nMSD}}
\bigr)
\to 1-\alpha$
as $(N,p)\to\infty$.
\end{proposition}

Convenient plug–in estimators $\widehat V_\Pi^{(1)}$ and $\widehat V_\Pi^{(2)}$, obtained by
replacing population quantities in the variance formulas with their empirical counterparts, is
given in  \EqSigmaplugin\; see Section S1.5 of the Supplement for details.

\subsection{Dataset alignability inference}
\label{subsec:33}
To determine the alignability between two datasets, under our framework, we test the null hypothesis  $
H_0:\ \Pi_r(\rho_1)=\Pi_r(\rho_2).$
From the estimation pipeline above, recall the normalized vectors
$\widehat\Pi_{r,k}\in\Delta^{r-1}$ and their respective limiting covariances
$\Sigma_{\Pi,k}$ for $k=1,2$.
Let $\Delta_\Pi:=\widehat\Pi_{r,1}-\widehat\Pi_{r,2}$, $
N_{\mathrm{eff}}
:=
\frac{N_1N_2}{N_1+N_2},
\alpha
:=
\frac{N_2}{N_1+N_2},
$
and
$
\Sigma_\Delta
:=
\alpha\Sigma_{\Pi,1}
+
(1-\alpha)\Sigma_{\Pi,2}.
$ Define the test statistic
$
T_\Pi
:=N_{\mathrm{eff}}\,\Delta_\Pi^\top\bigl(\Sigma_\Delta\bigr)^{+}\Delta_\Pi$.
Here we use the Moore--Penrose inverse
because by definition $\Sigma_\Delta\mathbf 1=\mathbf 0$ and $\mathrm{rank}(\Sigma_\Delta)=r-1$.
By Theorem~\ref{thm:wald-null-corrected} below, under $H_0$ we have
$T_\Pi\xrightarrow{d}\chi^2_{\,r-1}$; thus for a target significance level
$\alpha_{\mathrm{sig}}$, an asymptotically valid test is to reject $H_0$ if
$T_\Pi>\chi^2_{\,r-1,\,1-\alpha_{\mathrm{sig}}}$, which is the upper $1-\alpha_{\textup{sig}}$ percentile of the $\chi^2_{r-1}$-distribution. In practice, we implement the test statistic $T_\Pi$ by replacing the unknown oracle covariance
$\Sigma_\Delta/N_{\mathrm{eff}}$ with the plug-in estimator
$\widehat V_{\Pi,1}+\widehat V_{\Pi,2}$. We now present the detailed theoretical guarantees for size and power of the proposed alignability test.

\color{black}

\begin{theorem}
\label{thm:wald-null-corrected}
Under Assumptions~\ref{asmp:noise}, \ref{asmp:signal}, and
\ref{asmp:signal-strength}, suppose that $\Sigma_\Delta$ is positive
definite on the simplex tangent space
$T\Delta^{r-1}
:=
\bigl\{
x\in\mathbb R^r:
\mathbf 1^\top x=0
\bigr\}.
$
Then, under
$H_0:\Pi_r(\rho_1)=\Pi_r(\rho_2)$, we have
$
T_\Pi
\xrightarrow{d}
\chi^2_{r-1}.
$
\end{theorem}

\begin{theorem}
\label{thm:wald-power-corrected}
Under the assumptions of
Theorem~\ref{thm:wald-null-corrected}, suppose that
$
\Delta
:=
\Pi_r(\rho_1)-\Pi_r(\rho_2)
\neq 0.
$
Then
$
\frac{T_\Pi}{N_{\mathrm{eff}}}
\xrightarrow{\mathbb P}
c_0
:=
\Delta^\top
\Sigma_\Delta^+
\Delta
>0.
$
Consequently, for any $\varepsilon\in(0,c_0)$,
$
\Pr\left\{
T_\Pi
\ge
N_{\mathrm{eff}}
(c_0-\varepsilon)
\right\}
\longrightarrow 1.
$
Moreover, the $p$-value decays as
$
\exp\left\{
-\frac{N_{\mathrm{eff}}}{2}
(c_0-\varepsilon)
+
o(N_{\mathrm{eff}})
\right\}.
$
\end{theorem}

 \begin{corollary}
\label{cor:sep}
Under the setting of
Theorem~\ref{thm:wald-power-corrected}, let
$\lambda_{\max}(\Sigma_\Delta)$ denote the largest nonzero eigenvalue
of $\Sigma_\Delta$. Then
$
c_0
=
\Delta^\top
\Sigma_\Delta^+
\Delta
\ge
\frac{\|\Delta\|_2^2}
{\lambda_{\max}(\Sigma_\Delta)}.
$
Hence, if
$\|\Pi_r(\rho_1)-\Pi_r(\rho_2)\|_2\ge C>0$, then
$
T_\Pi
\gtrsim
N_{\mathrm{eff}}
\frac{C^2}
{\lambda_{\max}(\Sigma_\Delta)}
$
with high probability.
\end{corollary}

Theorem~\ref{thm:wald-null-corrected} establishes the asymptotic null distribution of the oracle test statistic, providing a principled test for the alignability between two datasets. 
The procedure accommodates heterogeneous noise covariances and unequal sample sizes, while the effective sample size $N_{\mathrm{eff}}$ determines how each dataset contributes. 
Theorem~\ref{thm:wald-power-corrected} establishes power consistency and shows that, under fixed alternatives, the $p$-value decays exponentially in the effective sample size.
By Corollary~\ref{cor:sep},  we have \(c_0 \ge \|\Delta\|_2^2/\lambda_{\max}( \Sigma_\Delta)\), which means larger separation \(\|\Delta\|_2\) between the two datasets yields a faster exponential decay and hence higher power.
In Section~\ref{sec:sim}, we compare our test with alternative methods and demonstrate that it consistently achieves superior performance.

\subsection{Extension to kernel-based manifold learning}
\label{subsec:kernel-mss}

In many applications, similarities between datasets may only become apparent after a nonlinear feature transformation. Our  framework naturally extends to such settings through the use of reproducing kernel Hilbert spaces (RKHS). In this part, we formulate the RKHS version of our dissimilarity measure, working  with Gram matrices to maintain consistency with the covariance-based perspective adopted throughout the paper. We define the key quantities and state the main results, while deferring detailed proofs to Section~S13 of the Supplement. The arguments require kernel spectral theory \citep{ding2024kernel,SmaleZhou2009Geometry} but are analogous to those in the linear setting.

Let $(\X,\mathcal B)$ be a measurable space and
$K:\X\times\X\to\R$ a Mercer kernel with reproducing kernel Hilbert space
$\mathcal H_K$ and feature map $\phi:\X\to\mathcal H_K$ satisfying $
  K(x,y) \;=\; \langle \phi(x),\phi(y)\rangle_{\mathcal H_K},$ and $
  \kappa^2 \;:=\; \sup_{x\in\X} K(x,x) \;<\;\infty.$ 
For dataset $k\in\{1,2\}$, let $\rho_k$ be the distribution of $X_{k}$ on $\X$ and define
the (centered) population kernel covariance operator
$
  M_{K,k}
  \;:=\;
  \E_{X\sim\rho_k}\Bigl[(\phi(X)-m_k)\otimes(\phi(X)-m_k)\Bigr],$ with $
  m_k := \E_{X\sim\rho_k}[\phi(X)],$
where $(u\otimes v)f := \langle f,v\rangle_{\mathcal H_K}\,u$.
By boundedness of $K$, each $M_{K,k}$ is self-adjoint, positive, and trace-class.
Let $(\mathcal D_{k,j})_{j\ge1}$ be the eigenvalues of $M_{K,k}$ arranged in nonincreasing
order. We assume a spiked structure for the leading $r$ components: for some fixed
$r\ge1$ and each $k$, the top $r$ eigenvalues satisfy
$  \mathcal D_{k,1}\ge\cdots\ge\mathcal D_{k,r} > \mathcal D_{k,r+1},$ and $
  \min_{1\le j\le r}
  \bigl(\mathcal D_{k,j}-\mathcal D_{k,j+1}\bigr)
  \;\ge\; \gamma_{k}>0.$

In analogy with Definition~\ref{def:Pi} and Definition~\ref{def:dist}, we define the
top-$r$ kernel principal-variance vector
$
  \mathcal D_{k}^{(r)} := (\mathcal D_{k,1},\ldots,\mathcal D_{k,r})^\top,$ $
  \Pi^{(K)}_r(\rho_k)
  := \frac{\mathcal D_{k}^{(r)}}{\mathbf 1^\top \mathcal D_{k}^{(r)}}
  \;\in\;\Delta^{\,r-1},$
and the {kernelized  spectral dissimilarity}
$
  d^{(K)}_r(\rho_1,\rho_2)
  \;:=\;
  \bigl\|\Pi^{(K)}_r(\rho_1)-\Pi^{(K)}_r(\rho_2)\bigr\|_2.
$
This dissimilarity compares the relative allocation of variance across the leading $r$
eigendirections of $M_{K,k}$ in the feature space $\mathcal H_K$.
It is invariant under unitary transformations of $\mathcal H_K$ and under a common
global rescaling of the kernel $K$.
For the linear kernel $K(x,y)=\langle x,y\rangle$, $d^{(K)}_r$ reduces to the manifold spectral
dissimilarity in Definition~\ref{def:dist}.

Given 
$\{x_{k,1},\ldots,x_{k,N_k}\}$ from $\rho_k$, let $K_k$ be the $N_k\times N_k$ kernel
matrix with entries $[K_k]_{ij}=K(x_{k,i},x_{k,j})$ and let
$\widetilde K_k$ be its standard centered version in feature space.
Define the normalized Gram matrix
$G_k := \frac{1}{N_k}\,\widetilde K_k.$
It is straightforward that the nonzero eigenvalues of $G_k$ coincide with those of the
empirical kernel covariance operator
$
  \widehat M_{K,k}
  := \frac{1}{N_k}\sum_{i=1}^{N_k}
     \bigl(\phi(x_{k,i})-m_k\bigr)\otimes
     \bigl(\phi(x_{k,i})- m_k\bigr),$
so that all spectral quantities of interest can be computed from $G_k$.
Let $\hat{\mathcal D}_{k,1} \ge\cdots\ge\hat{\mathcal D}_{k,r}>0$ be the leading $r$ eigenvalues of $G_k$, and define the empirical kernel spectral profile
$
  \widehat\Pi^{(K)}_{r,k}
  := \frac{(\hat{\mathcal D}_{k,1},\ldots,\hat{\mathcal D}_{k,r})}
          {\sum_{j=1}^r \hat{\mathcal D}_{k,j}}
  \;\in\;\Delta^{\,r-1},
$
which leads to the plug-in estimator  $
  \widehat d^{(K)}_r
  := \bigl\|\widehat\Pi^{(K)}_{r,1}-\widehat\Pi^{(K)}_{r,2}\bigr\|_2.$

We work in the regime where $N_k\to\infty$ with fixed $r$.
Under the bounded-kernel assumption, the spectral perturbation bound
$\|\widehat M_{K,k}-M_{K,k}\|_{\op}=O_{\mathbb P}(N_k^{-1/2})$ holds, and the spiked
eigenvalues of $M_{K,k}$ are stable.
Consequently, for each fixed $j\le r$, we have
$\hat{\mathcal{D}}_{k,j} \xrightarrow{\ \mathbb{P}\ } \mathcal{D}_{k,j},$ $
\widehat\Pi^{(K)}_{r,k} \xrightarrow{\ \mathbb{P}\ } \Pi^{(K)}_r(\rho_k),$ and finally $
\widehat d^{(K)}_r \xrightarrow{\ \mathbb{P}\ } d^{(K)}_r(\rho_1,\rho_2)$.

\section{Simulation and real data applications}

\label{sec:sim}

\subsection{Comparison with existing methods using simulated data}

We evaluate the proposed methods under our manifold signal-plus-noise model (\ref{model}),
where \(X_i\) has i.i.d.\ \(\mathcal N(0,1)\) entries. We let \(\Sigma_i=\mathrm{diag}(\sigma^{(i)}_1,\ldots,\sigma^{(i)}_p)\) be piecewise constant over four blocks (sizes \(\approx p/3,p/6,p/6\), remainder) with levels \((3,4,5,6)\) for \(i=1\) and \((2.5,3,6,4.5)\) for \(i=2\). We fix \(n=r=3\), generate two independent samples of sizes \(N_1=N_2=N\) (so \(N_{\mathrm{eff}}=N/2\)), set \(p=100\), and assess our methods by repeating each simulation setting \(n_{\mathrm{rep}}\) times.
For signal generation, we draw \(V\in\mathbb{R}^{p\times r}\) with orthonormal columns by QR of a standard Gaussian matrix. For dataset \(i\), fix semi–axes \(D_i=(D_{i,1},D_{i,2},D_{i,3})\) and sample independent latent columns on the ellipsoid surface:
\(
x_j^{(i)}=\mathrm{diag}(D_i)\,v_j\) where \( v_j\sim\mathrm{Unif}(\mathbb S^{r-1}).
\)
Set the deterministic signal \(S_{\samp,i}=\sqrt r\,V\,[x^{(i)}_1,\ldots,x^{(i)}_N]\), so that the population signal covariance equals \(V\mathrm{diag}(D_i^2)V^\top\). The two datasets share \(V\) but have independent columns. For evaluation of our proposed test, at the null we take \(D_1=D_2=(7,6,5)\); in power analysis we increase \(D_{2,1}\) anisotropically.
We first conduct simulation studies to validate the theoretical predictions about the proposed statistical test.

\emph{Validation 1: Null calibration.}
Set $D_1=D_2=(7,6,5)$, $N_1=N_2=N=1500$, and
$n_{\mathrm{rep}}=800$.
Under the null hypothesis $H_0:\Pi_1=\Pi_2$, we compare the empirical distribution of the feasible statistic $T_\Pi$ with the oracle $\chi^2_{r-1}$ reference, where $r-1=2$.  The empirical distribution of $T_\Pi$ aligns closely with the $\chi^2_2$ reference distribution, as shown in Figure~S3. The goodness-of-fit test gives a $p$-value of $0.867$, and the estimated size at level $\alpha=0.05$ is $0.044$. The empirical quantiles are also close to the oracle-reference quantiles in Table~\MainTableone.

\emph{Validation 2: Power under anisotropic separation.}
We fix $N_1=N_2=N=3000$ and increase the first semi-axis of dataset~2 by a factor
$c\in\{1.00,1.05,1.10,1.20,1.30,1.50\}$:
$D_2=\sqrt{\operatorname{diag}(c,1,1)}\,D_1$.
For each $c$, we report the separation $\|\Delta\Pi\|$, an oracle-based noncentral chi-square approximation, and the empirical power of our procedure at level $0.05$ over $n_{\mathrm{rep}}=600$ replications. The comparison in Supplementary Table~\MainTabletwo\ and Figure~S4 is a numerical assessment of how closely our procedure follows the oracle-based power approximation. Power rises monotonically with \(\|\Delta\Pi\|\); even for modest anisotropy (\(c=1.10\)) the test attains power \(\approx\) 0.53 to 0.55 and quickly approaches one once \(c\ge 1.2\).

Next we compare the proposed test (nMSD test) with several related methods and summarize two representative scenarios using the same data generator.
In \emph{Scenario A (noise only)}, the signal manifold is fixed across datasets while the noise level is scaled: $
D_2 = D_1,$ $\Sigma_2 = \alpha\,\Sigma_1,$ $\alpha \in [1,2].$ 
In \emph{Scenario B (signal reweighting)}, we anisotropically amplify the first signal axis (with scale removed) while keeping the noise profile fixed: $
D_2 = \sqrt{\mathrm{diag}(c,1,1)}\,D_1,$ $c \ge 1,$ $\Sigma_2 = \Sigma_1.$
All other aspects (dimension, rank, and sampling scheme) follow the main Monte Carlo design.

\begin{figure}[t]
  \centering

  \begin{minipage}{0.4\linewidth}
    \centering
    \includegraphics[width=\linewidth]{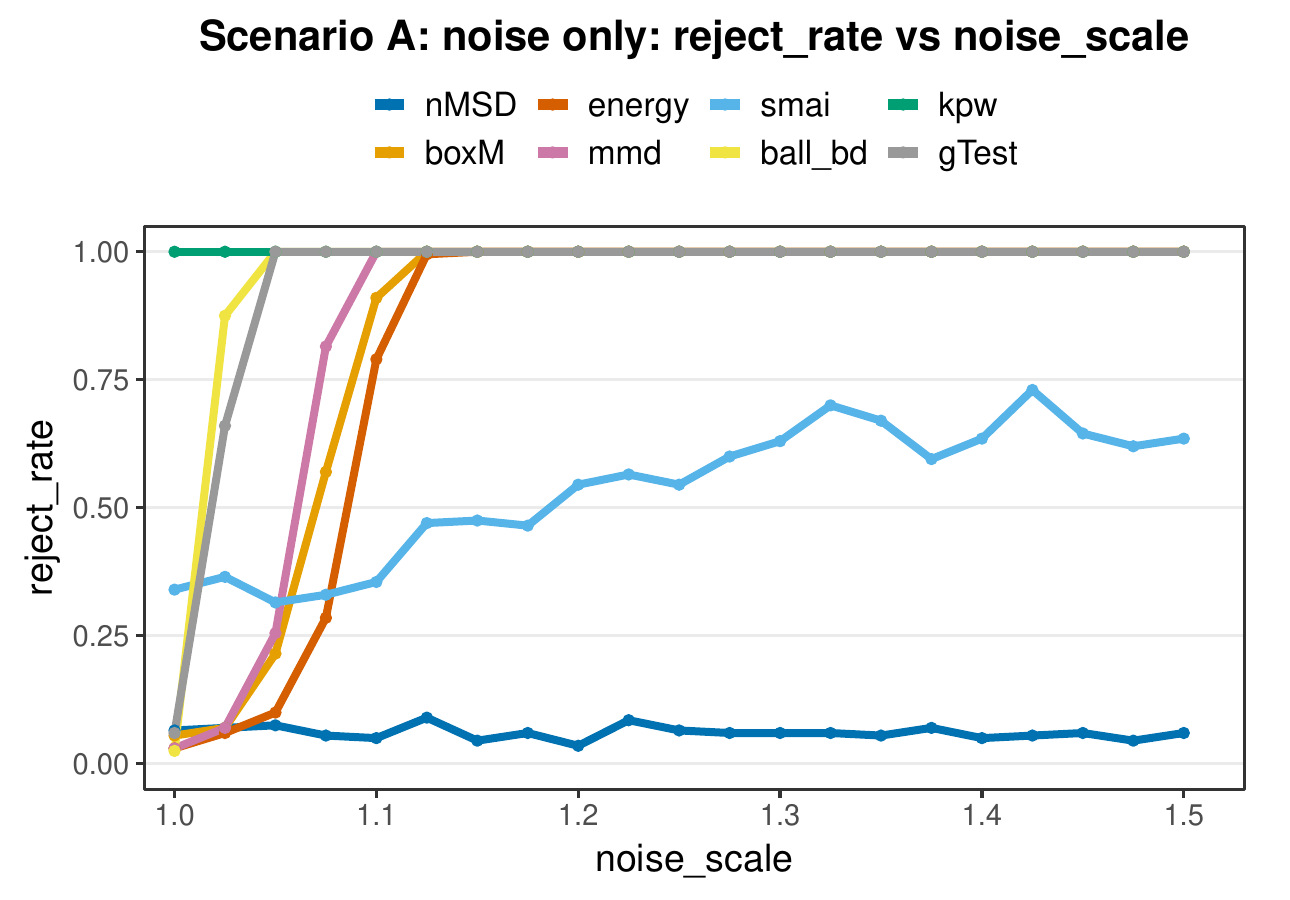}
    \\ \small (a)
  \end{minipage}
  \begin{minipage}{0.4\linewidth}
    \centering
    \includegraphics[width=\linewidth]{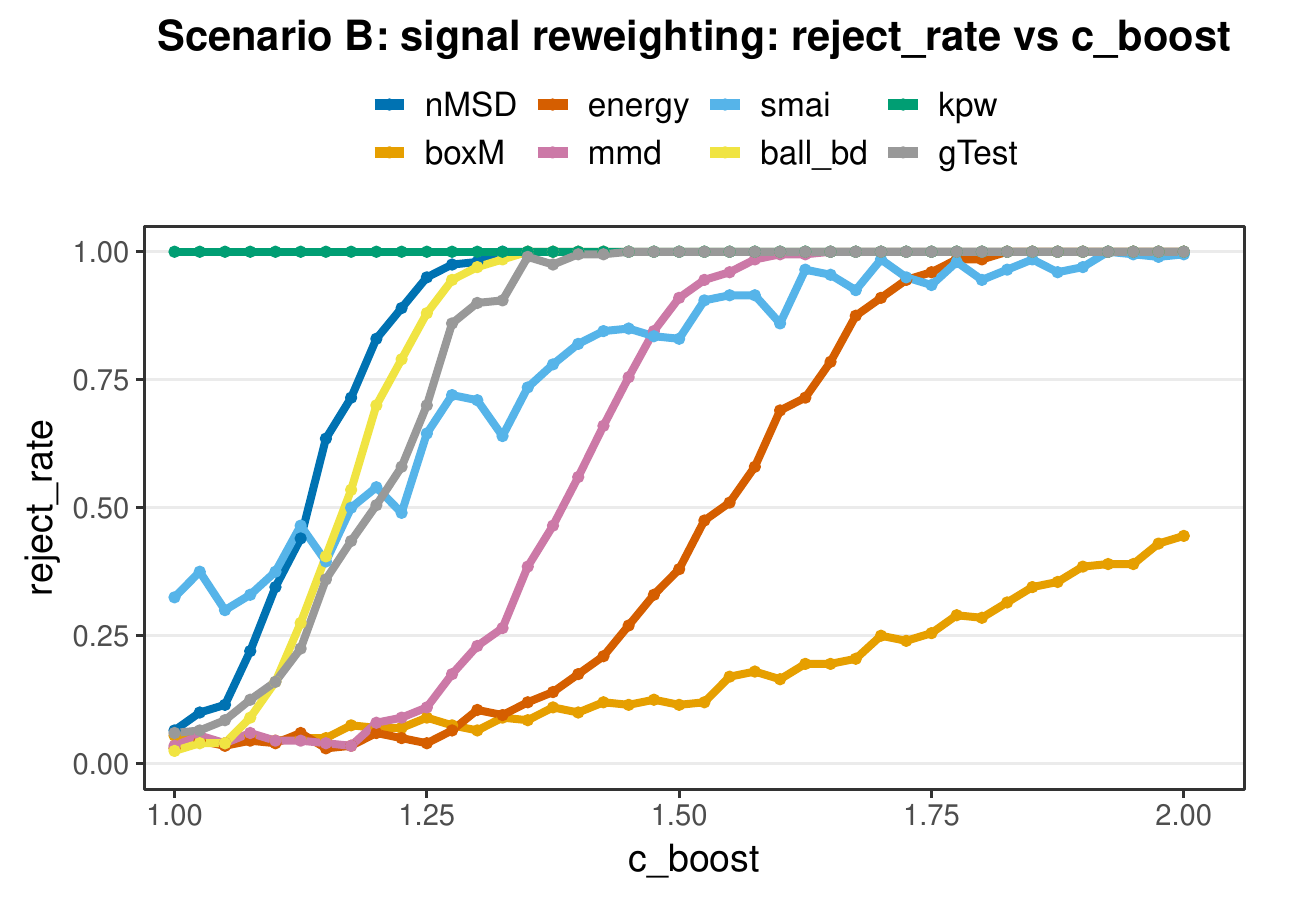}
    \\ \small (b)
  \end{minipage}
   \centering

  \caption{\textbf{Simulation summary across methods.} 
 Panel (a) shows the rejection rates for identical signals under varying noise scales, and Panel (b) shows the rejection rates under fixed noise with increasing signal separation.}
  \label{fig:sim-methods-summary}
\end{figure}

Figure~\ref{fig:sim-methods-summary} compares seven competitors---Box’s \(M\) \citep{box1949M}, Energy distance \citep{energydistance1},  SMAI \citep{smai}, MMD (Gaussian kernel) \citep{kernelmmd1}, Ball Divergence (BD) \citep{Pan2018BallDiv}, Kernel Projected Wasserstein (KPW) \citep{wasserstein1}, and a graph–based  test (gTest) \citep{graphtest1}---against the proposed nMSD test. Detailed comparisons of mean \(p\)-values are provided in Supplementary Figure \Mainfigp.

Panel~(a) varies only the noise profile while keeping the signal manifold fixed; Panel~(b) increases the anisotropy of the first signal axis (\(c-boost\)) while keeping the noise profile fixed. In the noise-only sweep (Panel~(a)), the nMSD test remains well calibrated: as shown in the Supplement, its mean \(p\)-values stay near \(0.5\), while its rejection rate in the main figure tracks the nominal size. By contrast, omnibus tests that target any covariance or distributional difference quickly move to very small mean \(p\)-values (Supplementary Figure \Mainfigp) and near-unit rejection (Panel~(a)) even under modest heteroskedasticity. The Type I error rate of SMAI also increases significantly as the noise disparity between the two groups grows.

When the signal proportion changes (Panel~(b)), the nMSD test shows a steep power ramp: as \(c\) increases, the mean \(p\)-value collapses (see Supplement) and the rejection rate approaches one by \(c\approx 1.3\text{--}1.4\), matching its design to detect reallocation of variance across the top-\(r\) eigenvectors. BD and gTest also gain power but require somewhat larger separation; MMD and KPW improve more gradually. Box's \(M\) and Energy rise only slowly with \(c\) in this setting: their behavior is driven by omnibus sensitivity rather than the proportional-signal alternative. Together, we observe that the nMSD test is stable under noise changes yet sharply responsive to signal-axis reallocation, whereas omnibus tests confound noise and signal, and purely geometric alignment such as SMAI can miss proportional shifts.

The intended null equality of signal energy proportions in the top-\(r\) directions despite differing noise—is respected by nMSD test, which is insensitive to pure noise changes yet responsive to reallocation of signal variance. Baselines that test broader hypotheses conflate noise and signal and thus over-reject when only noise differs; purely geometric alignments tend to be conservative for axis reweighting.

\section{Real data analysis} \label{sec:real}

\subsection{$\Pi_r(\rho)$ captures a hepatocyte zonation-like continuum}

We first used the integrated GepLiver single-cell RNA-seq dataset \citep{Li2023GepLiver} to test
whether the proposed method captures a known biological structure in hepatocyte populations.
The dataset contains 409,775 cells from 349 liver samples, including 57,022
hepatocytes (a liver cell type) from 294 samples. In this analysis, we focused on the HCC malignant hepatocytes, which consist of $n=29,505$ cells and formed the
main malignant hepatocyte population. Because hepatocyte gene-expression programs vary
continuously along the porto-central axis of the liver lobule
\citep{Halpern2017,Aizarani2019}, we used hepatocyte zonation as a cell-level biological label for subsequent evaluations.
The zonation marker panel was curated a priori from published hepatocyte-zonation studies.
For each cell, log-normalized expression of curated pericentral and
periportal marker genes was standardized across hepatocytes. The zonation
score was defined as
$
s_i =
\frac{1}{|C|}\sum_{g \in C} z_{ig}
-
\frac{1}{|P|}\sum_{g \in P} z_{ig},$
where $C$ and $P$ denote the pericentral and periportal marker sets,
respectively. Higher scores indicate a more pericentral-like state, whereas
lower scores indicate a more periportal-like state.

\begin{figure}[t]
  \centering

  \begin{minipage}[t]{0.4\linewidth}
    \centering
    \includegraphics[
      width=\linewidth
    ]{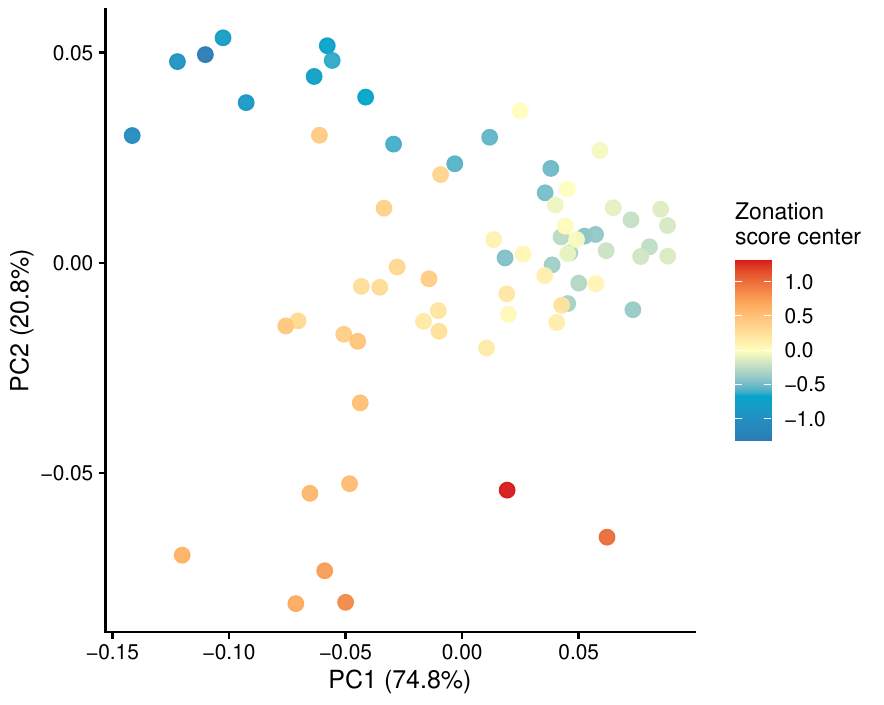}
    \caption{\textbf{PCA embedding of top-6 principal variance profiles computed from 70 zonation-ordered hepatocyte bins.}
}
    \label{fig:k70}
  \end{minipage}
  \begin{minipage}[t]{0.4\linewidth}
    \centering
    \includegraphics[
      width=\linewidth
    ]{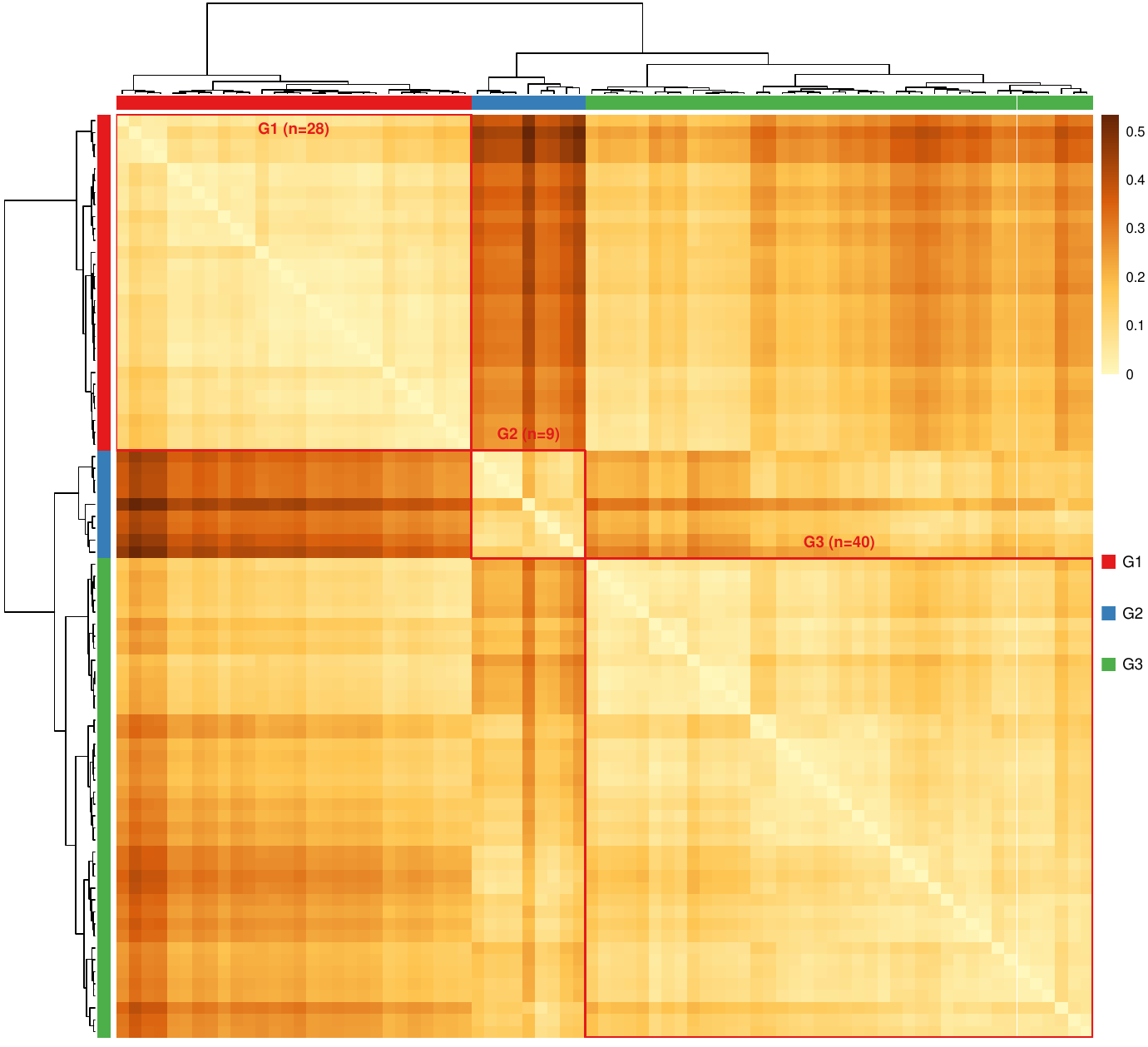}
    \caption{\textbf{Clustering results.}}
    \label{fig:cluster}
  \end{minipage}

\end{figure}

Cells were ordered by $s_i$ and divided into 70 non-overlapping bins of equal size. For each bin, we estimated the top-$r$ normalized principal variances $\Pi_r(\rho)$ and stacked them into an $r\times 70$ matrix, with $r=6$ chosen according to Remark~2.1. We then obtained a PCA embedding of the 70 bins based on their top-$r$ principal variance profiles. As shown in Figure~\ref{fig:k70}, the bin-level profiles exhibited a continuous trajectory in the low-dimensional embedding that closely followed the median zonation score, which represents the underlying biological ordering of the bins. These results demonstrate that the proposed method can successfully recover the latent population-level hepatocyte zonation structure.

\begin{figure}
  \centering
  \includegraphics[width=0.9\linewidth]{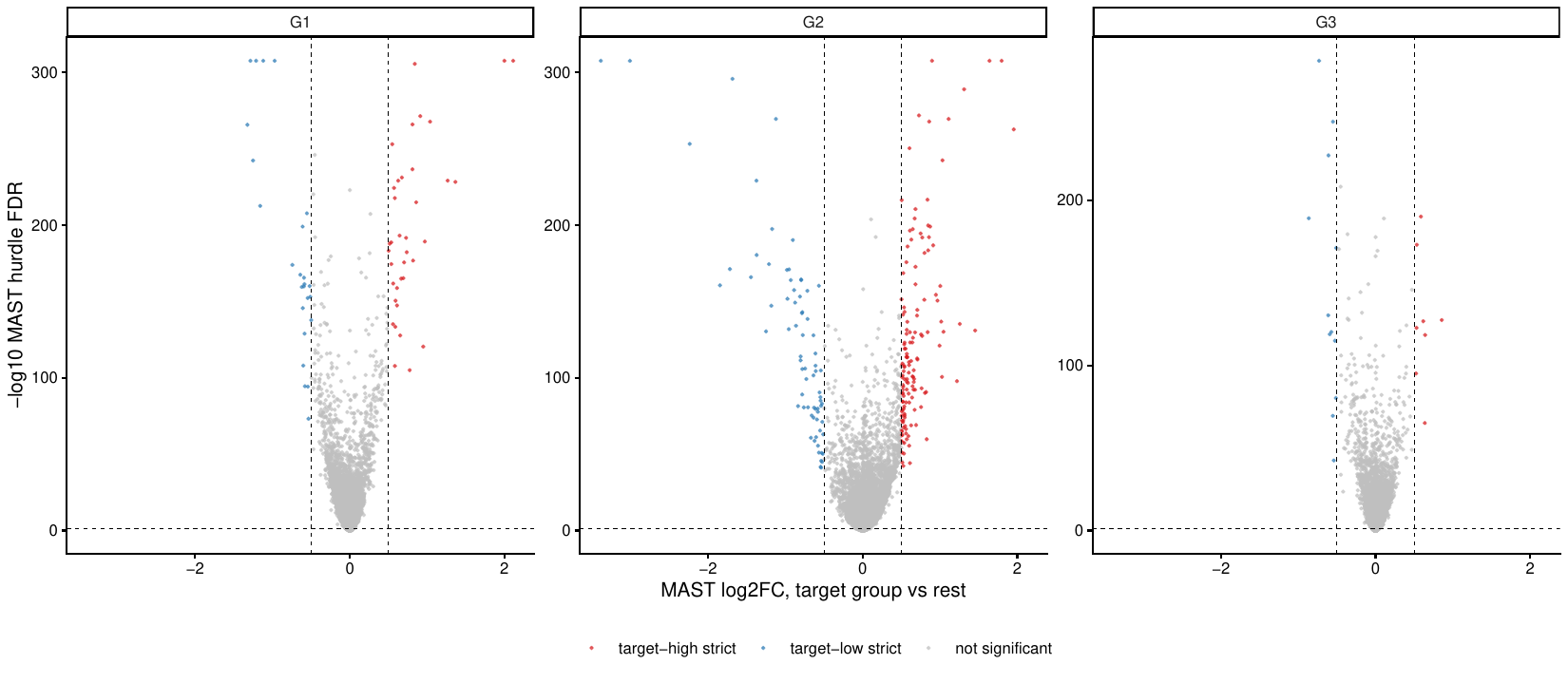}
  \caption{one-versus-rest differential expression for the three MSD-defined hepatocyte groups.}
  \label{fig:volcano}
\end{figure}
\begin{figure}
  \centering
  \includegraphics[width=0.9\linewidth]{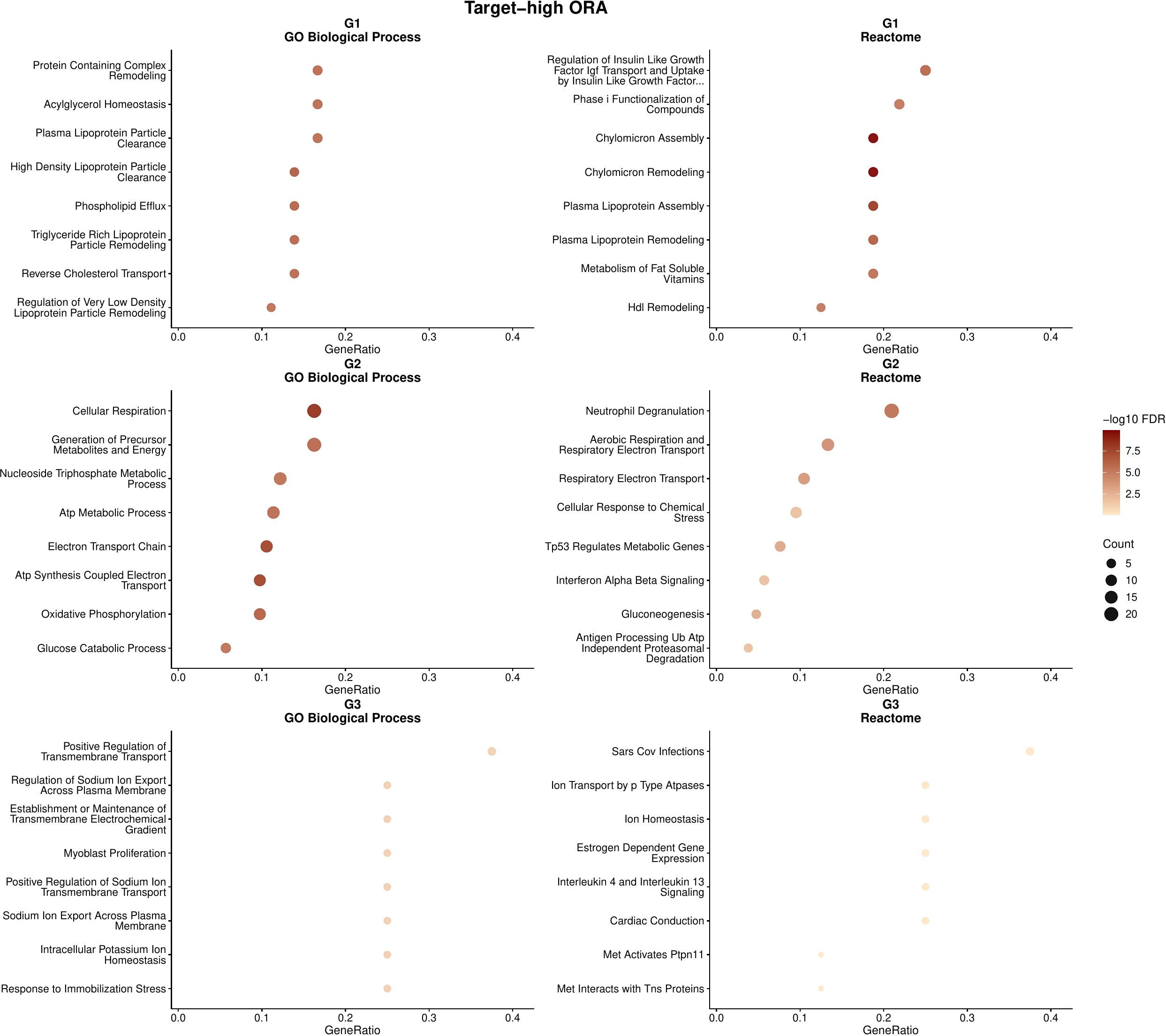}
  \caption{Reactome and GO Biological Process pathway enrichment for the three MSD-defined hepatocyte groups.}
  \label{fig:hcc_pathway}
\end{figure}

\subsection{MSD identifies candidate hepatocyte subtypes in HCC}

For the HCC subtype analysis, we retained samples containing at least 100
hepatocytes. This yielded 77 HCC samples from three GepLiver cohorts:
SC12/GSE156625 ($n=16$), SC13/SRP318499 ($n=8$), and
SC14/HRA001748 ($n=53$). The final dataset contained 28,308 hepatocytes,
with 108--984 cells per sample and a median of 311 cells per sample.

We applied MSD to characterize sample-level differences in hepatocyte
population structure. A common set of 750 highly variable genes was used
across all samples, and each sample was represented by a rank-6 normalized
principal variance profile. Hierarchical clustering of Euclidean distances
between these profiles using Ward's criterion identified three groups:
G1 ($n=28$), G2 ($n=9$), and G3 ($n=40$)
(Figure~\ref{fig:cluster}). The groups showed little dependence on cohort
and were not explained by differences in sample-level mean or median
zonation scores.

To functionally characterize the MSD-defined groups, we performed
one-versus-rest differential-expression analysis. This identified 62, 215,
and 20 strict marker genes for G1, G2, and G3, respectively
(Figure~\ref{fig:volcano}). Reactome and Gene Ontology Biological Process
(GO BP) over-representation analyses were performed separately for genes
expressed at higher and lower levels in each group. Results for
higher-expression genes are shown in Figure~\ref{fig:hcc_pathway}, with
lower-expression results provided in Supplementary
Figure~\Mainpathwaylow.

G1 showed consistent enrichment of lipoprotein assembly, transport, and
remodeling pathways, indicating retention of specialized hepatocyte
secretory and lipid-handling functions. This interpretation is consistent
with the established role of primary hepatocytes in lipoprotein metabolism
\citep{Ling2013LipoproteinMetabolism} and with differentiated
hepatocyte-like programs previously described in HCC
\citep{Hoshida2009HCCSubclasses}.

G2 was characterized by increased mitochondrial respiration, ATP
production, and glucose-related metabolism. In contrast, genes expressed
at lower levels in G2 were strongly enriched for lipoprotein-processing
and plasma-protein pathways. G2 therefore represents a bioenergetically
active state with reduced expression of specialized hepatocyte secretory
functions, consistent with the metabolic heterogeneity and reprogramming
reported in HCC \citep{Bidkhori2018MetabolicSubtypes,TCGA2017HCC}.

G3 showed a substantially weaker positive transcriptional signal. Its
higher-expression genes were modestly associated with ion transport and
homeostasis, but no GO BP or Reactome term reached the FDR threshold of
0.05. Genes expressed at lower levels were enriched for acute-phase,
inflammatory, and antigen-presentation-related processes. G3 therefore
represents a less clearly resolved candidate state, characterized mainly
by reduced acute-phase and immune-associated expression.

Overall, MSD identified three biologically interpretable HCC hepatocyte
groups. The strongest contrast separated a lipoprotein- and
secretory-function-associated state in G1 from a bioenergetic state with
reduced specialized hepatocyte functions in G2, whereas G3 showed a
weaker and less distinct functional program. These downstream
transcriptional differences support the biological relevance of the
sample-level population structure captured by MSD.

\color{black}

\section*{Supplementary Material}
The supplementary material \citep{ChenMa-supp} contains technical proofs of the main theorems, some examples and additional discussions.

\section*{Acknowledgment}

This work was carried out primarily while H.C. was an undergraduate student at Weiyang College, Tsinghua University.

\bibliographystyle{plainnat}
\bibliography{reference}

@article{Ling2013LipoproteinMetabolism,
  author  = {Ling, Ji and Lewis, Jamie and Douglas, Donna and Kneteman, Norman M. and Vance, Dennis E.},
  title   = {Characterization of lipid and lipoprotein metabolism in primary human hepatocytes},
  journal = {Biochimica et Biophysica Acta (BBA) - Molecular and Cell Biology of Lipids},
  year    = {2013},
  volume  = {1831},
  number  = {2},
  pages   = {387--397},
  doi     = {10.1016/j.bbalip.2012.08.012}
}

@article{Li2023GepLiver,
  author  = {Li, Ziteng and Zhang, Hena and Li, Qin and Feng, Wanjing and Jia, Xiya and Zhou, Runye and Huang, Yi and Li, Yan and Hu, Zhixiang and Hu, Xichun and Zhu, Xiaodong and Huang, Shenglin},
  title   = {GepLiver: an integrative liver expression atlas spanning developmental stages and liver disease phases},
  journal = {Scientific Data},
  year    = {2023},
  volume  = {10},
  pages   = {376},
  doi     = {10.1038/s41597-023-02257-1}
}

@article{Hoshida2009HCCSubclasses,
  author  = {Hoshida, Yujin and Nijman, Sebastian M. B. and Kobayashi, Masahiro and Chan, Jennifer A. and Brunet, Jean-Philippe and Chiang, Derek Y. and Villanueva, Augusto and Newell, Philippa and Ikeda, Kenji and Hashimoto, Masaji and Watanabe, Goro and Gabriel, Stacey and Friedman, Scott L. and Kumada, Hiromitsu and Llovet, Josep M. and Golub, Todd R.},
  title   = {Integrative transcriptome analysis reveals common molecular subclasses of human hepatocellular carcinoma},
  journal = {Cancer Research},
  year    = {2009},
  volume  = {69},
  number  = {18},
  pages   = {7385--7392},
  doi     = {10.1158/0008-5472.CAN-09-1089}
}

@article{Bidkhori2018MetabolicSubtypes,
  author  = {Bidkhori, Gholamreza and Benfeitas, Rui and Klevstig, Martina and Zhang, Cheng and Nielsen, Jens and Uhl{\'e}n, Mathias and Bor{\'e}n, Jan and Mardinoglu, Adil},
  title   = {Metabolic network-based stratification of hepatocellular carcinoma reveals three distinct tumor subtypes},
  journal = {Proceedings of the National Academy of Sciences of the United States of America},
  year    = {2018},
  volume  = {115},
  number  = {50},
  pages   = {E11874--E11883},
  doi     = {10.1073/pnas.1807305115}
}

@article{TCGA2017HCC,
  author  = {{Cancer Genome Atlas Research Network}},
  title   = {Comprehensive and integrative genomic characterization of hepatocellular carcinoma},
  journal = {Cell},
  year    = {2017},
  volume  = {169},
  number  = {7},
  pages   = {1327--1341.e23},
  doi     = {10.1016/j.cell.2017.05.046}
}

@article{Halpern2017,
  author  = {Halpern, Keren Bahar and Shenhav, Rom and
             Matcovitch-Natan, Orit and others},
  title   = {Single-cell spatial reconstruction reveals global division
             of labour in the mammalian liver},
  journal = {Nature},
  year    = {2017},
  volume  = {542},
  pages   = {352--356},
  doi     = {10.1038/nature21065}
}

@article{Aizarani2019,
  author  = {Aizarani, Nadim and Saviano, Antonio and Sagar and others},
  title   = {A human liver cell atlas reveals heterogeneity and
             epithelial progenitors},
  journal = {Nature},
  year    = {2019},
  volume  = {572},
  pages   = {199--204},
  doi     = {10.1038/s41586-019-1373-2}
}

@article{Linderman2022ALRA,
  author  = {Linderman, George C. and Zhao, Jun and Roulis, Manolis and
             Bielecki, Piotr and Flavell, Richard A. and Nadler, Boaz and
             Kluger, Yuval},
  title   = {Zero-preserving imputation of single-cell {RNA}-seq data},
  journal = {Nature Communications},
  year    = {2022},
  volume  = {13},
  number  = {1},
  pages   = {192},
  doi     = {10.1038/s41467-021-27729-z},
  url     = {https://doi.org/10.1038/s41467-021-27729-z}
}

@article{Aparicio2020RMT,
  author  = {Aparicio, Luis and Bordyuh, Mykola and Blumberg, Andrew J. and
             Rabadan, Raul},
  title   = {A Random Matrix Theory Approach to Denoise Single-Cell Data},
  journal = {Patterns},
  year    = {2020},
  volume  = {1},
  number  = {3},
  pages   = {100035},
  doi     = {10.1016/j.patter.2020.100035},
  url     = {https://doi.org/10.1016/j.patter.2020.100035}
}

@article{Verma2020Manifold,
  author  = {Verma, Archit and Engelhardt, Barbara E.},
  title   = {A robust nonlinear low-dimensional manifold for single cell
             {RNA}-seq data},
  journal = {BMC Bioinformatics},
  year    = {2020},
  volume  = {21},
  number  = {1},
  pages   = {324},
  doi     = {10.1186/s12859-020-03625-z},
  url     = {https://doi.org/10.1186/s12859-020-03625-z}
}

@Manual{panc8data,
  title  = {panc8.SeuratData: Eight Pancreas Datasets Across Five Technologies},
  author = {{Satija Lab}},
  year   = {2019},
  note   = {R package version 3.0.2},
}

@article{dai2018principal,
  title={PRINCIPAL COMPONENT ANALYSIS FOR FUNCTIONAL DATA ON RIEMANNIAN MANIFOLDS AND SPHERES},
  author={Dai, Xiongtao and M{\"u}ller, Hans-Georg},
  journal={The Annals of Statistics},
  volume={46},
  number={6B},
  pages={3334--3361},
  year={2018},
}

@article{huckemann2006principal,
  title={Principal component analysis for Riemannian manifolds, with an application to triangular shape spaces},
  author={Huckemann, Stephan and Ziezold, Herbert},
  journal={Advances in Applied Probability},
  volume={38},
  number={2},
  pages={299--319},
  year={2006},
}

@article{donoho2023screenot,
  title={ScreeNOT: Exact MSE-optimal singular value thresholding in correlated noise},
  author={Donoho, David and Gavish, Matan and Romanov, Elad},
  journal={The Annals of Statistics},
  volume={51},
  number={1},
  pages={122--148},
  year={2023},
  publisher={Institute of Mathematical Statistics}
}

@misc{ChenMa-supp,
  author = {Hongrui Chen and Rong Ma},
  year   = {2026},
  title  = {Supplement to ``Statistical Inference of Manifold Similarity and Alignability across Noisy High-Dimensional Datasets''},
}

@article{jiang2022universal,
  title   = {A Universal Test on Spikes in a High-Dimensional Generalized Spiked Model and Its Applications},
  author  = {Jiang, Dandan},
  journal = {Statistica Sinica},
  volume  = {33},
  year    = {2023},
  pages   = {1749--1770},

}

@article{luo2020unsupervised,
  title={Unsupervised domain adaptation via discriminative manifold propagation},
  author={Luo, You-Wei and Ren, Chuan-Xian and Dai, Dao-Qing and Yan, Hong},
  journal={IEEE transactions on pattern analysis and machine intelligence},
  volume={44},
  number={3},
  pages={1653--1669},
  year={2020},
}

@article{DingMa2023,
  title   = {Learning Low-Dimensional Nonlinear Structures from High-Dimensional Noisy Data:
             An Integral Operator Approach},
  author  = {Ding, Xiucai and Ma, Rong},
  journal = {The Annals of Statistics},
  volume  = {51},
  number  = {4},
  pages   = {1744--1769},
  year    = {2023},
}

@article{bao2021singular,
  title={SINGULAR VECTOR AND SINGULAR SUBSPACE DISTRIBUTION FOR THE MATRIX DENOISING MODEL},
  author={Bao, Zhigang and Ding, Xiucai and Wang, Ke},
  journal={The Annals of Statistics},
  volume={49},
  number={1},
  pages={370--392},
  year={2021},
  publisher={JSTOR}
}

@book{bai2010spectral,
  title={Spectral analysis of large dimensional random matrices},
  author={Bai, Zhidong and Silverstein, Jack W},
  volume={20},
  year={2010},
  publisher={Springer}
}

@book{erdHos2017dynamical,
  title={A dynamical approach to random matrix theory},
  author={Erd{\H{o}}s, L{\'a}szl{\'o} and Yau, Horng-Tzer},
  volume={28},
  year={2017},
  publisher={American Mathematical Soc.}
}

@article{yang2022limiting,
  title={Limiting distribution of the sample canonical correlation coefficients of high-dimensional random vectors},
  author={Yang, Fan},
  journal={Electronic Journal of Probability},
  volume={27},
  pages={1--71},
  year={2022},
}

@article{bao2019canonical,
  title={Canonical correlation coefficients of high-dimensional Gaussian vectors: Finite rank case},
  author={Bao, Zhigang and Hu, Jiang and Pan, Guangming and Zhou, Wang},
  journal={The Annals of Statistics},
  volume={47},
  number={1},
  pages={612--640},
  year={2019},
}

@article{liu2013gaussian,
  title={GAUSSIAN GRAPHICAL MODEL ESTIMATION WITH FALSE DISCOVERY RATE CONTROL},
  author={LIU, WEIDONG},
  journal={The Annals of Statistics},
  volume={41},
  number={6},
  pages={2948--2978},
  year={2013}
}

@article{xia2018multiple,
  title={Multiple testing of submatrices of a precision matrix with applications to identification of between pathway interactions},
  author={Xia, Yin and Cai, Tianxi and Cai, T Tony},
  journal={Journal of the American Statistical Association},
  volume={113},
  number={521},
  pages={328--339},
  year={2018},
}

@article{xia2015testing,
  title={Testing differential networks with applications to the detection of gene-gene interactions},
  author={Xia, Yin and Cai, Tianxi and Cai, T Tony},
  journal={Biometrika},
  volume={102},
  number={2},
  pages={247--266},
  year={2015},
}

@article{cai2013two,
  title={Two-sample covariance matrix testing and support recovery in high-dimensional and sparse settings},
  author={Cai, Tony and Liu, Weidong and Xia, Yin},
  journal={Journal of the American Statistical Association},
  volume={108},
  number={501},
  pages={265--277},
  year={2013},
}

@article{li2023eigenvalue,
  title={Eigenvalue distribution of a high-dimensional distance covariance matrix with application},
  author={Li, Weiming and Wang, Qinwen and Yao, Jianfeng},
  journal={Statistica Sinica},
  volume={33},
  number={1},
  pages={149--168},
  year={2023},
}

@article{yan2024inference,
  title={Inference for heteroskedastic PCA with missing data},
  author={Yan, Yuling and Chen, Yuxin and Fan, Jianqing},
  journal={The Annals of Statistics},
  volume={52},
  number={2},
  pages={729--756},
  year={2024},
}

@inproceedings{mestre2011asymptotic,
  title={Asymptotic analysis of a consistent subspace estimator for observations of increasing dimension},
  author={Mestre, Xavier and Vallet, Pascal and Loubaton, Philippe and Hachem, Walid},
  booktitle={2011 IEEE Statistical Signal Processing Workshop (SSP)},
  pages={677--680},
  year={2011},
  organization={IEEE}
}

@article{wu2018think,
  title={THINK GLOBALLY, FIT LOCALLY UNDER THE MANIFOLD SETUP: ASYMPTOTIC ANALYSIS OF LOCALLY LINEAR EMBEDDING},
  author={Wu, Hau-Tieng and Wu, Nan},
  journal={The Annals of Statistics},
  volume={46},
  number={6B},
  pages={3805--3837},
  year={2018},
}

@article{ding2024kernel,
  title={How do kernel-based sensor fusion algorithms behave under high-dimensional noise?},
  author={Ding, Xiucai and Wu, Hau-Tieng},
  journal={Information and Inference: A Journal of the IMA},
  volume={13},
  number={1},
  pages={iaad051},
  year={2024},
}

@article{ding2025kernel,
  title={Kernel spectral joint embeddings for high-dimensional noisy datasets using duo-landmark integral operators},
  author={Ding, Xiucai and Ma, Rong},
  journal={Journal of the American Statistical Association},
  pages={1--28},
  year={2025},
}

@article{meilua2024manifold,
  title={Manifold learning: What, how, and why},
  author={Meil{\u{a}}, Marina and Zhang, Hanyu},
  journal={Annual Review of Statistics and Its Application},
  volume={11},
  number={1},
  pages={393--417},
  year={2024},
}

@article{donoho2014minimax,
  title={Minimax risk of matrix denoising by singular value thresholding},
  author={Donoho, David and Gavish, Matan},
  journal={The Annals of Statistics},
  volume={42},
  number={6},
  year={2014},
}

@article{bengio2013representation,
  title={Representation learning: A review and new perspectives},
  author={Bengio, Yoshua and Courville, Aaron and Vincent, Pascal},
  journal={IEEE transactions on pattern analysis and machine intelligence},
  volume={35},
  number={8},
  pages={1798--1828},
  year={2013},
}

@article{fefferman2016testing,
  title={Testing the manifold hypothesis},
  author={Fefferman, Charles and Mitter, Sanjoy and Narayanan, Hariharan},
  journal={Journal of the American Mathematical Society},
  volume={29},
  number={4},
  pages={983--1049},
  year={2016}
}

@article{LinPanZhaoZhou2024,
  title        = {Asymptotic Distribution of Spiked Eigenvalues in the Large Signal-Plus-Noise Model},
  author       = {Lin, Zeqin and Pan, Guangming and Zhao, Peng and Zhou, Jia},
  journal      = {arXiv preprint arXiv:2401.11672},
  year         = {2024}
}

@article{energydistance2,
  author  = {Chu, Lynna and Dai, Xiongtao},
  title   = {Manifold energy two-sample test},
  journal = {Electronic Journal of Statistics},
  year    = {2024},
  volume  = {18},
  pages   = {145--166},

}

@article{box1949M,
  author  = {Box, G. E. P.},
  title   = {A General Distribution Theory for a Class of Likelihood Criteria},
  journal = {Biometrika},
  year    = {1949},
  volume  = {36},
  number  = {3--4},
  pages   = {317--346}
}

@article{kernelmmd1,
  author  = {Gretton, Arthur and Borgwardt, Karsten M. and Rasch, Malte J. and Sch{\"o}lkopf, Bernhard and Smola, Alexander J.},
  title   = {A Kernel Two-Sample Test},
  journal = {Journal of Machine Learning Research},
  year    = {2012},
  volume  = {13},
  pages   = {723--773}
}

@article{kernelmmd3,
  author  = {Balasubramanian, Krishnakumar and Li, Tong and Yuan, Ming},
  title   = {On the Optimality of Kernel-Embedding Based Goodness-of-Fit Tests},
  journal = {Journal of Machine Learning Research},
  year    = {2021},
  volume  = {22},
  pages   = {1--45}
}

@article{Pan2018BallDiv,
  author  = {Pan, Wenliang and Tian, Yuan and Wang, Xueqin and Zhang, Heping},
  title   = {Ball Divergence: Nonparametric Two Sample Test},
  journal = {The Annals of Statistics},
  year    = {2018},
  volume  = {46},
  number  = {3},
  pages   = {1109--1137},
}

@article{energydistance1,
  author  = {Rizzo, Maria L. and Sz{\'e}kely, G{\'a}bor J.},
  title   = {Energy distance},
  journal = {Wiley Interdisciplinary Reviews: Computational Statistics},
  year    = {2016},
  volume  = {8},
  pages   = {27--38},
  doi     = {10.1002/wics.1375}
}

@article{graphtest1,
  author  = {Chen, Hao and Chen, Xu and Su, Yi},
  title   = {A weighted edge-count two-sample test for multivariate and object data},
  journal = {arXiv preprint arXiv:1604.06515},
  year    = {2016}
}

@article{graphtest2,
  author  = {Chen, Hao and Friedman, Jerome H.},
  title   = {A new graph-based two-sample test for multivariate and object data},
  journal = {Journal of the American statistical association},
  year    = {2016}
}

@article{interpoint2,
  author  = {Li, Jun},
  title   = {Asymptotic normality of interpoint distances for high-dimensional data with applications to the two-sample problem},
  journal = {Biometrika},
  year    = {2018},
  volume  = {105},
  number  = {3},
  pages   = {529--546},
}

@article{rmt1,
  author  = {Ding, Xiucai and Yang, Fan},
  title   = {Spiked separable covariance matrices and principal components},
  journal = {The Annals of Statistics},
  year    = {2021},
  volume  = {49},
  number  = {2},
  pages   = {1113--1138},
}

@inproceedings{wasserstein1,
  author    = {Wang, Jie and Gao, Rui and Xie, Yao},
  title     = {Two-Sample Test with Kernel Projected Wasserstein Distance},
  booktitle = {Proceedings of the 25th International Conference on Artificial Intelligence and Statistics (AISTATS)},
  series    = {Proceedings of Machine Learning Research},
  volume    = {151},
  pages     = {8003--8021},
  year      = {2022},
  publisher = {PMLR}
}

@article{smai,
  author  = {Ma, R. and Sun, E. D. and Donoho, D. and Zou, J.},
  title   = {Principled and interpretable alignability testing and integration of single-cell data},
  journal = {Proc Natl Acad Sci U S A},
  year    = {2024},
  volume  = {121},
  number  = {10},
  pages   = {e2313719121},
  month   = {March},
}

@article{SmaleZhou2009Geometry,
  title   = {Geometry on Probability Spaces},
  author  = {Smale, Steve and Zhou, Ding-Xuan},
  journal = {Constructive Approximation},
  volume  = {30},
  number  = {3},
  pages   = {311--323},
  year    = {2009},

}

@article{Boysen_2009,
   title={Consistencies and rates of convergence of jump-penalized least squares estimators},
   volume={37},
   number={1},
   journal={The Annals of Statistics},
   author={Boysen, Leif and Kempe, Angela and Liebscher, Volkmar and Munk, Axel and Wittich, Olaf},
   year={2009},
}

@article{Weinmann2014L1Potts,
author = {Weinmann, Andreas and Storath, Martin and Demaret, Laurent},
title = {The L1-Potts Functional for Robust Jump-Sparse Reconstruction},
journal = {SIAM Journal on Numerical Analysis},
volume = {53},
number = {1},
pages = {644-673},
year = {2015},


}

@book{Janssens2007CME,
  title     = {Computational Materials Engineering: An Introduction to Microstructure Evolution},
  editor    = {Janssens, Koenraad G. F. and Raabe, Dierk and Kozeschnik, Ernst and Miodownik, Mark A. and Nestler, Britta},
  year      = {2007},
  publisher = {Academic Press},
  address   = {Burlington, MA},
  isbn      = {978-0-12-369468-3},
  url       = {https://www.sciencedirect.com/book/9780123694683/computational-materials-engineering}
}

@article{Leeb2021MatrixDenoising,
  author  = {Leeb, William E.},
  title   = {Matrix Denoising for Weighted Loss Functions and Heterogeneous Signals},
  journal = {SIAM Journal on Mathematics of Data Science},
  year    = {2021},
  volume  = {3},
  number  = {3},
  pages   = {987--1012},
  doi     = {10.1137/20M1319577}
}

\end{document}



\def\spacingset#1{\renewcommand{\baselinestretch}%
{#1}\small\normalsize} \spacingset{1}
\if1\blind
{
  \title{\bf Supplement to ``Statistical Inference for Manifold Similarity and Alignability across Noisy High-Dimensional Datasets"}
  \author{Hongrui Chen\textsuperscript{1,2} and Rong Ma\textsuperscript{2}\thanks{Corresponding author: \href{mailto:rongma@hsph.harvard.edu}{rongma@hsph.harvard.edu}}}
  \date{\textsuperscript{1} Weiyang College, Tsinghua University, Beijing, China \\ 
        \textsuperscript{2} Department of Biostatistics, Harvard T.H. Chan School of Public Health, Boston, MA, USA}
          \date{}  
  \maketitle
  \tableofcontents
} \fi

\if0\blind
{
  \bigskip
  \bigskip
  \bigskip
  \begin{center}
    {\LARGE\bf Supplement to ``Statistical Inference for  Similarity and Alignability between Noisy High-Dimensional Datasets"}
    \date{}

  \end{center}
  \medskip
} \fi

\abstract{This Supplementary Material provides the comprehensive mathematical foundation, rigorous technical proofs, and extended empirical evidence supporting the Manifold Spectrometrics (MS) framework. We begin by giving a systematic outline of our theoretical approach, clarifying the core proof strategies and the formal construction of the proposed MS estimators. Building upon this primary roadmap, the document provides an in-depth treatment of critical methodological components: the convergence properties of the robust Newton iteration for outlier map inversion, the mean-squared consistency of the Potts-smoothed noise variance estimator, and a precise characterization of the supercriticality and separation conditions for signal spikes under heterogeneous noise environments. Furthermore, we establish the asymptotic normality results by deriving the joint central limit theorems (CLTs) for outlier eigenvalues and the normalized manifold spectral dissimilarity (nMSD), complemented by consistency proofs for the dissimilarity measure in general nonlinear settings via kernel-based manifold learning. To validate the practical utility and robustness of our inferential procedures, the supplement concludes with expanded simulation results and an empirical investigation into block-structured variance patterns observed in real-world single-cell RNA-seq datasets. Taken together with the main manuscript, this document ensures a self-contained and rigorous exposition of the proposed statistical methodology.}

\tableofcontents

\section{Technical Details of the Proposed Methods}

 \label{tech}

In this section, we present the theoretical foundations of our main results, focusing on the  derivation of the  estimators and the corresponding asymptotic variance estimation. Detailed proofs of the main theorems are provided in the following sections. All the analysis is conducted under \MainAsmpNoise, \MainAsmpSignal\ and \MainAsmpSignalStrength.

In Section~\ref{subsec:41}, we first establish delocalization properties of the signal representation and derive the outlier eigenvalue CLT conditional on the signal matrix. 
In Section~\ref{subsec:42}, under the Haar distribution assumption, we decouple the rank-$r$ secular equations and recover the signal strengths in closed form. 
Sections~\ref{subsec:variance} and \ref{Deltamethod} derive the variance decomposition of the outlier eigenvalues into noise- and signal-related components, connected through the outlier map $\theta$, and uses this result to obtain the asymptotic distribution of $\widehat{d}^{\,2}$ and the normalized vector $\widehat{\Pi}$. 
Finally, Section~\ref{subsec:45} constructs practical plug-in
approximations for the residual moments and covariance components used
in implementation.

\subsection{Conditional eigenvalue CLT for sample covariances}
\label{subsec:41}

\begin{lemma}[Delocalization]
\label{lem:B-deloc}

Let \(\SigmaS=\Sigma+\frac{S_\samp S_\samp^\top}{N}\) and let \(\psi_k\) be the \(k\)-th eigenvector of \  \( \SigmaS\) associated with an outlier \(\pssig_k\).
Set $B=\frac{S_\samp^\top\psi_k}{\sqrt{N}}\in\R^N$.
Under \MainAsmpNoise, \MainAsmpSignal and \MainAsmpSignalStrength, with probability $1-o(1)$, we have $\|B\|_\infty \;\lesssim \,\frac{\|D_\samp\|_{\mathrm{op}}\sqrt{n\log N}}{\sqrt N}$.
\end{lemma}
\vspace{-0.5em}
This lemma could be proved with the delocalization of $V_\samp$ in \MainPropmatspec and the low rank property of $S_\samp$ in \Mainprolowrank. Based on this lemma we could derive the asymptotic normality of $\lambda_k$ using the recent random matrix theory results \citep{LinPanZhaoZhou2024}. 

Under \MainAsmpNoise,\MainAsmpSignal and \MainAsmpSignalStrength, the deterministic outlier map and its derivative are
\vspace{-0.5em}
\begin{equation}\label{eq:theta-def}
\theta(\xi) = \xi + \frac{1}{N} \operatorname{tr} \biggl( \frac{\xi \Sigma}{\xi - \Sigma} \biggr), \quad
\theta'(\xi) = 1 - \frac{1}{N} \operatorname{tr} \biggl( \frac{\Sigma^2}{(\xi - \Sigma)^2} \biggr).
\end{equation}
Fix $r$ and let $\pssig_1>\cdots>\pssig_r$ be separated outliers of
$\SigmaS=\Sigma+\frac{S_\samp S_\samp^\top}{N}$ with associated unit eigenvectors
$\psi_1,\dots,\psi_r$.  For $k\le r\le n$ write
\(
a_k := \psi_k^\top \Sigma \psi_k, \quad
b_k
:=
\psi_k^\top S_{\mathrm{samp}}S_{\mathrm{samp}}^\top\psi_k/N
=
\pssig_k-a_k,
\quad
\theta'_k:=\theta'(\pssig_k).
\)
and denote by $\kappa_3:=\E[(X_{i\mu})^3]$ and
$\kappa_4:=\E[(X_{i\mu})^4]-3$ the third and fourth cumulants
of the standardized noise entries.  For vectors $u,v\in\R^p$ set
$M_{2,2}(u,v):=\sum_{i=1}^p u_i^2 v_i^2$.
\begin{lemma}[Conditional CLT for sample covariance outlier eigenvalues \citep{LinPanZhaoZhou2024}]
\label{lem:clt-lambda-con}
Under \MainAsmpNoise, \MainAsmpSignal and \MainAsmpSignalStrength and conditional on $S_\samp$, we have
\begin{equation}\label{eq:cond-clt}
\sqrt{N}\,\bigl(\lambda_k-\theta(\pssig_k)\bigr)
\;=\; \Phi_k \;+\; \Theta_k \;+\; L_k \;+\; o_{\mathbb P}(1),
\end{equation}
where $\E[\Phi_k\mid S_\samp]=\E[\Theta_k\mid S_\samp]=0$, and $L_k\equiv 0$ if $\kappa_3=0$.
Moreover, for each fixed $k$, conditional on $S_\samp$ the term $\Phi_k$ is
asymptotically normal as $N\to\infty$.

The conditional second moments admit the block decomposition
\vspace{-0.5em}
\begin{align}
\Var(\Phi_k\mid S_\samp)
&= V^{(G)}_{kk}
\;+\;\kappa_4\,(\theta'_k)^2\,
   M_{2,2}\!\bigl(\Sigma^{1/2}\psi_k,\Sigma^{1/2}\psi_k\bigr)
\;+\;o(1),
\label{eq:phi-variance}\\
\Cov(\Phi_k,\Theta_k\mid S_\samp)
&= W_{kk} \;+\; o(1),
\qquad
W_{kk}=\kappa_3\,(\theta'_k)^2\,a_k\,b_k,
\label{eq:phi-theta-cov}\\
\Var(\Theta_k\mid S_\samp)
&= 4\,(\theta'_k)^2\,a_k\,b_k,
\label{eq:theta-variance}
\end{align}
with the Gaussian block
\begin{equation}\label{eq:VG-kk}
V^{(G)}_{kk}
\;=\;
2\,(\theta'_k)^{2}a_k^{2}
\;+\;2\,{\pssig_k}^{2}\theta'_k
\;-\;2\,{\pssig_k}^{2}(\theta'_k)^{2}.
\end{equation}
Consequently,
\begin{equation}\label{eq:var-lambda-conditional}
\Var(\lambda_k\mid S_\samp) 
= \frac{1}{N} \Bigl( V^{(G)}_{kk} + 2\,W_{kk} + 4\,(\theta'_k)^2\,a_k\,b_k 
 + \kappa_4\,(\theta'_k)^2\,M_{2,2}\!\bigl(\Sigma^{1/2}\psi_k,\Sigma^{1/2}\psi_k\bigr) \Bigr) + o\!\Bigl(\tfrac{1}{N}\Bigr).
\end{equation}
In particular, if the noise is symmetric ($\kappa_3=0$) then $W_{kk}=0$; if it is
Gaussian ($\kappa_3=\kappa_4=0$) then only $V^{(G)}_{kk}$ and
$4(\theta'_k)^2 a_k b_k$ remain.
\end{lemma}

This lemma shows that the gap between the sample eigenvalue $\lambda$ and its corresponding limit $\theta$ could be further divided into three parts.  
To eventually derive the asymptotic normality of the nMSD estimator, we need the right-hand side of (\ref{eq:cond-clt}) to be asymptotically normal. The first term $\Phi$ has been proved to be asymptotically normal and the third part is zero under \MainAsmpNoise. The second term, however, is not asymptotically normal in general. Under our latent manifold model, we establish a delocalization property of $\frac{S_\samp^\top\psi_k}{\sqrt N}$ using Lemma~\ref{lem:B-deloc} and \MainPropmatspec, which ensures the asymptotic normality of $\Theta_k$.

\subsection{Solving secular equations}
\label{subsec:42}
Write
$
  \frac{1}{N} S_{\samp} S_{\samp}^\top = U_r D_r^2 U_r^\top,
$
where $D_r = \diag(d_1,\dots,d_r)$ collects the unknown signal strengths and 
$U_r \in \mathbb{R}^{p \times r}$ is Haar distributed on 
$\{U \in \mathbb{R}^{p \times r} : U^\top U = I_r\}$, independently of $\Sigma$.
Set $G_{U_r}(\xi) := U_r^\top A(\xi) U_r \in \mathbb{R}^{r \times r}$ and define the residual $R(\xi) := G_{U_r}(\xi) - g(\xi) I_r$.

Using the outlier map \eqref{eq:theta-def}, we can invert the sample spikes $\lambda_k$ to obtain estimates of the conditional population spikes $\pssig_k$. Our ultimate target, however, is the signal strengths $d_k^2$. These are linked to $\pssig_k$ through
\(
\Sigma + \frac{S_\samp S_\samp^\top}{N}= \SigmaS,
\)
because $d_k^2$ are the eigenvalues of $\frac{S_\samp S_\samp^\top}{N}$, whereas $\pssig_k$ are the eigenvalues of $\SigmaS$. For any $\pssig_j\notin\spec(\Sigma)$, the matrix determinant lemma (Sylvester’s identity) gives:
\[
\begin{aligned}
    &\det\bigl(\Sigma + U_r D_r^2 U_r^\top - \xi_j^{\mathrm{S}} I_p\bigr) 
    = \det(\Sigma - \xi_j^{\mathrm{S}} I_p) \det\bigl(I_r + D_r^2 U_r^\top(\Sigma - \xi_j^{\mathrm{S}} I_p)^{-1} U_r\bigr) \\
    &\quad = \det(\Sigma - \xi_j^{\mathrm{S}} I_p) \det\bigl(I_r + D_r^2 G_{U_r}(\xi_j^{\mathrm{S}})\bigr).
\end{aligned}
\]
Hence $\pssig_j$ is an eigenvalue of $\Sigma+U_r D_r^2 U_r^\top$ if and only if
\(
\det\!\big(I_r + D_r^2\,G_{U_r}(\pssig_j)\big)=0, \quad j = 1,2,\dots,r.
\)
Our next lemma shows that the matrix $G_{U_r}(\pssig_j)$ can be well approximated by a diagonal 
matrix $g(\pssig_j) I_r$ independent of $U_r$, thereby providing a means to simplify the secular 
system and derive its closed-form solution.

\begin{lemma}[Concentration of $G_{U_r}(\xi)$]\label{lem:isotropy}
Under \MainAsmpNoise, \MainAsmpSignal and \MainAsmpSignalStrength, with $U_r\sim\mathrm{Haar}$ on $\{U \in \mathbb{R}^{p \times r} : U^\top U = I_r\}$ and independent of $\Sigma$, we have $\E\!\big[G_{U_r}(\xi)\big]\;=\;g(\xi)\,I_r.$
Moreover, as $(N,p)\to\infty$ with fixed $r$ and $\min_{1\le i\le p}\,|\sigma_i-\xi|\ge c>0$, we have
\begin{equation}\label{eq:Rop-bound}
\big\|R(\xi)\big\|_{\op}
\;
=\;
O_{\P}\!\Big(\frac{\sqrt{s_2(\xi)}}{p}\sqrt r
\;+\;\|A(\xi)\|_{\op}\sqrt{\frac{r}{p}}\Big)
\;
 =\;o_{\P}(1),
\end{equation}
where $\|\cdot\|_{\op}$ denotes the operator norm.

\end{lemma}
\begin{proposition}[Exact solution under Haar prior]\label{prop:decouple}
Assume the setting of Lemma~\ref{lem:isotropy} with fixed $r$.
Set $g_j:=g(\pssig_j)$ and assume $|g_j|\ge g_{\min}\ge c_0>0$. Let $D_r=\diag(d_1,\dots,d_r)$ solve the secular system \MainEqRankRSecular. If, in addition, the numbers $\{-1/g_j\}_{j=1}^r$ are pairwise separated by at least $\Delta>0$, then there exists a unique permutation $\pi$ of $\{1,\dots,r\}$ such that
\begin{equation}\label{dg}
d_{\pi(j)}^2 \;=\; -\frac{1}{g_j} \;+\; o_{\mathbb P}(1),\qquad j=1,\dots,r.
\end{equation}
\end{proposition}

For each $j=1,\dots,r$, write $G_j := G_{U_r}(\pssig_j)$ and $R_j := R(\pssig_j)$.
When $G_j$ is replaced by its conditional expectation $g_j I_r$, 
\MainEqRankRSecular reduce to 
$\prod_{\ell=1}^r (1 + g_j d_\ell^2) = 0$ for each $j$. 
The unique multiset solution is therefore 
$\{d_\ell^2\}_{\ell=1}^r = \{-1 / g_j\}_{j=1}^r$. 
Lemma~\ref{lem:isotropy} extends this mean-field result to a high-probability statement, establishing that 
$\max_j \|R_j\|_{\mathrm{op}} = o_{\mathbb{P}}(1)$, 
which in turn yields $d_{\pi(j)}^2 = -1 / g_j + o_{\mathbb{P}}(1)$ for some permutation $\pi$, as formalized in Proposition~\ref{prop:decouple}. 

\subsection{Unconditional eigenvalue CLT for sample covariances} \label{subsec:variance}

A key feature of our data-generating model is that the randomness of the final estimators arises from both the noise and the signal-generating processes. In this section, we first derive the population-level contribution to the asymptotic covariance of the sample eigenvalues  $\lambda_k$ attributable to signal randomness, and then establish the unconditional eigenvalue CLT for the sample covariance.

\subsubsection{Delineating uncertainty due to noiseless samples}
\label{subsec:signal-resampling-eig}

In this subsection we isolate the variability that comes from repeatedly
generating noiseless signal samples $S_1,\dots,S_N$, while keeping the noise
covariance $\Sigma$ fixed.
We treat the empirical signal covariance
$\widehat M = N^{-1}\sum_{i=1}^N S_i S_i^\top$ as a random
perturbation of its population counterpart $M$, and use first–order
eigenvalue perturbation theory to describe how a simple spike of
$\widetilde\Sigma=\Sigma+M$ changes when $M$ is replaced by $\widehat M$.
This shows that, to leading order, the shift of each outlier eigenvalue is
governed by a quadratic form of the perturbation taken along its population
eigenvector.
We then derive explicit formulas for the covariance matrix of these quadratic
forms, which will later serve as the ``noiseless–sample'' component in the
overall asymptotic covariance of the sample spikes.

Throughout, remember that \(S_i\in\R^p\) denote a random column with
\(\E\|S_i\|^4<\infty\), population second moment is represented by \(M=\E[S_iS_i^\top]\), and the sample second moment can be calculated by
\(\widehat M=\frac1N\sum_{i=1}^N S_i S_i^\top\). We define the centered fluctuation
\(\Delta_M:=\widehat M-M\).

Let \(\widetilde\Sigma = \Sigma + M\) and suppose it has a simple outlier eigenvalue \(\xi_k\) with associated unit eigenvector \(\popsi_k\).
Denote by \(\lambda_k(A)\) the \(k\)-th eigenvalue of a symmetric matrix \(A\), and define the eigengap at \(\xi_k\) by
$  \gamma_k := \min\{\xi_{k-1} - \xi_k,\ \xi_k - \xi_{k+1}\} > 0.$

\begin{lemma}[First-order sensitivity of a simple spike]
\label{lem:hadamard}
For any symmetric perturbation \(E\in\R^{p\times p}\), we have $
\lambda_k(\widetilde\Sigma+E)-\lambda_k(\widetilde\Sigma)
\;=\;\popsi_k^\top E\,\popsi_k\;+\;R_k(E),$ and $
|R_k(E)|\ \le\ C\,\frac{\|E\|_{\op}^2}{\gamma_k},$
for a universal constant \(C>0\). In particular, if \(\|E\|_{\op}=o(1)\),
then the spike shift satisfies
\(
\lambda_k(\widetilde\Sigma+E)-\lambda_k(\widetilde\Sigma)
=\popsi_k^\top E\,\popsi_k+o(\|E\|_{\op}).
\)
\end{lemma}

This is the classical Hadamard first--variation formula for simple eigenvalues,
together with a second--order remainder bounded by the eigengap. It follows from
analytic perturbation theory or by a resolvent expansion plus the
Davis--Kahan control of eigenprojections; the gap \(\gamma_k\) ensures
differentiability of \(\lambda_k\) at \(\widetilde\Sigma\).
Applying Lemma~\ref{lem:hadamard} with \(E=\Delta_M\) shows that, conditional on
\((\Sigma,M)\) so that \(\popsi_k\) is deterministic and since
\(\|\Delta_M\|_{\op}=O_p(N^{-1/2})\), we have
\begin{equation}\label{eq:delta-sigma-first-order}
\delta\xi_k := \lambda_k(\Sigma + \widehat M) - \lambda_k(\Sigma + M) 
= \tilde{\psi}_k^\top \Delta_M \tilde{\psi}_k + o_p(N^{-1/2}).
\end{equation}
Thus, to first order the spike responds only to the quadratic form of the
perturbation \(\Delta_M\) along its own eigenvector \(\popsi_k\).
The next lemma computes the covariance of those quadratic forms. It is stated
for general fixed evaluation directions and then specialized to \(v_k=\popsi_k\).

\begin{lemma}[Quadratic-form covariance]
\label{lem:signal-resampling-general-new}
Let \(v_1,\dots,v_m\in\R^p\) be fixed unit vectors, independent of the sample
\(\{S_i\}_{i=1}^N\). Define
$
\delta_k\;:=\;v_k^\top\Delta_M\,v_k
\;=\;\frac1N\sum_{i=1}^N\left\{(v_k^\top S_i)^2-\E[(v_k^\top S_i)^2]\right\}.$
Then \(\E[\delta_k]=0\) and the covariance matrix
\(\Gamma^{(\mathrm{sig})}=(\Gamma^{(\mathrm{sig})}_{kj})\) of
\((\delta_1,\dots,\delta_m)\) satisfies
\vspace{-0.5em}
\begin{equation}\label{eq:Gamma-general-raw-new}
\Gamma^{(\mathrm{sig})}_{kj} = \Cov(\delta_k, \delta_j) 
= \frac{1}{N} \biggl\{ \mathbb{E} \bigl[ (v_k^\top S)^2 (v_j^\top S)^2 \bigr] 
- \mathbb{E} \bigl[ (v_k^\top S)^2 \bigr] \mathbb{E} \bigl[ (v_j^\top S)^2 \bigr] \biggr\}.
\end{equation}
Writing the fourth raw--moment tensor \(T^{(4)}:=\E[S^{\otimes4}]\) and
\(b_k:=v_k^\top M v_k\), this equals
\begin{equation}\label{eq:Gamma-general-tensor-new}
\Gamma^{(\mathrm{sig})}_{kj}
=\frac1N\big\langle v_k\!\otimes v_k\!\otimes v_j\!\otimes v_j,\ T^{(4)}\big\rangle
-\frac1N\,b_k\,b_j .
\end{equation}
If, in addition,
\(\Kc:=\mathrm{Cum}_4(S)\) denotes the fourth central cumulant tensor, then
\begin{equation}\label{eq:Gamma-cumulant-new}
\quad
\Gamma^{(\mathrm{sig})}_{kj}
=\frac1N\Big\{\,2\,(v_k^\top M\,v_j)^2+\Kc[v_k,v_k,v_j,v_j]\Big\}.
\quad
\end{equation}

\end{lemma}
Combining \eqref{eq:delta-sigma-first-order} with Lemma~\ref{lem:signal-resampling-general-new}
specialized to \(v_k=\popsi_k\),  conditional on \((\Sigma,M)\), we have
\[
\sqrt N\,\delta\xi_k = \sqrt N\,\popsi_k^\top\Delta_M\popsi_k+o_p(1), \quad
\Cov\big(\popsi_k^\top\Delta_M\popsi_k,\ \popsi_j^\top\Delta_M\popsi_j\big) = \Gamma^{(\mathrm{sig})}_{kj}+o(N^{-1}).
\]
with \(\Gamma^{(\mathrm{sig})}\) given by
\eqref{eq:Gamma-cumulant-new}, upon substituting \(v_k=\popsi_k\).
To first order, the spike \(\xi_k\) depends only on the quadratic form of \(\Delta_M\) along its own eigenvector \(\popsi_k\). Evaluating along any other direction \(v \neq \popsi_k\) would correspond to the perturbation of an unrelated linear functional and would not capture the spike’s first-order behavior. In practice, since neither the true signal eigenvector nor the conditional population eigenvector is observable, we estimate the variance using plug-in sample eigenvectors. 

\subsubsection{Unconditional CLT and total asymptotic variance}
\label{lambda:variance}

Recall that $(\pssig_k,\psi_k)$ refer to the conditional population eigenpairs of $\SigmaS$, in contrast to the population pairs $(\xi_k,\popsi_k)$ used in the previous subsection.
Fix $r$ and let $\pssig_1>\cdots>\pssig_r$ be separated outliers of
$\SigmaS=\Sigma+\frac{S_\samp S_\samp^\top}{N}$ with unit eigenvectors
$\psi_1,\dots,\psi_r$.  For $k,j\le r$ write
$
A_{kj}:=\psi_k^\top\Sigma\psi_j,$ $
B_{kj}:=\psi_k^\top S_\samp S_\samp^\top\psi_j/N,$ $
a_k:=A_{kk},$ $b_k:=B_{kk}=\pssig_k-a_k,
$ and $\theta'_k:=\theta'(\pssig_k)$.
Let the standardized noise entries have third and fourth cumulants
$\kappa_3:=\E[(X_{i\mu})^3]$ and
$\kappa_4:=\E[(X_{i\mu})^4]-3$, and set
$M_{2,2}(u,v):=\sum_{i=1}^p u_i^2 v_i^2$ for $u,v\in\R^p$.
The next proposition concerns the unconditional CLT for the eigenvalues of the sample covariance matrix.

\begin{proposition}
[Joint CLT of outlier eigenvalues]\label{pro:joint-clt-uncond}
Under \MainAsmpNoise, \MainAsmpSignal and \MainAsmpSignalStrength, fix $r$ and let $\xi_1>\cdots>\xi_r$ be separated outliers of $\widetilde\Sigma=\Sigma+M$ with unit eigenvectors $\popsi_1,\dots,\popsi_r$.Let
$\boldsymbol\pssig=(\pssig_1,\dots,\pssig_r)^\top$
denote the conditional spikes, set
$\theta'_k:=\theta'(\pssig_k)$
for the conditional covariance block, and set
$\dot\theta_k:=\theta'(\xi_k)$ and
$\dot\Theta:=\operatorname{diag}(\dot\theta_1,\dots,\dot\theta_r)$
for the signal-resampling block.

Assume a fixed outlier gap and the moment conditions stated above.

Then we have
$
\sqrt N\bigl(\boldsymbol\lambda-\theta(\boldsymbol\xi)\bigr)\ \xrightarrow{d}\ \mathcal N_r\!\bigl(0,\ V_*\bigr),$
where the limiting covariance $V_*$ admits the decomposition
\[
V_*
=
\mathbb E\!\bigl[V^{\rm cond}\bigr]
+
N\dot\Theta\,
\Gamma^{(\mathrm{sig})}\,
\dot\Theta^\top.
\]

Here $\nabla\theta(\boldsymbol\pssig)=\diag(\theta'_1,\dots,\theta'_r)$, and
\begin{enumerate}
\item for the {conditional noise block} $V^{\rm cond}=(V^{\rm cond}_{kj})$, which is the conditional limiting covariance given $S_\samp$, we have:
\[
V^{\rm cond}_{kj} 
= V^{(G)}_{kj} + 2\,W_{kj} + 4\,\theta'_k\theta'_j\,A_{kj}\,B_{kj} 
+ \kappa_4\,\theta'_k\theta'_j\,M_{2,2}\!\bigl(\Sigma^{1/2}\psi_k,\Sigma^{1/2}\psi_j\bigr),
\]
with $A_{kj}=\psi_k^\top\Sigma\psi_j$ and $B_{kj}=\psi_k^\top S_\samp S_\samp^\top\psi_j/N$, we also have:
\[
V^{(G)}_{kk} = 2(\theta'_k)^2 a_k^2+2{\pssig_k}^{\,2}\theta'_k-2{\pssig_k}^{\,2}(\theta'_k)^2, 
W_{kj}=\kappa_3\,\theta'_k\theta'_j\,A_{kj}B_{kj},
\]
and $a_k=A_{kk}$.
\item for the signal-related variance block
$N\dot\Theta\Gamma^{(\mathrm{sig})}\dot\Theta^\top$, we have:
\[
\begin{aligned}
&\Gamma^{(\mathrm{sig})}_{kj} 
= \Cov\Bigl( \popsi_k^\top(\widehat M-M)\popsi_k, 
\quad \popsi_j^\top(\widehat M-M)\popsi_j \Bigr) \\
&= \frac{1}{N} \Bigl\{ \E\bigl[ (\popsi_k^\top S)^2 (\popsi_j^\top S)^2 \bigr] 
 - \E[(\popsi_k^\top S)^2] \E[(\popsi_j^\top S)^2] \Bigr\} \quad + o\bigl(\tfrac{1}{N}\bigr).
\end{aligned}
\]
so that the total variance satisfies:
\[
\operatorname{Var}(\boldsymbol\lambda)
=
\frac1N\mathbb E[V^{\rm cond}]
+
\dot\Theta\,
\Gamma^{(\mathrm{sig})}\,
\dot\Theta^\top
+
o\!\left(\frac1N\right).
\]

\end{enumerate}
\end{proposition}

\subsection{CLT of the nMSD estimator}
\label{Deltamethod}

Recall that by \MainAlgPi, we have
\begin{equation}
\label{eq:d2-hat}
\quad
\widehat d_j^{\,2}
= -\Biggl(\frac{1}{p}\sum_{i=1}^p\frac{1}{\Sigma_{ii}-\widehat{\pssig}_j}\Biggr)^{-1},
\qquad j=1,\dots,r.
\quad
\end{equation}
Under the setting of Proposition~\ref{pro:joint-clt-uncond}, suppose that
$\theta'(\xi_j)\neq0$ for $j=1,\dots,r$. Applying the multivariate delta method to the composition
$\lambda_j\mapsto\theta^{-1}(\lambda_j)\mapsto-1/g(\theta^{-1}(\lambda_j))$
gives

\[
\sqrt{N}\,\bigl(\widehat d^{\,2}-d^{2}\bigr)\ \xrightarrow{d}\
\mathcal N\!\bigl(0,\,\Sigma_{d^2}\bigr),
\]
where $\Sigma_{d^2}=\Gamma\,V_*\,\Gamma $,  $d^2=(d_1^2,\dots,d_r^2)^\top$ and $\Gamma=\diag(\Gamma_1,\dots,\Gamma_r)$ with
\begin{equation}\label{Gamma_j}
\Gamma_j
=
\left[
\mathrm D\bigl(\varphi\circ\theta^{-1}\bigr)
\bigl(\theta(\boldsymbol\xi)\bigr)
\right]_{jj}
=
\underbrace{
\frac{s_2(\xi_j)}{p\,g(\xi_j)^2}
}_{\varphi'(\xi_j)}
\cdot
\underbrace{
\frac{1}{\theta'(\xi_j)}
}_{(\theta^{-1})'(\theta(\xi_j))}.
\end{equation}
where $\mathrm D$ denotes the Jacobian matrix,
and $V_*$ is defined in subsection~\ref{lambda:variance}.
To see this, note that for each $j$, 
\begin{equation}
\label{eq:vard2}
\operatorname{Var}\!\bigl(\widehat d_j^{\,2}\bigr)
=
\left(
\frac{s_2(\xi_j)}
{p\,g(\xi_j)^2}
\right)^{\!2}
\frac{\operatorname{Var}(\lambda_j)}
{\bigl(\theta'(\xi_j)\bigr)^2}
+
o\!\left(\frac{1}{N}\right).
\end{equation}
Let  $s=\mathbf 1^\top d^2$, and $\Pi=d^2/s\in\Delta^{r-1}$.
With
\begin{equation}
\label{eq:varJ}
J^\Pi_{kj}=\frac{\partial \Pi_k}{\partial d_j^2}
=\frac{\delta_{kj}\,s-d_k^2}{s^2},
\end{equation}
and Writing
$\widehat{\boldsymbol\xi}:=\theta^{-1}(\boldsymbol\lambda)$
and evaluating the relevant Jacobians at the deterministic population spikes
$\boldsymbol\xi$, we have
$\mathrm D\theta^{-1}\bigl(\theta(\boldsymbol\xi)\bigr)
=
\operatorname{diag}\!\bigl(
1/\theta'(\xi_1),\dots,1/\theta'(\xi_r)
\bigr)$
and
$\mathrm D\varphi(\boldsymbol\xi)
=
\operatorname{diag}\!\bigl(
\varphi'(\xi_1),\dots,\varphi'(\xi_r)
\bigr)$,
where
$\varphi'(s)=s_2(s)/(p\,g(s)^2)$.

On the other hand, under the setting of Proposition~\ref{pro:joint-clt-uncond}, let $\widehat\Pi=\widehat d^{\,2}/(\mathbf 1^\top\widehat d^{\,2})$, $\Pi=d^{2}/(\mathbf 1^\top d^{2})$ and $s=\mathbf 1^\top d^2>0$. Then we have
\begin{equation}
\label{eq:Pi-CLT}
\quad
\sqrt{N}\,\bigl(\widehat\Pi-\Pi\bigr)\ \xrightarrow{d}\
\mathcal N\!\bigl(0,\ \Sigma_\Pi\bigr),
\end{equation}
on $T_\Pi\Delta^{r-1}:=\{u\in\mathbb R^r:\mathbf 1^\top u=0\}$, where
\begin{equation}
\label{eq:SigmaPi-mult}
\Sigma_\Pi
=J^\Pi\,\Gamma\,V_*\,\Gamma\,(J^\Pi)^\top
\;,
\Sigma_\Pi\mathbf 1=\mathbf 0,
\quad
\operatorname{rank}(\Sigma_\Pi)\le r-1.
\end{equation}
With the above CLT for $\widehat\Pi$, we can then immediately obtain the CLT for the nMSD estimator.

\subsection{Practical plug-in covariance approximation}
\label{subsec:45}
We now explain how to turn the theoretical variance expressions into
practical estimators.  Throughout this subsection we suppress the dataset
index $k$ and work with a fixed sample covariance $Q=YY^\top/N$.

Let $(\lambda_j,\widehat\psi_j)$, $j=1,\dots,r$, be the outlier eigenpairs
of $Q$, where $\widehat\psi_j$ is the $j$th sample eigenvector. In practice
we use $\widehat\psi_j$ to approximate both the population outlier
eigenvector $\popsi_j$ of $\widetilde\Sigma=\Sigma+M$ and the conditional
outlier eigenvector $\psi_j$ of $\Sigma+\widehat M$. This approximation does not affect consistency of the principal variance profile estimator itself.

\paragraph*{(1) Conditional (noise) block}
Define the plug-in analogues of
$A_{kj}$, $B_{kj}$, and $M_{2,2}$ by
\[
\widehat A_{kj} := \widehat\psi_k^\top\,\widehat\Sigma\,\widehat\psi_j, \quad
\widehat B_{kj} := \widehat\psi_k^\top\,\widehat M\,\widehat\psi_j, \quad
\widehat M_{2,2}(k,j) := \sum_{a=1}^p \widehat\psi_{k,a}^{\,2}\,\widehat\psi_{j,a}^{\,2}\,\widehat\Sigma_{aa}^2.
\]
where $\widehat M$ is the empirical signal covariance (for example, the
rank-$r$ fit) and $\widehat\psi_{k,a}$ denotes the $a$th coordinate of
$\widehat\psi_k$.

To estimate the noise cumulants, we form residuals by projecting $Y$
onto the orthogonal complement of the spike subspace spanned by
$\widehat\psi_1,\dots,\widehat\psi_r$:
\[
  \widehat P := \sum_{k=1}^r \widehat\psi_k \widehat\psi_k^\top,
  \quad
  \widehat R_i := Y_i - \widehat P Y_i,\quad i=1,\dots,N,
\]
and compute coordinatewise empirical moments
\[
\hat m_{2,a} = \frac1N\sum_{i=1}^N \widehat R_{a,i}^2, \quad
\hat m_{3,a} = \frac1N\sum_{i=1}^N \widehat R_{a,i}^3, \quad
\hat m_{4,a} = \frac1N\sum_{i=1}^N \widehat R_{a,i}^4.
\]
The third and fourth standardized cumulants are then estimated by
\[
  \widehat\kappa_3
  := \frac{1}{p}\sum_{a=1}^p \frac{\hat m_{3,a}}{\hat m_{2,a}^{3/2}},
  \qquad
  \widehat\kappa_4
  := \frac{1}{p}\sum_{a=1}^p \Bigl(\frac{\hat m_{4,a}}{\hat m_{2,a}^{2}}-3\Bigr).
\]

Plugging these quantities into the conditional covariance formula yields
\[
  \widehat V^{\rm cond}_{kj}
  = \widehat V^{(G)}_{kj}
    + \widehat\kappa_4\,\hat\theta'_k\hat\theta'_j\,\widehat M_{2,2}(k,j)
    + 2\,\widehat W_{kj}
    + 4\,\hat\theta'_k\hat\theta'_j\,\widehat A_{kj}\,\widehat B_{kj},
\]
where
\(
\widehat V^{(G)}_{kk} = 2(\hat\theta'_k)^2 \widehat A_{kk}^{\,2} + 2\,\widehat\pssig_k^{\,2}\,\hat\theta'_k  - 2\,\widehat\pssig_k^{\,2}\,(\hat\theta'_k)^{2}, \quad
\widehat W_{kj} = \widehat\kappa_3\,\hat\theta'_k\hat\theta'_j\,\widehat A_{kj}\,\widehat B_{kj}.
\)
\paragraph*{(2) Signal-related covariance block}
The signal-related covariance $\Gamma^{(\mathrm{sig})}$ was defined
theoretically in Lemma~\ref{lem:signal-resampling-general-new} in terms
of the population eigenvectors $\popsi_k$ and the true signal columns
$S_i$.
In practice we approximate $\popsi_k$ by $\widehat\psi_k$ and estimate the signal–sampling covariance via
\[
  \widehat\Gamma^{(\mathrm{sig})}_{kj}
  := \frac{2}{N}\,\widehat B_{kj}^{\,2}
     + \frac{1}{N}\Bigl(\widehat{\mathcal K}_Y - \widehat{\mathcal K}_\varepsilon\Bigr)
        [\widehat\psi_k,\widehat\psi_k,\widehat\psi_j,\widehat\psi_j],
\]
where $\widehat{\mathcal K}_Y$ and $\widehat{\mathcal K}_\varepsilon$ denote the empirical fourth–order cumulant tensors of $Y$ and of the estimated noise, respectively.

\paragraph*{(3) Variance of spikes}
Combining the conditional and signal–sampling parts, the plug-in
covariance matrix of the sample spikes is
\[
  \widehat{\Var}(\boldsymbol\lambda)
  = \frac{1}{N}\,\widehat V^{\rm cond}
    + \widehat{\nabla\theta}\,
      \widehat\Gamma^{(\mathrm{sig})}\,
      \widehat{\nabla\theta}^\top,\quad
  \widehat{\nabla\theta} := \diag(\hat\theta'_1,\dots,\hat\theta'_r),
\]
so in particular
\(
  \widehat{\Var}(\lambda_j)
  = \frac{1}{N}\,\widehat V^{\rm cond}_{jj}
    + \hat\theta_j'^{\,2}\,\widehat\Gamma^{(\mathrm{sig})}_{jj}.
\)
Under symmetric Gaussian noise ($\kappa_3=\kappa_4=0$) this simplifies to
\(
  \widehat V^{\rm cond}_{kj}
  = \widehat V^{(G)}_{kj}
    + 4\,\hat\theta'_k\hat\theta'_j\,\widehat A_{kj}\,\widehat B_{kj},
\)
and if moreover $\widehat A_{kj}\approx\widehat B_{kj}\approx 0$ for
$k\neq j$, then $\widehat V^{\rm cond}$ is approximately diagonal.

\paragraph*{(4) Variance of signal strengths and spectral profiles}
Define the empirical counterparts of $g$ and $s_2$ by replacing
$\Sigma$ with $\widehat\Sigma$ in their definitions:
\[
  \widehat g(s)
  := \frac{1}{p}\sum_{i=1}^p \frac{1}{\widehat\sigma_i - s},
  \qquad
  \widehat s_2(s)
  := \sum_{i=1}^p \frac{1}{(\widehat\sigma_i - s)^2}.
\]
The delta–method derivative for the mapping
$\lambda_j \mapsto \widehat d_j^{\,2}$ is estimated by
\[
  \widehat\Gamma_j
  := \frac{\widehat s_2(\widehat\pssig_j)}
          {p\,\widehat g(\widehat\pssig_j)^2\,\hat\theta'_j},
  \qquad
  \widehat\Gamma := \diag(\widehat\Gamma_1,\dots,\widehat\Gamma_r).
\]
Let $\widehat V_*$ denote the practical plug-in approximation to the
oracle asymptotic covariance $V_*$ in
Proposition~\ref{pro:joint-clt-uncond}, constructed using the empirical
quantities described above.
Then the covariance matrices of $\widehat d^{\,2}$ and of the normalized
profile $\widehat\Pi$ are estimated by
\begin{equation}\label{eq:SigmaPi-plugin}
  \widehat\Sigma_{d^2}
  = \widehat\Gamma\,\widehat V_*\,\widehat\Gamma,\quad
  \widehat\Sigma_{\Pi}
  = J^\Pi(\widehat d^{\,2})\,
    \widehat\Gamma\,\widehat V_*\,
    \widehat\Gamma\,J^\Pi(\widehat d^{\,2})^\top, \quad
  \widehat V_\Pi = \frac{\widehat\Sigma_\Pi}{N},
\end{equation}
The matrices in \eqref{eq:SigmaPi-plugin} are the practical plug-in
covariance approximations used in implementation. They are obtained by
replacing the population noise levels, population eigenvectors, and
other oracle quantities by their empirical counterparts. The Wald
theorems stated in the main paper concern the oracle covariance, whereas the accuracy of
this feasible substitution is evaluated numerically in the simulation
studies.

\section{Two examples for nMSD}
Here we provide two examples demonstrating some basic properties of $d_r(\rho_1,\rho_2)$.

    \begin{example}[Scale and rotation invariance]
Let $(\M,\mathbb P)$ be a compact manifold with a probability measure $\mathbb P$
(bounded density w.r.t.\ volume). Let $\iota_1:\M\to\R^3$ be an embedding and set
\[
\rho_1:=\mathbb P\circ\iota_1^{-1},\qquad
M_1:=\E_{S\sim\rho_1}[SS^\top],
\]
whose eigenvalues are \(\mathcal D^{(1)}=(4,1,1)\) (this is just a choice of
coordinates/eigenbasis for $M_1$). For any orthogonal \(R\in\mathbb O(3)\) and any
scalar \(c>0\), define another embedding
\[
\iota_2:=c\,R\,\iota_1,\qquad
\rho_2:=\mathbb P\circ\iota_2^{-1}.
\]
Then
\[
\begin{aligned}
&M_2 = \E_{S\sim\rho_2}[SS^\top] \\
&= \E_{x\sim\mathbb P}\big[(cR\iota_1(x))(cR\iota_1(x))^\top\big] \\
&= c^2\,R\,M_1\,R^\top,
\end{aligned}
\]
so the eigenvalues satisfy \(\mathcal D^{(2)}=c^2\,\mathcal D^{(1)}\).
Taking \(c=\sqrt{2}\) gives \(\mathcal D^{(2)}=(8,2,2)\).

Normalize the top-3 principal variances:
\[
\begin{aligned}
&\Pi_3(\rho_1)=\frac{(4,1,1)}{6}=\Big(\tfrac23,\tfrac16,\tfrac16\Big),\\
&\Pi_3(\rho_2)=\frac{(8,2,2)}{12}=\Big(\tfrac23,\tfrac16,\tfrac16\Big),
\end{aligned}
\]
hence
\(
d_3(\rho_1,\rho_2)=\bigl\|\Pi_3(\rho_1)-\Pi_3(\rho_2)\bigr\|_2=0.
\)
This shows $d_r$ depends on $\rho=\mathbb P\circ\iota^{-1}$ only through the
\emph{relative} principal variances and is invariant under ambient rotations and global
rescaling of the embedding.
\end{example}

\begin{example}
Consider uniform distributions supported on \(2\)-dimensional manifolds isometrically embedded in \(\R^3\) by $\iota \equiv\textup{id}$: i) a \emph{unit sphere} \(\mathcal M_1=\mathbb{S}^2\) with a uniform distribution $\mathbb{P}_1$, whose pushforward measure $\rho_1$ has a population covariance containing three equal eigenvalues $\mathcal D^{(1)} \;=\;(\tfrac13,\;\tfrac13,\;\tfrac13)$; and ii) an \emph{ellipsoid} \(\mathcal M_2\) elongated from $\mathcal{M}_1$ along the first axis, with a uniform distribution $\mathbb{P}_2$, whose pushforward measure $\rho_2$ has a population covariance with eigenvalues $\mathcal D^{(2)} \;=\;(\tfrac43,\;\tfrac13,\;\tfrac13).$
Then we compute
\[
\begin{aligned}
&\Pi_3(\rho_1)
=\frac{(\tfrac13,\tfrac13,\tfrac13)}{\tfrac13+\tfrac13+\tfrac13}
=(\,\tfrac13,\;\tfrac13,\;\tfrac13\,),\\
&\Pi_3(\rho_2)
=\frac{(\tfrac43,\tfrac13,\tfrac13)}{\tfrac43+\tfrac13+\tfrac13}
=(\,\tfrac{4}{6},\;\tfrac{1}{6},\;\tfrac{1}{6}\,),
\end{aligned}
\]
so that
\[
\begin{aligned}
&d_3(\rho_1,\rho_2)
=\bigl\|\Pi_3(\rho_1)-\Pi_3(\rho_2)\bigr\|_2 \\
&=\sqrt{\Bigl(\tfrac{4}{6}-\tfrac{1}{3}\Bigr)^{2}
      +2\Bigl(\tfrac{1}{6}-\tfrac{1}{3}\Bigr)^{2}}
\approx 0.408.
\end{aligned}
\]
This value quantifies the extent to which the ellipsoid concentrates variance along its leading principal direction relative to the uniform, isotropic variance distribution of the sphere.
\end{example}

\begin{figure}
  \centering
  \includegraphics[width=0.4\linewidth]{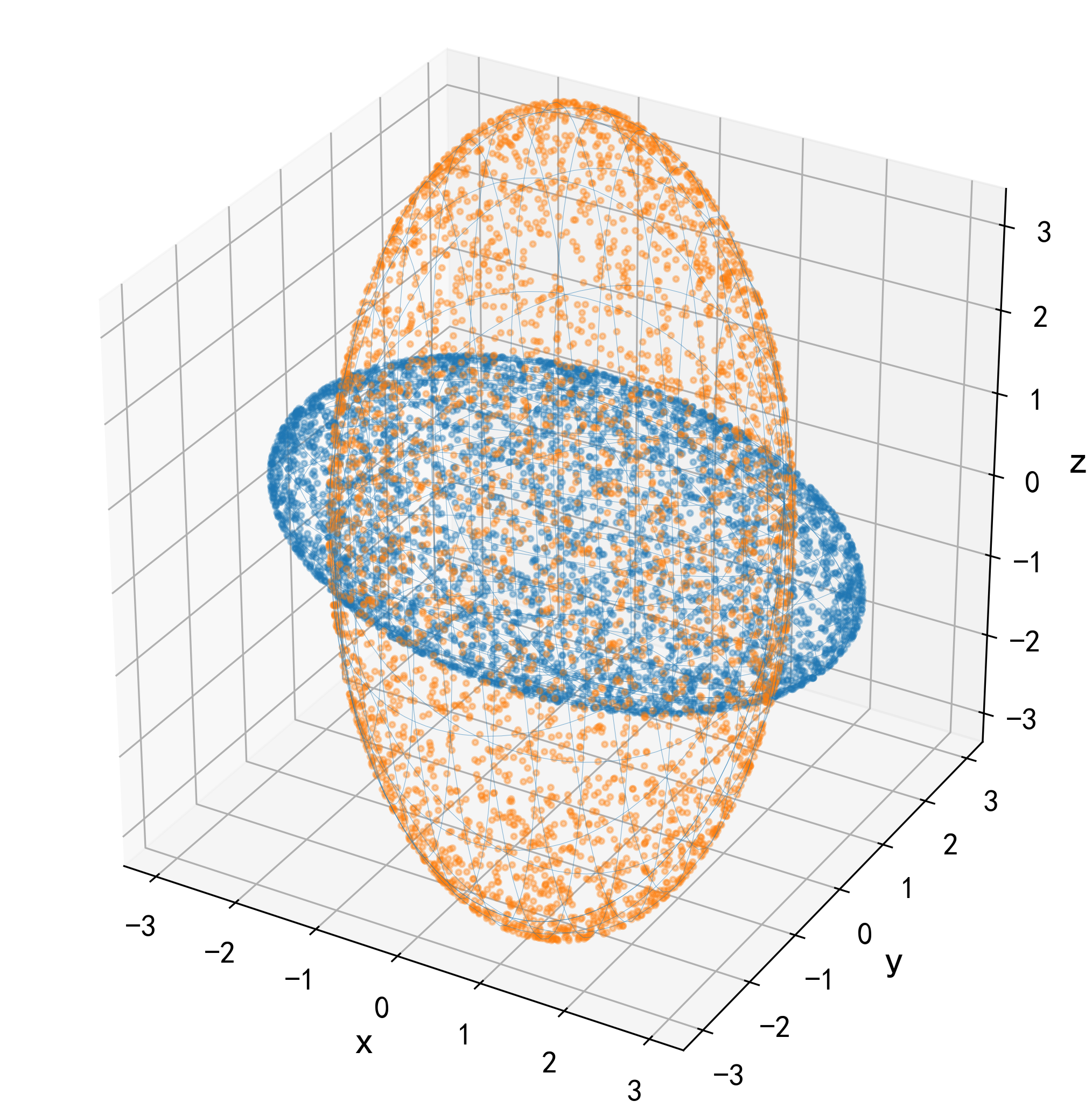}
  \caption{Rotation/scale invariance of the normalized principal–variance vector.
  Blue/orange point clouds are samples from $\rho_1$ and its rotated–rescaled version $\rho_2$ in $\mathbb{R}^3$;
  wireframes are the corresponding ellipsoids. Although $M_2=c^2\,R M_1 R^\top$, the normalized vectors
  $\Pi_3(\rho_1)$ and $\Pi_3(\rho_2)$ coincide, hence $d_3(\rho_1,\rho_2)=0$.}
  \label{fig:rot-scale}
\end{figure}

\section{Robust Newton iteration algorithm for solving the outlier inversion map}

Step~2 of \MainAlgPi \ requires, for each sample spike
$\lambda_{j,k}$, solving the scalar equation
\[
  \theta\!\bigl(s;\widehat\Sigma_k,\varphi_k\bigr)
  \;=\; \lambda_{j,k},
  \qquad s>\max_i \widehat\sigma_{i,k},
\]
where $\widehat\Sigma_k=\diag(\widehat\sigma_{1,k},\dots,\widehat\sigma_{p,k})$
is the estimated diagonal noise covariance for dataset $k$ and
$\varphi_k=p/N_k$ is the aspect ratio.
For fixed $(\widehat\Sigma_k,\varphi_k)$, the outlier map
$s\mapsto\theta(s;\widehat\Sigma_k,\varphi_k)$ is continuously differentiable
and strictly monotone on $(\max_i\widehat\sigma_{i,k},\infty)$, and the image
contains each observed outlier $\lambda_{j,k}$.
Hence the equation admits a unique solution
$s^\star_{j,k}>\max_i\widehat\sigma_{i,k}$.

In practice we compute $s^\star_{j,k}$ by a damped Newton iteration applied to
the function
\[
  f_{j,k}(s)
  := \theta\!\bigl(s;\widehat\Sigma_k,\varphi_k\bigr)-\lambda_{j,k},
  \qquad
  f'_{j,k}(s) = \theta'\!\bigl(s;\widehat\Sigma_k,\varphi_k\bigr),
\]
restricted to the domain $s>\max_i\widehat\sigma_{i,k}$.
A convenient implementation is as follows.

\medskip
\noindent\emph{Initialization.}
Set
\[
  s^{(0)}_{j,k}
  := \max_i\widehat\sigma_{i,k}+\varepsilon,
\]
where $\varepsilon>0$ is a small numerical offset (e.g.\ $\varepsilon =
10^{-3}\,\mathrm{median}_i(\widehat\sigma_{i,k})$).
This guarantees $s^{(0)}_{j,k}$ lies in the admissible domain and that
$\theta(s^{(0)}_{j,k};\widehat\Sigma_k,\varphi_k)$ is well defined.

\medskip
\noindent\emph{Newton step with damping.}
Given an iterate $s^{(t)}_{j,k}$, form the full Newton update
\[
  \tilde s^{(t+1)}_{j,k}
  := s^{(t)}_{j,k}
     - \frac{f_{j,k}\!\bigl(s^{(t)}_{j,k}\bigr)}
            {f'_{j,k}\!\bigl(s^{(t)}_{j,k}\bigr)}
  = s^{(t)}_{j,k}
     - \frac{\theta\!\bigl(s^{(t)}_{j,k};\widehat\Sigma_k,\varphi_k\bigr)
              -\lambda_{j,k}}
            {\theta'\!\bigl(s^{(t)}_{j,k};\widehat\Sigma_k,\varphi_k\bigr)}.
\]
To ensure robustness, we apply a simple backtracking line search:
starting from step size $\eta=1$, repeatedly replace
$\eta\leftarrow\eta/2$ until the damped proposal
\[
  s^{(t+1)}_{j,k}
  := s^{(t)}_{j,k}
     - \eta\,
       \frac{\theta\!\bigl(s^{(t)}_{j,k};\widehat\Sigma_k,\varphi_k\bigr)
                -\lambda_{j,k}}
            {\theta'\!\bigl(s^{(t)}_{j,k};\widehat\Sigma_k,\varphi_k\bigr)}
\]
satisfies both $s^{(t+1)}_{j,k}>\max_i\widehat\sigma_{i,k}$ and
\[
  \bigl|f_{j,k}(s^{(t+1)}_{j,k})\bigr|
  \le \bigl|f_{j,k}(s^{(t)}_{j,k})\bigr|.
\]
If at some iteration the derivative
$\bigl|\theta'\!\bigl(s^{(t)}_{j,k};\widehat\Sigma_k,\varphi_k\bigr)\bigr|$
falls below a small threshold, we also reduce $\eta$ to avoid overly large
steps.  In all our experiments, at most a handful of backtracking reductions
is needed.

\medskip
\noindent\emph{Stopping rule and output.}
The iteration is terminated when either
\[
  \bigl|\theta\!\bigl(s^{(t)}_{j,k};\widehat\Sigma_k,\varphi_k\bigr)
        -\lambda_{j,k}\bigr|
  \le \texttt{tol}
  \qquad\text{or}\qquad
  t\ge t_{\max},
\]
for a prescribed numerical tolerance \texttt{tol} (e.g.\ $10^{-8}$) and
maximum number of iterations $t_{\max}$ (e.g.\ $20$).
The final iterate is then recorded as
\[
  \widehat{\pssig}_{j,k}
  := s^{(t)}_{j,k}.
\]

\medskip
Under the regularity condition $\theta'(\pssig_{j,k})\neq0$ used in the main
text, the above damped Newton scheme converges locally at a linear rate, and
once the iterates enter a small neighborhood of $\pssig_{j,k}$ the convergence
is quadratic.  In particular, with the initialization described above and
moderate $t_{\max}$, the residual
$\bigl|\theta(\widehat{\pssig}_{j,k};\widehat\Sigma_k,\varphi_k)-\lambda_{j,k}\bigr|$
is negligible at numerical precision.

\section{Consistency of the Potts-smoothed noise estimator}

In this section we fix one dataset and drop the index $k$ in the
signal-plus-noise model~\MainEqSignalPlusNoise. Thus
\[
  Y \;=\; \Sigma^{1/2}X + S_{\samp}\in\R^{p\times N},
\]
where $S_{\samp}=[S_1,\dots,S_N]$ is the signal component and
$\Sigma=\diag(\sigma_1,\dots,\sigma_p)$ is the diagonal,
block-heteroskedastic noise covariance from \MainAsmpNoise.
The random matrix $X=[x_{i\mu}]$ has independent entries satisfying
\MainAsmpNoise(ii).
We assume the signal $S_i$ follows the manifold model in
\MainAsmpSignal.  Denote
\[
\begin{aligned}
  & M \;:=\; \frac1N\E[S_{\samp}S_{\samp}^\top]=\; \E[S_i S_i^\top], \\
  & C \;:=\; \frac1N\E[YY^\top] \;=\; M+\Sigma,
\end{aligned}
\]
and let $\sigma=(\sigma_1,\dots,\sigma_p)^\top$ be the diagonal of $\Sigma$.

Under \MainAsmpSignal, the embedded manifold
$\iota(\M)\subset\R^p$ has Euclidean dimension
$\dim(\iota(\M))=n\le p$.
By the orthogonal change of coordinates described below
\MainAsmpSignal, the signal covariance $M$ has rank $n$; we
call $n$ the \emph{signal rank}.
In what follows we assume $n$ is known and we always set the working
rank to be $r=n$.

Let $Q:=YY^\top/N$ be the sample covariance of $Y$, and let $\widehat C$ be
the best rank-$n$ approximation of $Q$ in Frobenius norm, i.e.
\[
  \widehat C \;=\; \widehat U_n \widehat\Lambda_n \widehat U_n^\top,
\]
where $\widehat \Lambda_n=\diag(\widehat \lambda_1,\dots,\widehat \lambda_n)$ contains the
largest $n$ eigenvalues of $Q$ and $\widehat U_n$ the corresponding
eigenvectors.
We define the \emph{factor analysis(FA) residual diagonal}
\[
  \widehat\sigma_i\ :=\ Q_{ii}-\widehat C_{ii},\qquad
  \widehat\sigma:=\mathrm{diag}(Q-\widehat C),
\]
and the Potts-regularized estimator
\[
\begin{aligned}
  \widehat\sigma^{\mathrm{pwc}}
  \in\arg\min_{\tau\in\mathbb R^p}
  \Big\{
    \sum_{i=1}^p(\widehat\sigma_i-\tau_i)^2 
    \; \\
    +\;\beta\,\#\{1\le i<p:\tau_{i+1}\neq\tau_i\}
  \Big\},
\end{aligned}
\]
with tuning parameter $\beta>0$ and jump-count functional
\[
  J(\tau):=\#\{1\le i<p:\tau_{i+1}\neq\tau_i\}.
\]

\subsection{Diagonal concentration}

\begin{lemma}[Diagonal concentration]\label{lem:diag}
Under \MainAsmpNoise \ and \MainAsmpSignal,
there exist constants $c,C>0$ such that for all $t>0$,
\[
\begin{aligned}
  \Pr\!\left(
    \max_{1\le i\le p}
    \left|
      \frac1N\sum_{s=1}^N\big(Y_{is}^2-\mathbb EY_{is}^2\big)
    \right|>t
  \right)\\
  \;\le\;
  2p\exp\!\Big(
    -cN\min\{t^2,\ t\}
  \Big).
\end{aligned}
\]
In particular, taking $t\asymp \sqrt{(\log p)/N}$,
\[
\begin{aligned}
  \max_{1\le i\le p}\big|Q_{ii}-C_{ii}\big|
    =O_\mathbb P\!\Big(\sqrt{\tfrac{\log p}{N}}\Big),\\
  \frac1p\sum_{i=1}^p(Q_{ii}-C_{ii})^2
    =O_\mathbb P\!\Big(\tfrac{\log p}{N}\Big).
\end{aligned}
\]
\end{lemma}
\begin{proof}
\MainAsmpNoise(ii) yields a uniform $(2+\delta)$-moment bound
for the entries of $X$. \MainAsmpSignal\ places the signal
$S_i$ on a compact manifold, so each coordinate of $S_i$ is almost surely
bounded. Since $Y_{is}=\sigma_i^{1/2}x_{is}+S_{is}$ is a sum of a
sub-Gaussian variable and a bounded variable, there exists a constant
$K>0$ such that $\|Y_{is}\|_{\psi_2}\le K$ uniformly in $i,s$.

Fix $i$. Then
$Z_{is}:=Y_{is}^2-\mathbb EY_{is}^2$
is sub-exponential with
$\|Z_{is}\|_{\psi_1}\le CK^2$
for some $C>0$.
Bernstein's inequality for i.i.d.\ sub-exponential variables gives
\[
\begin{aligned}
\Pr\!\left(
\left|
\frac1N\sum_{s=1}^N Z_{is}
\right|>t
\right)
&\le
2\exp\!\left(
-cN\min\{t^2/K^4,\ t/K^2\}
\right)\\
&\le
2\exp\!\left(
-c'N\min\{t^2,\ t\}
\right)
\end{aligned}
\]
for suitable constants $c,c'>0$ depending only on $K$.
A union bound over $i=1,\dots,p$ yields the stated probability bound.
Taking $t\asymp\sqrt{(\log p)/N}$ gives the two consequences.
\end{proof}

\subsection{Rank-\texorpdfstring{$n$}{n} approximation error}
\begin{lemma}[Rank-$n$ approximation error]\label{lem:rankn}
Under \MainAsmpNoise\ and \MainAsmpSignal, let $\widehat C$ be the best
rank-$n$ approximation to $Q$, with fixed $n$. Then
\[
\frac{1}{p}
\sum_{i=1}^p
(\widehat C_{ii}-M_{ii})^2
=
O_{\mathbb P}\!\left(\frac{n}{N}\right).
\]
\end{lemma}

\begin{proof}
Since $\widehat C$ has rank at most $n$ and
$\|\widehat C\|_{\mathrm{op}}\le\|Q\|_{\mathrm{op}}$, we have
\[
\|\widehat C\|_F^2
\le
n\|\widehat C\|_{\mathrm{op}}^2
\le
n\|Q\|_{\mathrm{op}}^2
=
O_{\mathbb P}(n),
\]
where $\|Q\|_{\mathrm{op}}=O_{\mathbb P}(1)$ follows from the
sub-Gaussian assumption and the proportional high-dimensional regime.

Likewise, since $M$ has rank at most $n$ and bounded operator norm,
\[
\|M\|_F^2
\le
n\|M\|_{\mathrm{op}}^2
=
O(n).
\]
Therefore,
\[
\begin{aligned}
\frac{1}{p}
\sum_{i=1}^p
(\widehat C_{ii}-M_{ii})^2
&\le
\frac{1}{p}
\|\widehat C-M\|_F^2\\
&\le
\frac{2}{p}
\left(
\|\widehat C\|_F^2+\|M\|_F^2
\right)\\
&=
O_{\mathbb P}\!\left(\frac{n}{p}\right)
=
O_{\mathbb P}\!\left(\frac{n}{N}\right),
\end{aligned}
\]
because $p/N$ is bounded away from zero and infinity.
\end{proof}

\subsection{Low-rank fitting and Potts denoising}

\begin{proposition}[FA residual diagonal]\label{prop:FA-residual}
Under \MainAsmpNoise\ and \MainAsmpSignal, let
$\widehat\sigma=\mathrm{diag}(Q-\widehat C)$ with fixed $r=n$.
Then
\[
  \frac1p\sum_{i=1}^p(\widehat\sigma_i-\sigma_i)^2
  =O_\mathbb P\!\Big(\frac{\log p}{N}\Big).
\]
\end{proposition}

\begin{proof}
Recall that $C_{ii}=M_{ii}+\sigma_i$, hence
\[
  \widehat\sigma_i-\sigma_i
  = (Q_{ii}-\widehat C_{ii}) - \sigma_i
  = (Q_{ii}-C_{ii}) - (\widehat C_{ii}-M_{ii}).
\]
Using $(a-b)^2\le 2a^2+2b^2$ and averaging over $i$,
\[
  \frac1p\sum_{i=1}^p(\widehat\sigma_i-\sigma_i)^2
  \le
  \frac{2}{p}\sum_{i=1}^p(Q_{ii}-C_{ii})^2
  +\frac{2}{p}\sum_{i=1}^p(\widehat C_{ii}-M_{ii})^2.
\]
Lemma~\ref{lem:diag} bounds the first term by
$O_\mathbb P((\log p)/N)$, and Lemma~\ref{lem:rankn} shows that the second
term is $O_\mathbb P(n/N)$.  Since $n$ is fixed while $\log p\to\infty$,
we have $n/N=o((\log p)/N)$, so the whole expression is
$O_\mathbb P((\log p)/N)$.
\end{proof}
\subsection{Proof of \MainThmNoiseConsistency}

\begin{proof}
Write $\widehat\sigma=\sigma+\varepsilon$, where
$\varepsilon_i=\widehat\sigma_i-\sigma_i$.
By Proposition~\ref{prop:FA-residual},
\[
  \frac1p\sum_{i=1}^p\varepsilon_i^2
  =O_\mathbb P\!\Big(\frac{\log p}{N}\Big).
\]

By optimality of $\widehat\sigma^{\mathrm{pwc}}$, we have
\[
  \sum_{i=1}^p(\widehat\sigma_i-\widehat\sigma^{\mathrm{pwc}}_i)^2
  +\beta J(\widehat\sigma^{\mathrm{pwc}})
  \;\le\;
  \sum_{i=1}^p(\widehat\sigma_i-\sigma_i)^2
  +\beta J(\sigma).
\]
Using $\widehat\sigma=\sigma+\varepsilon$ and setting
$\delta:=\widehat\sigma^{\mathrm{pwc}}-\sigma$, this becomes
\[
  \|\delta-\varepsilon\|_2^2 + \beta J(\widehat\sigma^{\mathrm{pwc}})
  \;\le\;
  \|\varepsilon\|_2^2 + \beta J(\sigma),
\]
hence
\[
  \|\delta-\varepsilon\|_2^2
  \;\le\;
  \|\varepsilon\|_2^2 + \beta J(\sigma).
\]
By the inequality $\|\delta\|_2^2\le 2\|\varepsilon\|_2^2
+2\|\delta-\varepsilon\|_2^2$, we obtain
\[
  \|\widehat\sigma^{\mathrm{pwc}}-\sigma\|_2^2
  =\|\delta\|_2^2
  \le
  2\|\varepsilon\|_2^2
  +2\|\delta-\varepsilon\|_2^2
  \le
  4\|\varepsilon\|_2^2
  +2\beta J(\sigma).
\]
\MainAsmpNoise(iv) states that $\sigma$ is block-constant
with $K_p$ blocks, so $J(\sigma)=K_p-1=o(p)$.
Dividing by $p$ and using $\beta\asymp(\log p)/N$,
\[
  \frac1p\|\widehat\sigma^{\mathrm{pwc}}-\sigma\|_2^2
  \le
  \frac{4}{p}\|\varepsilon\|_2^2
  +2\beta\frac{J(\sigma)}{p}
  =
  O_\mathbb P\!\Big(\frac{\log p}{N}\Big)
  +o\!\Big(\frac{\log p}{N}\Big),
\]
which yields the claimed bound.
\end{proof}

\subsection{Empirical advantages of the proposed variance estimators} 
\begin{figure}[t]
  \centering
  \includegraphics[width=0.65\linewidth]{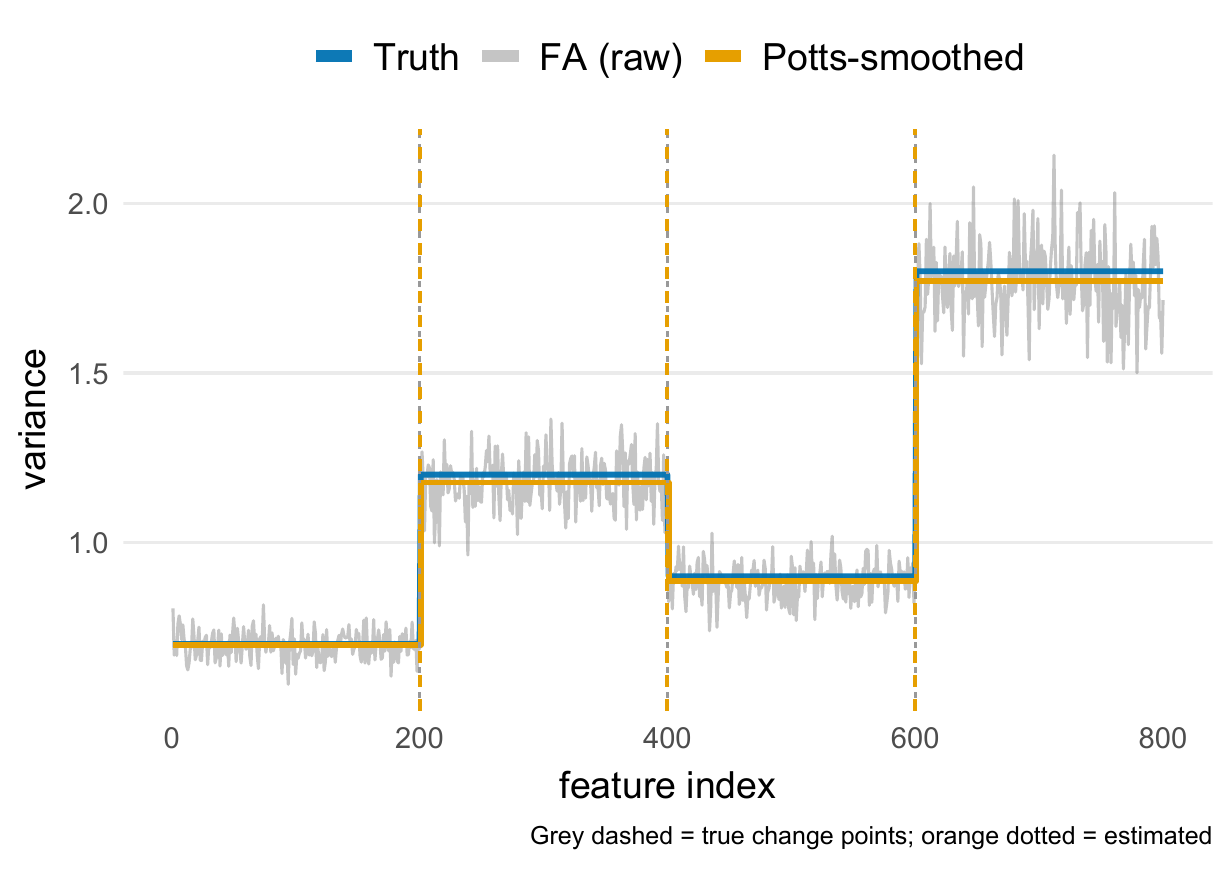}
  \caption{Significant improvement by the Potts segmentation step.
  Black: true noise variances; gray: initial variance estimators $\{\widehat\sigma_i\}$; blue: piecewise-smoothed estimators $\{\widehat\sigma_i^{\textup{pwc}}\}$.
  Vertical dashed lines mark the true change points. Our data-driven procedure selects $\beta=0.205$, leading to MSE $=0.0004$ for the piecewise-smoothed estimators, and MSE $=0.006$ for the initial variance estimators.}
  \label{fig:noise-example}
\end{figure}

Empirically, we find that the Potts segmentation step in our algorithm significantly improves the initial variance estimates $\widehat{\sigma}$. To demonstrate this, we generate a high-dimensional dataset containing low-rank signals and piecewise heteroskedastic noise following model (1) in our main text. The diagonal noise profile consists of four contiguous segments with distinct variances arranged along a fixed feature order. We apply Algorithm 1 and compare the performance of the initial estimators $\widehat{\sigma}$ with the piecewise-smoothed estimators $\widehat{\sigma}^{\textup{pwc}}$. As shown in Figure~\ref{fig:noise-example}, using a penalty of $\beta = 0.205$, the Potts-smoothed estimators achieve an MSE of $0.0004$ against the ground truth, representing a substantial improvement over the raw residual-based estimators (MSE = $0.006$).

\section{Further explanation about \MainAsmpSignalStrength}
\label{supp:supercritical}

In the main text, \MainAsmpSignalStrength \ requires the leading
$r$ eigenvalues of the conditional population covariance $\Sigma^S$ to be
“supercritical and mutually separated”. In this section we give a precise
formulation of these conditions.

\subsection{Noise bulk edge}

Recall that the noise covariance matrix is denoted by
\(\Sigma = \mathrm{diag}(\sigma_1,\ldots,\sigma_p)\) and that
\MainAsmpNoise \ implies the empirical distribution of the
diagonal entries of $\Sigma$ converges weakly to a probability measure
\(\nu\) on $[0,\infty)\).
Let $X\in\mathbb{R}^{p\times N}$ have i.i.d.\ standardized entries, and
consider the \emph{noise-only} sample covariance matrix
\[
Q^{\mathrm{noise}}
:= \frac{1}{N}\,\Sigma^{1/2} X X^\top \Sigma^{1/2}.
\]
Under \MainAsmpNoise, the empirical spectral distribution of
$Q^{\mathrm{noise}}$ converges almost surely to a deterministic probability
measure $\mu_{\mathrm{noise}}$ on $[0,\infty)$, whose support is a compact
interval. We denote by
\[
  \lambda_+^{\mathrm{noise}}
  := \sup\mathrm{supp}(\mu_{\mathrm{noise}})
\]
the right edge of the noise-only limiting spectrum and refer to it as the
(noise) bulk edge.

Equivalently, the bulk edge can be characterized in terms of the outlier map
$\theta$ in equation~\eqref{eq:theta-def} : let $m_{\mathrm{noise}}$ be the Stieltjes
transform of $\mu_{\mathrm{noise}}$, and define
\[
\begin{aligned}
  &\theta(s)
  \;=\;
  s \;+\; \varphi\,\int \frac{s\,t}{s-t}\,\nu(dt),\\
  &\theta'(s)
  \;=\;
  1 \;-\; \varphi\,\int \frac{t^2}{(s-t)^2}\,\nu(dt).
\end{aligned}
\]
Then the critical value $s_{\mathrm{crit}}$ solving $\theta'(s_{\mathrm{crit}})=0$
corresponds to the bulk edge via
\(
\lambda_+^{\mathrm{noise}} = \theta(s_{\mathrm{crit}}).
\)
We do not need the explicit form of $\mu_{\mathrm{noise}}$ in what follows;
only the existence of $\lambda_+^{\mathrm{noise}}$ and the local regularity
of $\theta$ near the spikes are used.

\subsection{Supercritical and mutually separated conditional spikes}

Let
\[
  \Sigma^S := \Sigma + \frac{1}{N} S_{\samp} S_{\samp}^\top
\]
be the conditional population covariance matrix given the noiseless signal
matrix $S_{\samp}=[S_1,\dots,S_N]$, and let
\[
  \pssig_1 \;\ge\; \pssig_2 \;\ge\; \cdots \;\ge\; \pssig_p
\]
denote its eigenvalues in nonincreasing order. We are interested in the
leading $r$ eigenvalues $\{\pssig_k\}_{k\le r}$, which correspond to the
signal spikes.

\begin{definition}[Supercritical and mutually separated spikes]
\label{def:supercritical-separated}
We say that the leading $r$ eigenvalues
$\{\pssig_k\}_{k\le r}$ of $\Sigma^S$ are \emph{supercritical and mutually
separated} if there exist constants $\delta>0$ and $\tau>0$, independent of
$(N,p)$, such that for all sufficiently large $(N,p)$ the following hold:
\begin{enumerate}
\item \textbf{Supercriticality above the noise bulk.}
Each spike lies strictly above the noise bulk edge with a fixed margin:
\[
  \pssig_k \;\ge\;
  \lambda_+^{\mathrm{noise}} + \delta,
  \qquad k=1,\dots,r.
\]
\item \textbf{Separation from non-spike eigenvalues.}
The remaining eigenvalues of $\Sigma^S$ stay within $\delta/2$ of the bulk:
\[
  \pssig_{r+1} \;\le\; \lambda_+^{\mathrm{noise}} + \frac{\delta}{2},
\]
where we adopt the convention $\pssig_{r+1}=0$ if $r=p$.
\item \textbf{Mutual separation of spikes.}
The spike eigenvalues are separated from each other by a fixed gap:
\[
  \min_{1\le k<\ell\le r} |\pssig_k-\pssig_\ell|
  \;\ge\; \tau.
\]
\end{enumerate}
\end{definition}

Condition (i) formalizes the “supercritical” regime for the spikes: each
$\pssig_k$ is located strictly to the right of the noise bulk edge, by at
least $\delta$. Condition (ii) ensures that there are no additional spikes
near the bulk edge. Condition (iii) guarantees that the spikes
$\{\pssig_k\}_{k\le r}$ are well separated from each other and from the
bulk, so that the outlier map $\theta$ is well defined, strictly monotone,
and admits a smooth inverse in a neighborhood of each $\pssig_k$.

These assumptions are directly analogous to the supercriticality threshold
and eigenvalue-gap conditions used in spiked covariance models with general
heteroskedastic noise; see, for example \citep{LinPanZhaoZhou2024}.
Under Definition~\ref{def:supercritical-separated}, the conditional sample
outlier eigenvalues $\{\lambda_k\}_{k\le r}$ converge to
$\{\theta(\pssig_k)\}_{k\le r}$, remain separated from the bulk spectrum,
and satisfy the conditional CLT and variance formulas stated in
Lemma \ref{lem:clt-lambda-con}.

\section{Additional numerical results}

\subsection{Additional figures from simulation studies} 

\begin{table*}[htbp]
\centering
\caption{Comparison of empirical quantiles and $\chi^2_2$ distribution.}
\label{tab:quantiles-null}
\tabcolsep=3pt 
\begin{tabular*}{\columnwidth}{@{\extracolsep{\fill}}lccccccc@{\extracolsep{\fill}}}
\toprule
$q$ & 0.01 & 0.10 & 0.25 & 0.50 & 0.75 & 0.90 & 0.99 \\
\midrule
$T_\Pi$   & 0.019 & 0.223 & 0.569 & 1.280 & 2.601 & 4.688  & 9.423 \\
$\chi^2_2$      & 0.020 & 0.211 & 0.575 & 1.386 & 2.773 & 4.605 & 9.210 \\
Diff  & $-0.001$ & 0.012 & $-0.007$ & $-0.106$ & $-0.106$ & 0.083 & 0.212 \\
\bottomrule
\end{tabular*}
\end{table*}
\begin{table*}[htbp]
\centering
\caption{Power and noncentrality (\(r=3\), \(N_1=N_2=3000\)).}
\label{tab:power}
\begin{tabular}{lrrrr}
\toprule
$c$ & $\|\Delta\Pi\|$ & $\lambda$ (theory) & Power (theory) & Power (emp.) \\
\midrule
1.00 & 0.00000 & \phantom{0}0.0000 & 0.0500 & 0.0417 \\
1.05 & 0.01488 & \phantom{0}1.4371 & 0.1725 & 0.1833 \\
1.10 & 0.02912 & \phantom{0}5.5051 & 0.5453 & 0.5233 \\
1.20 & 0.05586 &           20.3010 & 0.9864 & 0.9833 \\
1.30 & 0.08050 &           42.3648 & 1.0000 & 1.0000 \\
1.50 & 0.12439 &          102.7337 & 1.0000 & 1.0000 \\
\bottomrule
\end{tabular}
\end{table*}

Table~S1 assesses the null calibration of the feasible statistic $T_\Pi$ by comparing its empirical quantiles with the oracle $\chi^2_2$ reference. Table~S2 compares the empirical power of the feasible procedure with an oracle-based noncentral chi-square approximation across varying signal separations. These comparisons provide numerical assessments of feasible calibration and power.

\begin{figure}[h!]
    \centering
    \includegraphics[height=0.6\linewidth,width=1\linewidth]{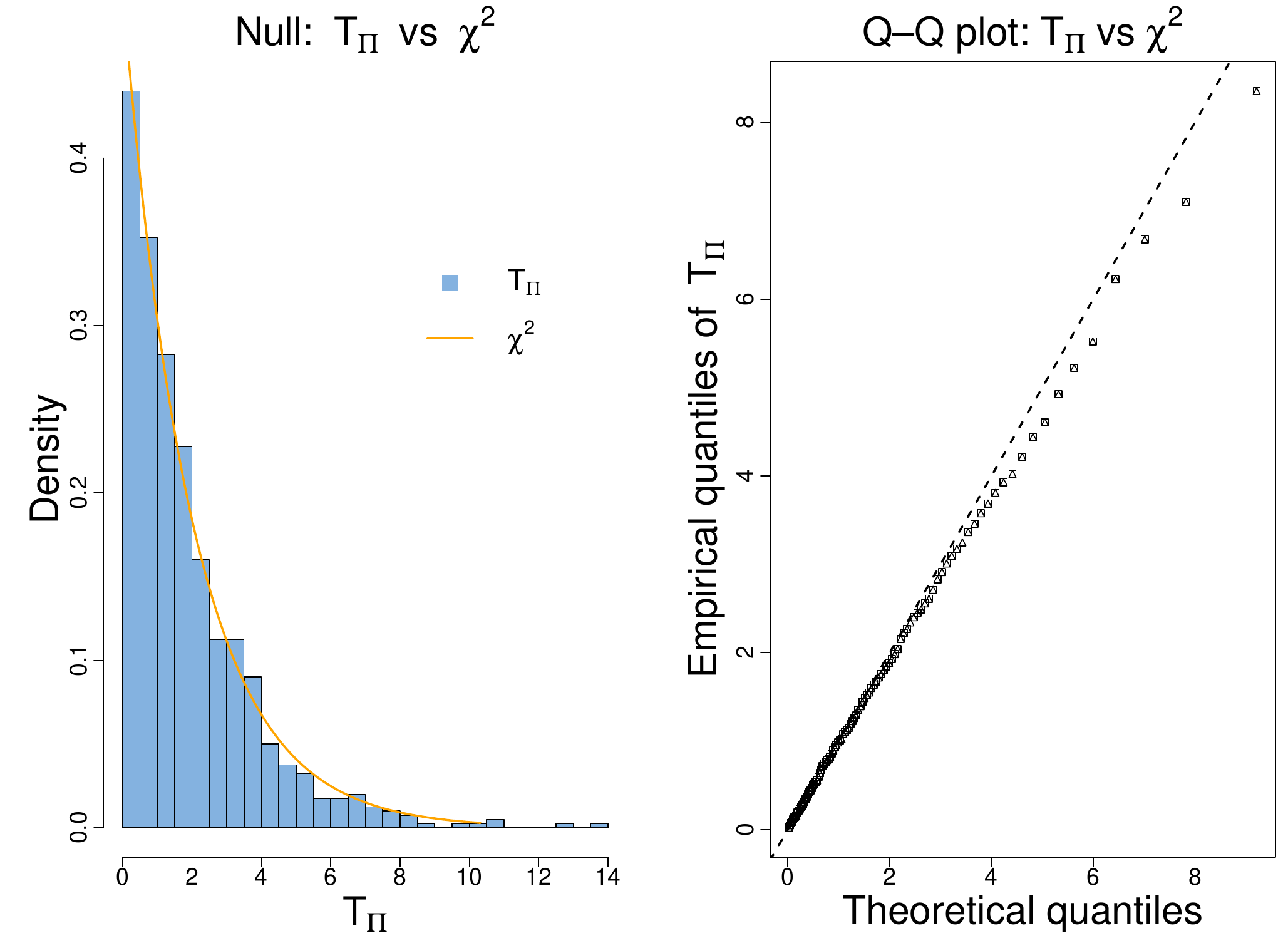}
      \caption{\textbf{Simulation: calibration and power.}
  Left: null calibration of \(T_\Pi\). Right: Goodness of fit Q-Q plot.}
    \label{fig:exp1}
\end{figure}

\begin{figure}[h!]
    \centering
    \includegraphics[height=0.6\linewidth,width=0.8\linewidth]{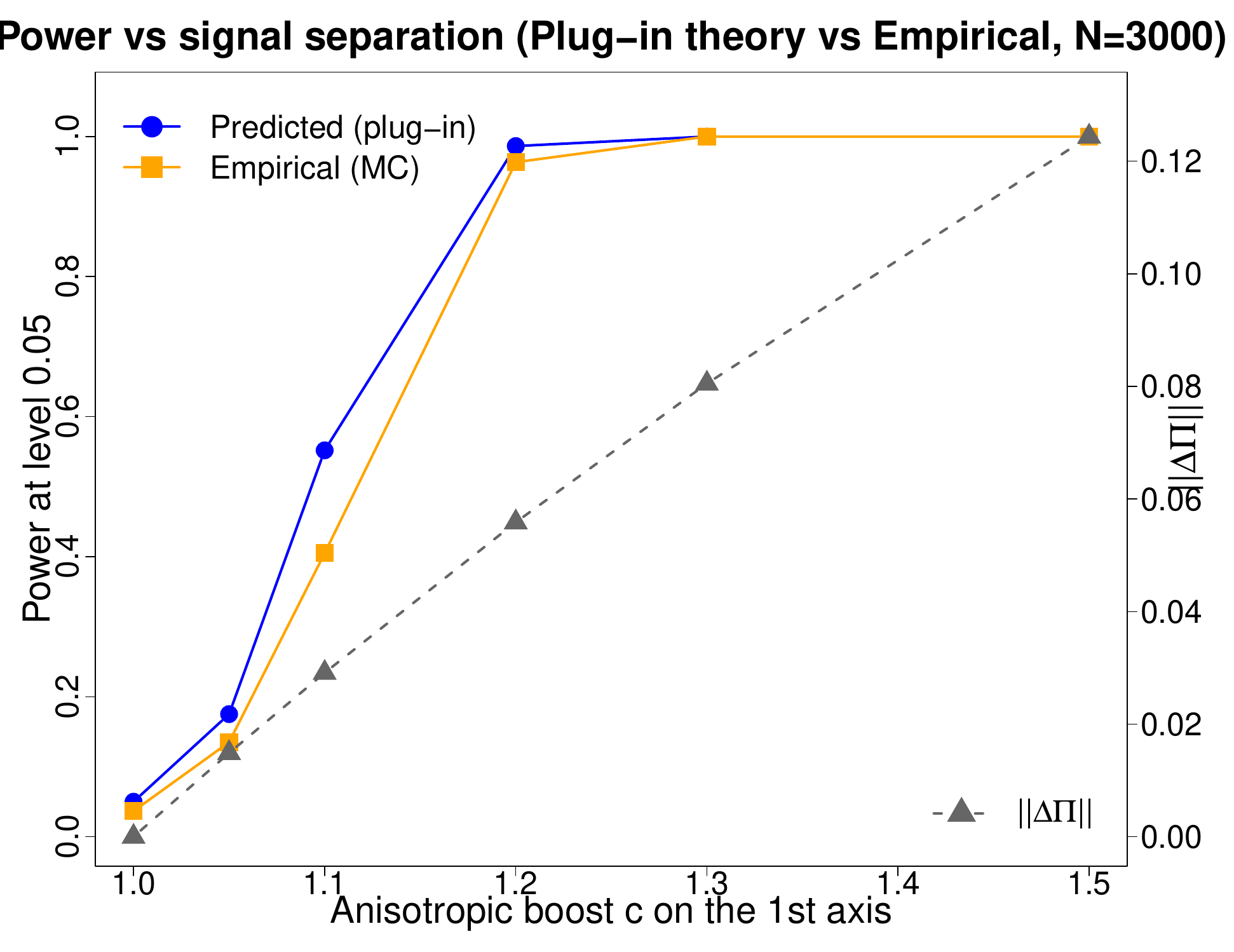}
    \caption{Power vs.\ anisotropic separation (\(r=3\), \(N_1=N_2=3000\)). Empirical power (orange) closely tracks the noncentral \(\chi^2\) prediction (blue) as \(\|\Delta\Pi\|\) increases with \(c\).}
     \label{fig:exp2}
\end{figure}

\begin{figure}[htbp]
  \centering

  \begin{minipage}{0.48\linewidth}
    \centering
    \includegraphics[width=\linewidth]{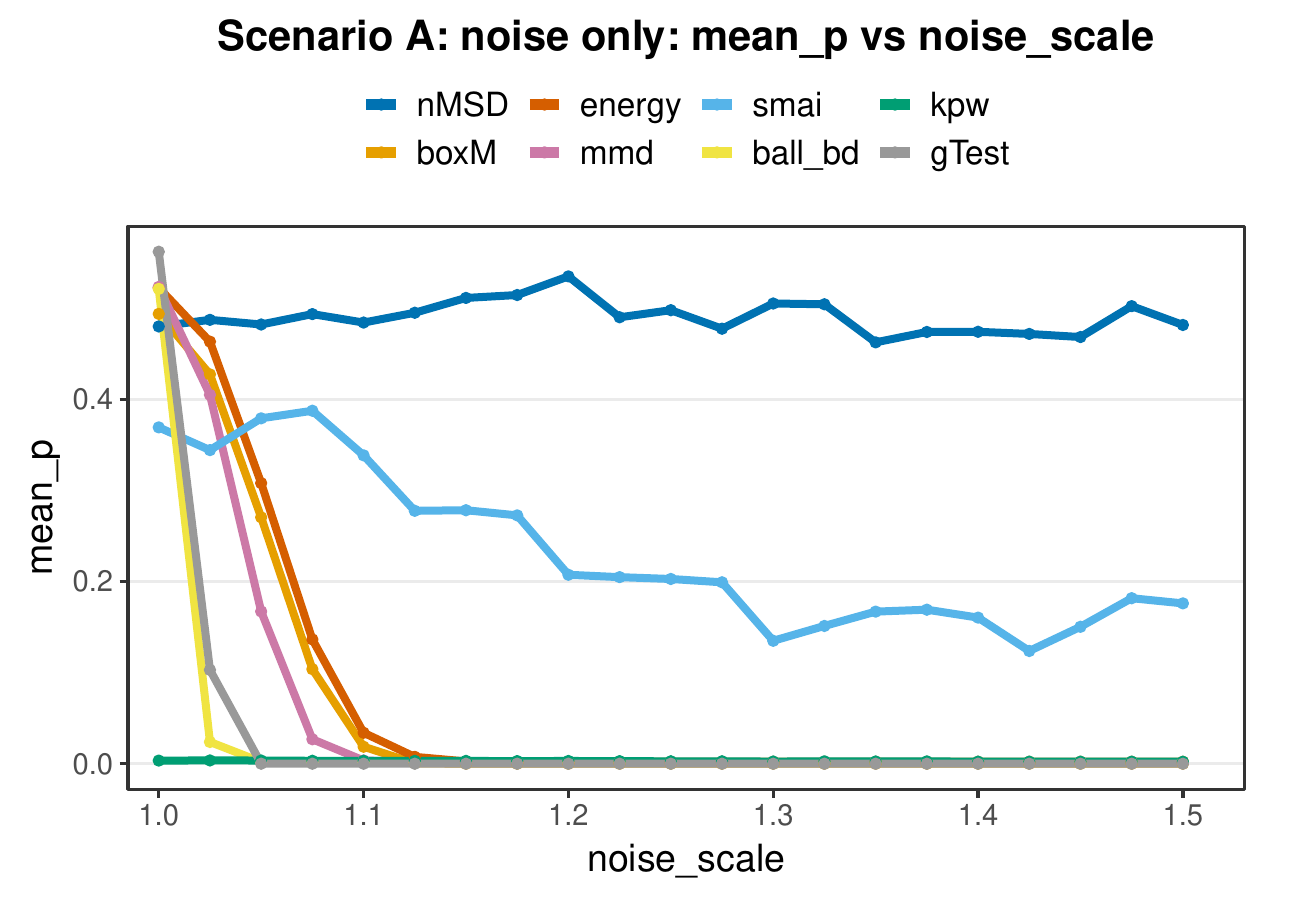}
    \\ \small (a)
  \end{minipage}
  \hfill
  \begin{minipage}{0.48\linewidth}
    \centering
    \includegraphics[width=\linewidth]{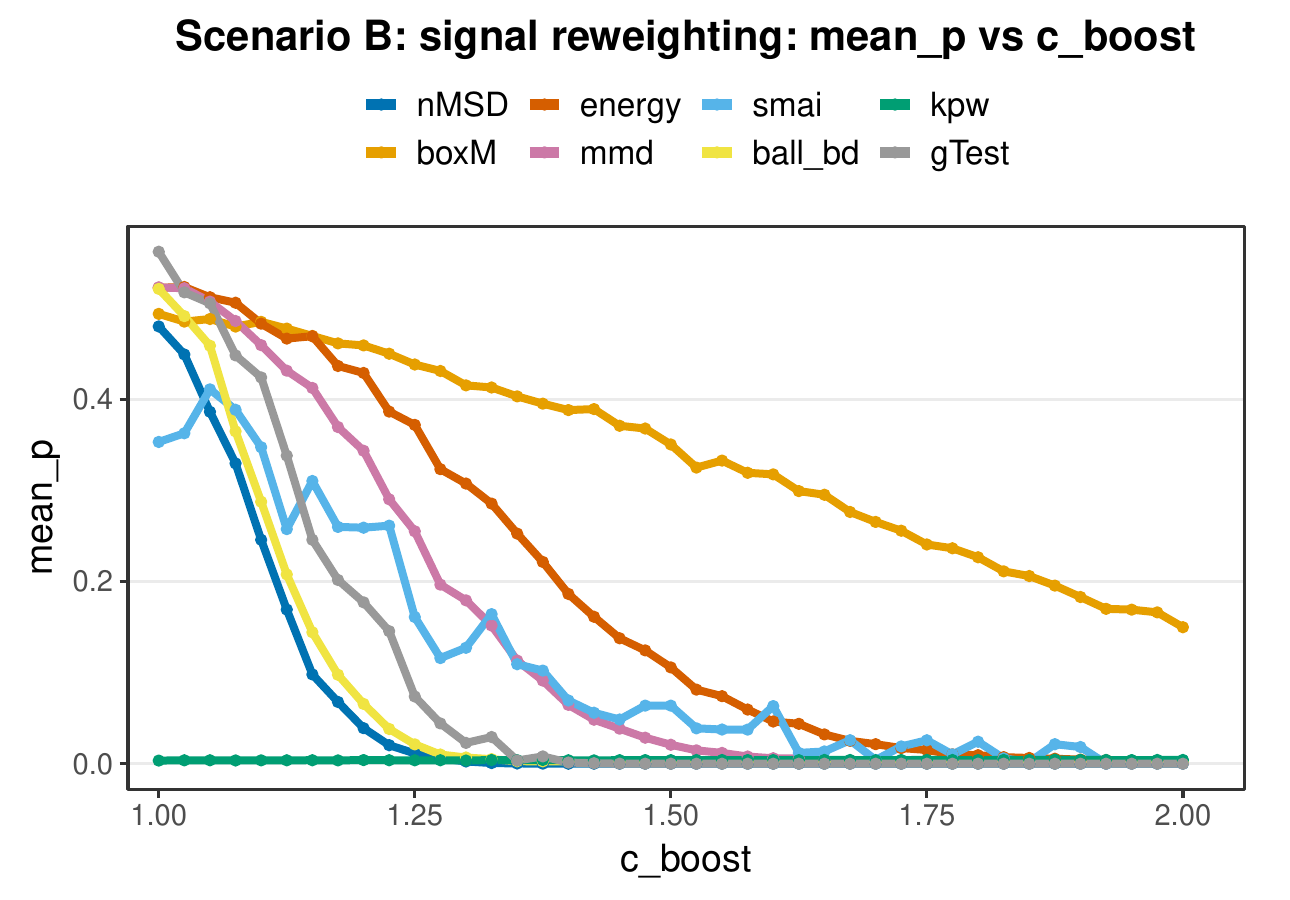}
    \\ \small (b)
  \end{minipage}

  \caption{\textbf{Simulation summary across methods.} 
 Panel (a) shows the mean p-value for identical signals under varying noise scales, and Panel (b) shows the mean p-value under fixed noise with increasing signal separation.}
  \label{fig:p}
\end{figure}

Figure \ref{fig:exp1} contains simulation results corresponding to our asymptotic null distribution. Figure \ref{fig:exp2} compares the theoretical power of our proposed test with the empirical power in simulations, showing their strong conformity. Figure~\ref{fig:p} contains simulation results corresponding to our asymptotic theory. Panel~(a) displays the mean $p$-values under the null distribution across varying noise scales. Panel~(b) compares the theoretical mean $p$-values of our proposed test with the empirical results in simulations, showing their strong conformity and validating our power analysis.

\subsection{Block-structured variance pattern in real biological data}
\label{sec:supp-panc8}

To illustrate the block-structured variance patterns that motivate our model, we analyze the \texttt{panc8} reference dataset distributed with the \texttt{SeuratData} package. The \texttt{panc8} object combines eight human pancreas single-cell RNA-seq datasets across five sequencing technologies (CEL-Seq, CEL-Seq2, Fluidigm C1, inDrops, and Smart-seq2) into a single harmonized Seurat object, with curated cell-type annotations for major endocrine, exocrine, vascular, immune, and stromal populations\citep{panc8data}.

For our analysis, we group cells by sequencing technology using the \texttt{tech} label, which yields five technology-specific subsets: CEL-Seq (\texttt{celseq}), CEL-Seq2 (\texttt{celseq2}), Fluidigm C1 (\texttt{fluidigmc1}), inDrops (\texttt{indrop}), and Smart-seq2 (\texttt{smartseq2}). Within each subset, we apply the standard Seurat preprocessing workflow used for this dataset: library-size normalization of counts, log-transformation of expression values, and identification of highly variable genes. Following the Seurat integration vignette, we then define a common set of 1{,}000 genes that are repeatedly variable across technologies. This shared feature set focuses the analysis on robust biological signals rather than platform-specific noise.

Using this shared set of variable genes, we construct a gene--gene correlation matrix for each technology. For a given technology, we restrict the expression matrix to the shared genes, remove genes with zero variance across cells, and rank the remaining genes by their empirical variance. From this ranking, we retain up to 1{,}000 of the most variable genes (or all available genes if fewer than 1{,}000 remain) and compute the Pearson correlation coefficient between the expression profiles of every pair of selected genes. The resulting correlation matrix is visualized as a heatmap after hierarchical clustering of genes. We use a yellow--orange color scale, ranging from light yellow for weak correlations to dark orange--brown for strong positive correlations, to highlight groups of highly correlated genes.

\begin{figure*}[t!] 
  \centering

  \begin{minipage}[t]{0.45\textwidth}
    \centering
    \includegraphics[width=0.8\linewidth]{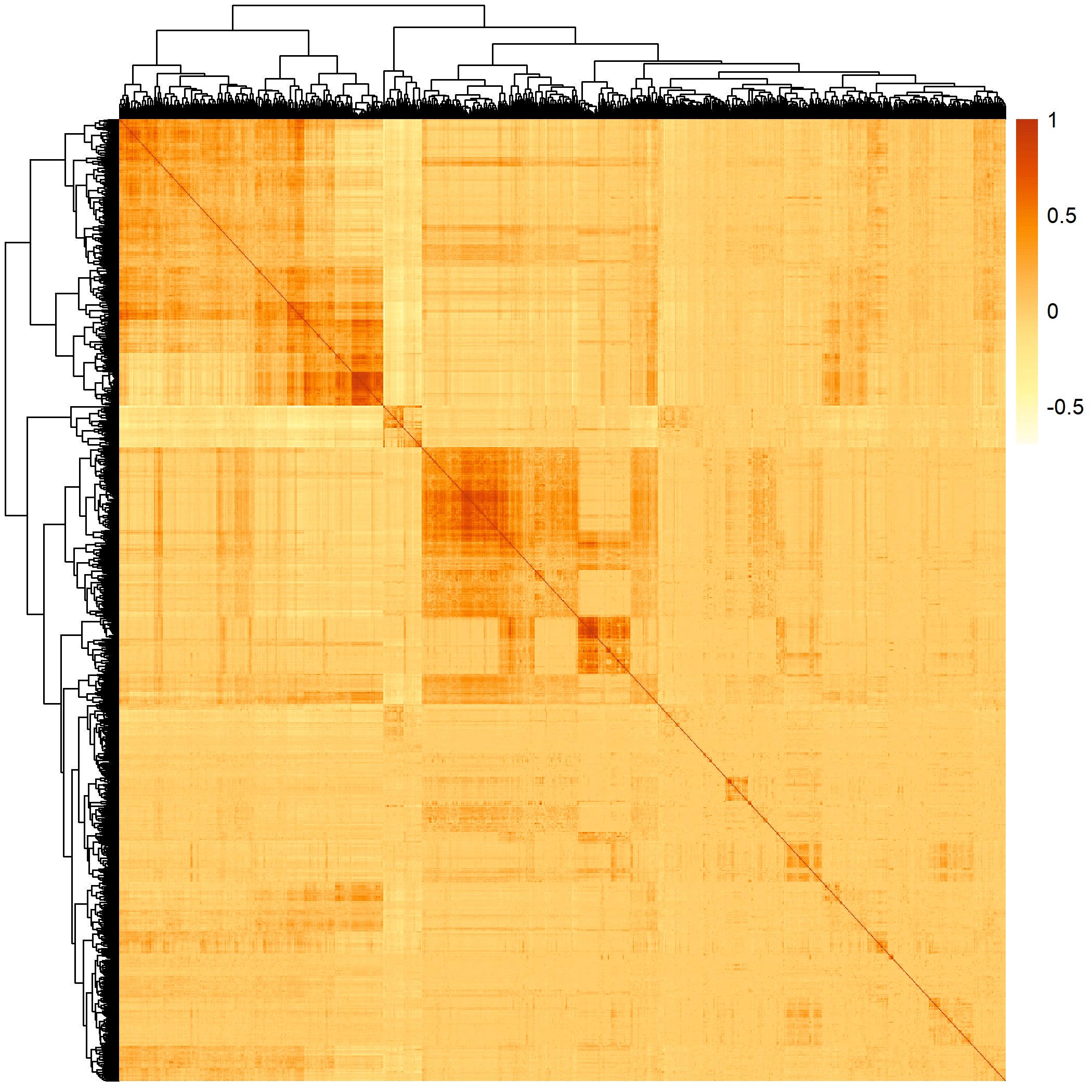}
    \centerline{\small (a)}

    \includegraphics[width=0.8\linewidth]{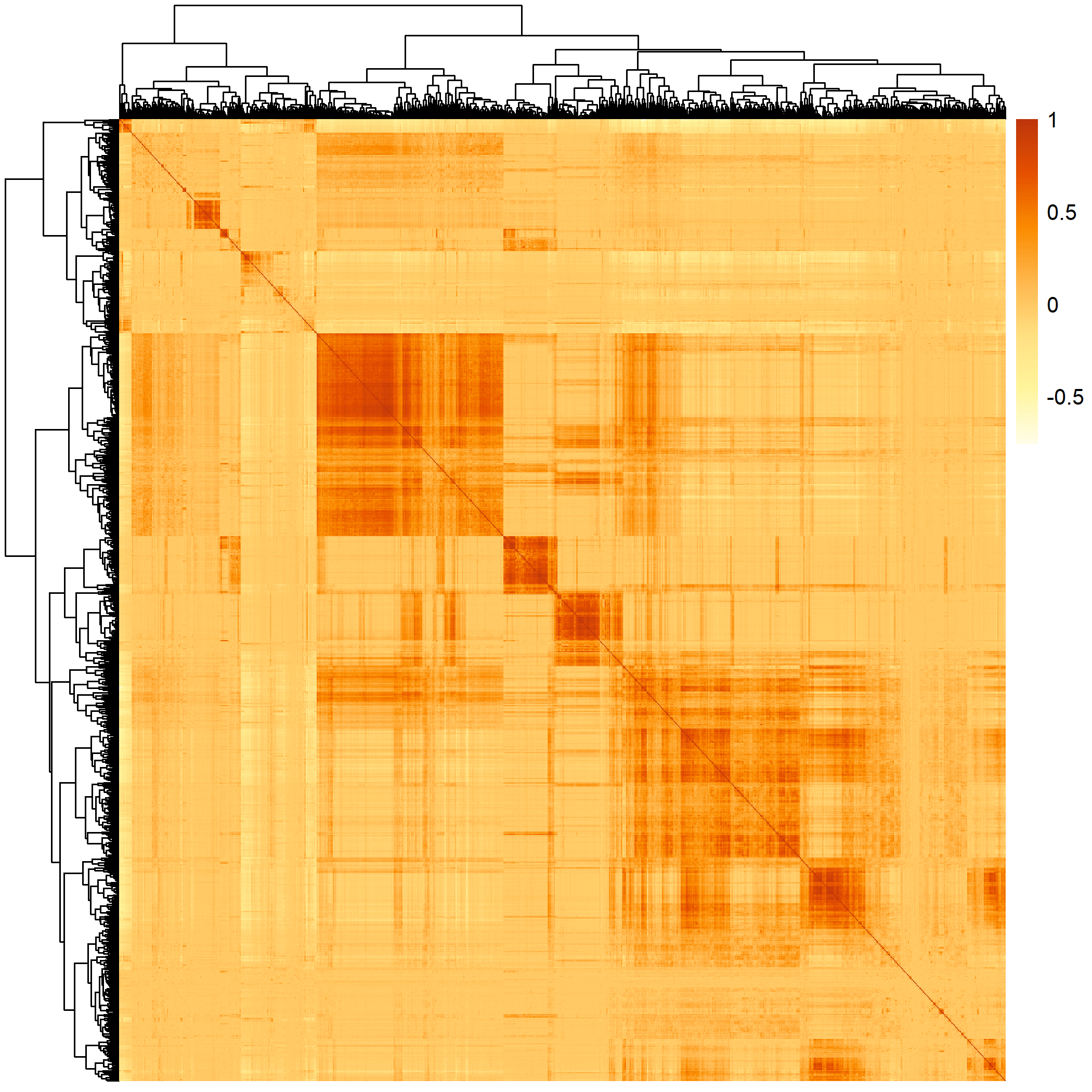}
    \centerline{\small (b)}
    
    \includegraphics[width=0.8\linewidth]{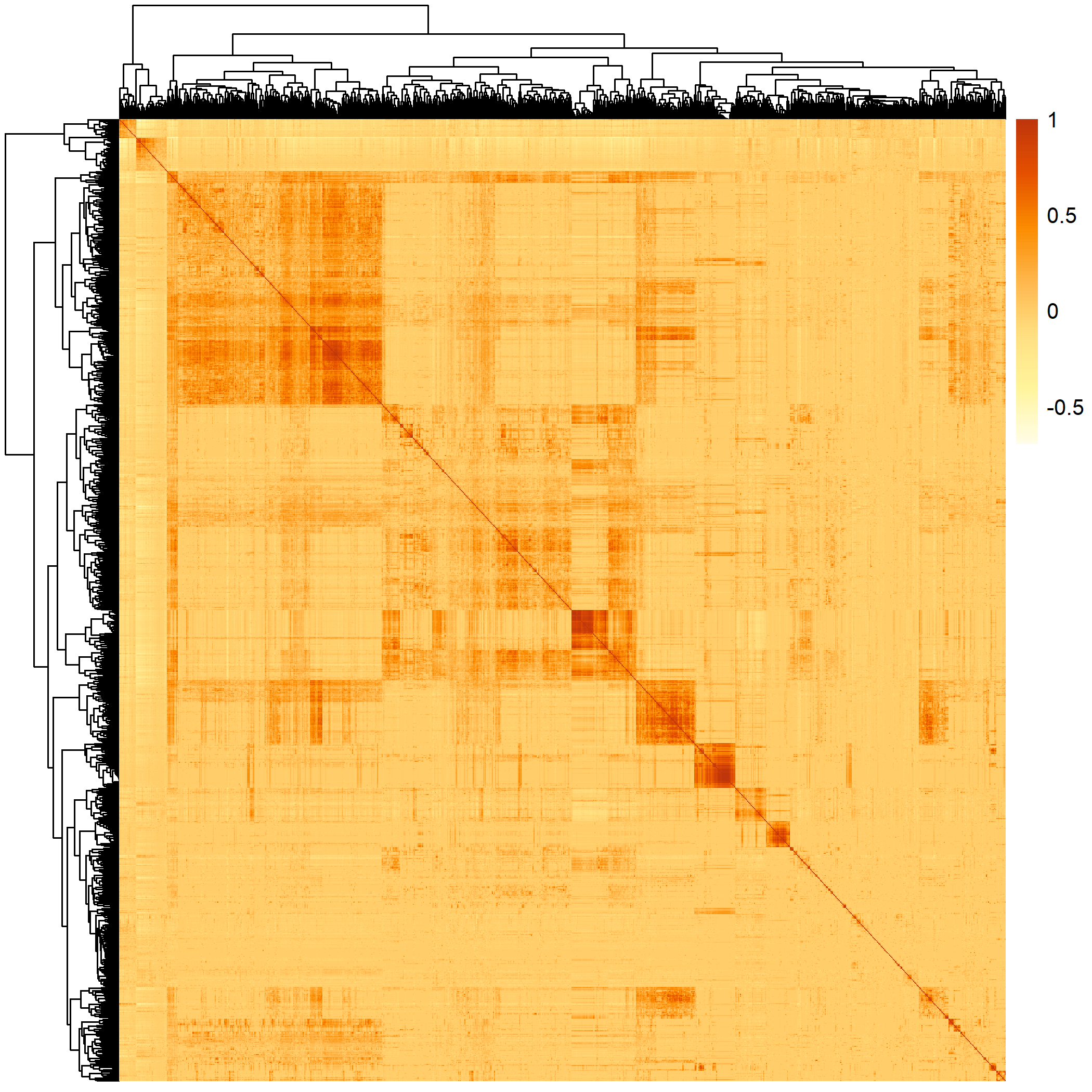}
    \centerline{\small (c)}
  \end{minipage}
  \hfill 
  \begin{minipage}[t]{0.45\textwidth}
    \centering
    \includegraphics[width=0.8\linewidth]{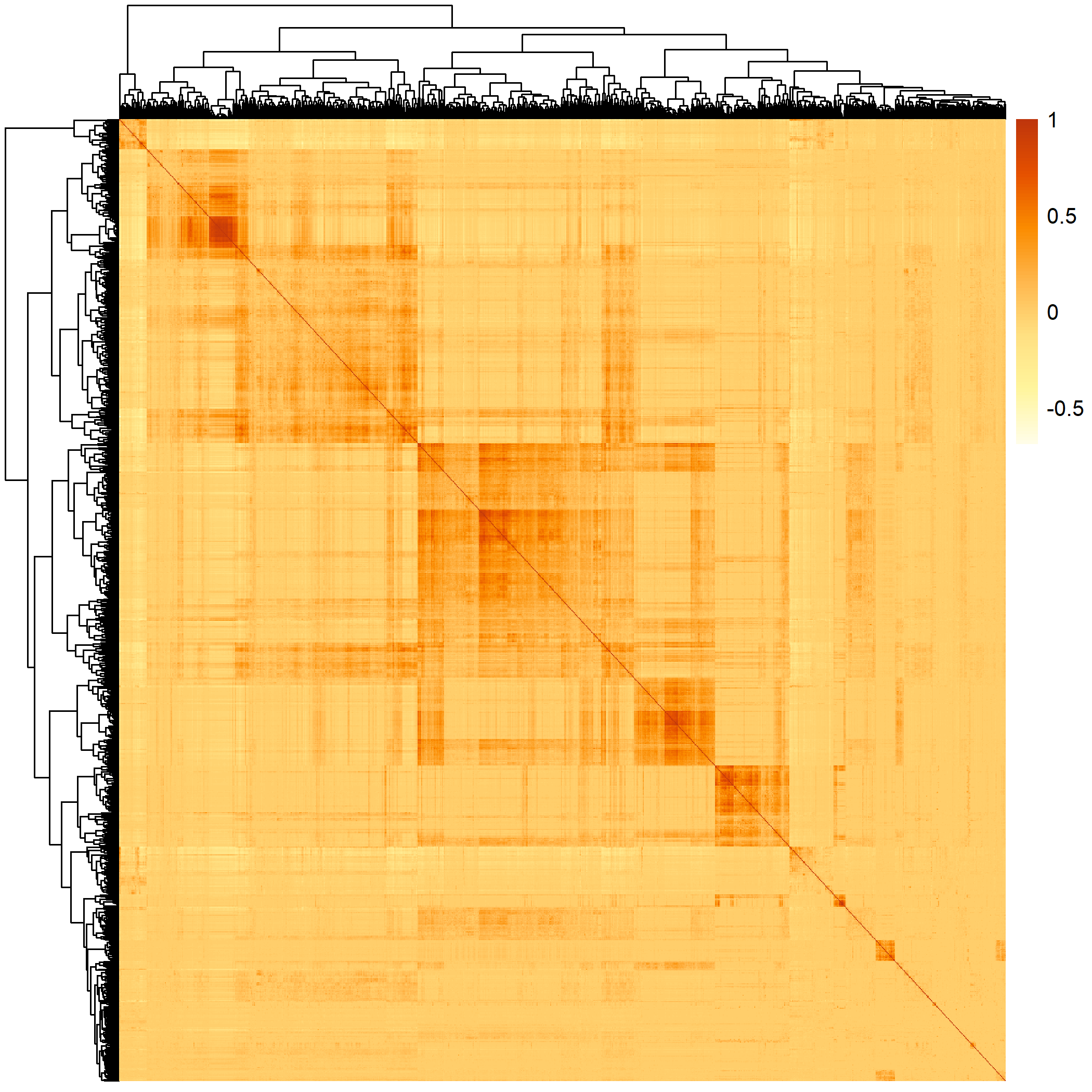}
    \centerline{\small (d)}

    \includegraphics[width=0.8\linewidth]{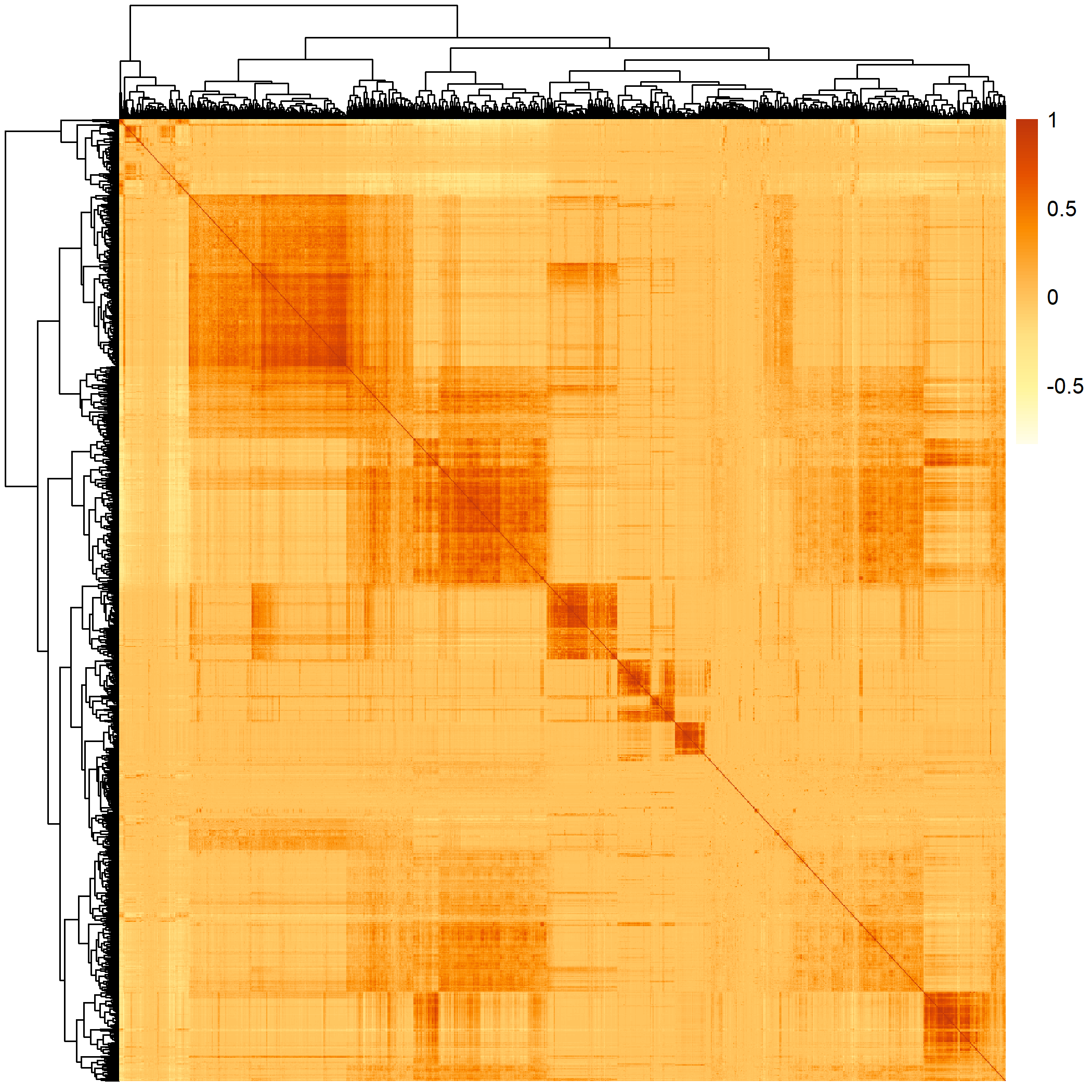}
    \centerline{\small (e)}

  \end{minipage}

  \vspace{1em}
  \caption{\textbf{Gene--gene Pearson correlation heatmaps for the five technology-specific subsets.} 
  Each panel is based on up to 1,000 shared, highly variable genes: (a) CEL-Seq, (b) CEL-Seq2, (c) Fluidigm C1, (d) inDrops, and (e) Smart-seq2. The block-like structures correspond to modular genetic programs.}
  \label{fig:panc8-heatmaps}
\end{figure*}

The five heatmaps in Figure~\ref{fig:panc8-heatmaps} (one for each technology) show clear block-like patterns: subsets of genes form tightly correlated clusters that appear as distinct blocks along the diagonal of the correlation matrix. These blocks are observed consistently across all five technologies, despite substantial differences in experimental protocol and sequencing chemistry. This provides empirical support for our assumption that noise variance can be modeled with a block-structured constraint, reflecting underlying genetic programs in which functionally related genes share similar variance and correlation profiles.

\section{An example for bias correction}

The advantage of our proposed approach in correcting the intrinsic bias of 
$\{\lambda_k\}$ through the inversion map $\theta$ can be readily illustrated 
with the following numerical example.

\begin{example}[Bias correction by the inversion map \texorpdfstring{$\theta^{-1}$}{theta^{-1}}]
\label{ex:theta-skip-bias}

We compare two pipelines for estimating the normalized signal strengths
\(\Pi=(\pi_1,\pi_2,\pi_3)\) from the sample covariance \(Q=YY^\top/N\):
(i) \emph{bias correction}—inverting the outlier map to obtain
\(\widehat{\xi}_j=\theta^{-1}(\lambda_j)\) and then set
\(\widehat{d}_j^2=-1/g(\widehat{\xi}_j)\); and
(ii) \emph{naive plug-in}—directly using $\lambda_j$ to estimate $\xi_j$ and setting \(\widehat{d}_j^2=-1/g(\lambda_j)\). Finally normalize $\{\widehat d_j^{\,2}\}$ to get \(\widehat\Pi\).
The only difference between the two is whether \(\theta^{-1}\) is applied.

We take \(p=800\), \(N=400\), \(r=3\), isotropic noise \(\Sigma=I_p\),
and signal directions \(V\in\mathbb O(p,r)\) drawn uniformly at random.
Signal samples lie on an ellipsoid in the signal subspace with axes
\(D=(4,3,2)\); each column is \(Z=\sqrt r\,V\,\mathrm{diag}(D)\,v\),
\(v\sim\mathrm{Unif}(\mathbb S^{r-1})\).
Observations are \(Y=X+S_{\mathrm{samp}}\) with
\(X_{ij}\stackrel{\mathrm{i.i.d.}}{\sim}\mathcal N(0,1)\).
We run \(500\) independent replications and use the oracle diagonal noise, so
any difference is due solely to \(\theta^{-1}\).

Figure~\ref{fig:ablation-pi} overlays the empirical densities of \(\pi_k\)
for \(k=1,2,3\) under the two pipelines, together with vertical lines at the
truth \(\Pi^\star=(\frac{16}{29},\frac{9}{29},\frac{4}{29})\).
The bias-corrected estimator is tightly centered at the ground truth for all three
components, whereas the naive plug-in estimator shows systematic shifts, which is predicted by the map $\theta(\cdot)$:
\(\pi_1\) is biased downward and \(\pi_3\) is biased upward, with a mild
upward bias for \(\pi_2\).
This experiment highlights the necessity of applying \(\theta^{-1}\) step even when
\(\Sigma\) is known.

\begin{figure}[t]
  \centering
  \includegraphics[width=.95\linewidth]{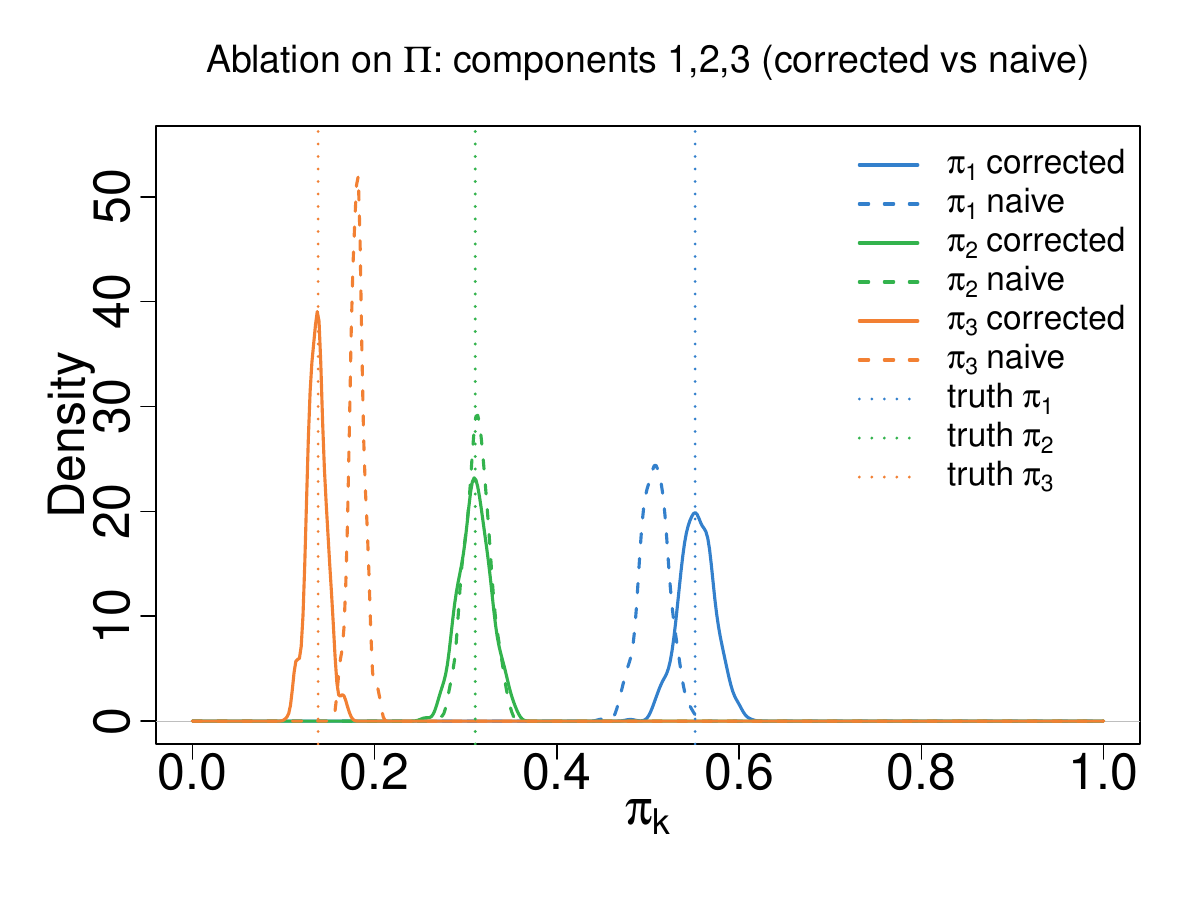}
  \caption{Ablation on \(\Pi\) (components \(1,2,3\)): corrected
  (solid) vs.\ naive (dashed). Dotted lines mark the truth
  \(\Pi^\star=(16,9,4)/29\). Skipping \(\theta^{-1}\) over-weights weaker
  components because sample spikes are inflated relative to the population
  locations.}
  \label{fig:ablation-pi}
\end{figure}
\end{example}

\section{Proof of \MainSubsecTwentyFour \ and Section \ref{subsec:41}}
\subsection{Proof of \MainPropPCA}
\begin{proof}
Under \MainAsmpSignal, the embedded manifold
$\iota(\M)\subset\R^p$ has Euclidean dimension $n$, so
$V:=\operatorname{span}(\iota(\M))$ is an $n$-dimensional subspace and
$\operatorname{rank}(M)\le n$, with $M$ vanishing on $V^\perp$.
Restrict $M$ to $V$; this is a symmetric positive semidefinite operator on an
$n$-dimensional space, so there exists an orthonormal basis
$\{\eta_1,\dots,\eta_n\}$ of $V$ consisting of eigenvectors with eigenvalues
$\mathcal D_1\ge\cdots\ge\mathcal D_n\ge0$.
Extend $\{\eta_1,\dots,\eta_n\}$ to an orthonormal basis of $\R^p$, and let
$R\in\mathbb O(p)$ be the orthogonal matrix whose first $n$ rows are
$\eta_1^\top,\dots,\eta_n^\top$.
For a generic noiseless sample $S\sim\rho$ set $Z:=RS$ and write
$Z=(f_1(S),\dots,f_n(S),0,\dots,0)$, which defines the coordinate functions
$f_1,\dots,f_n$.
Then $\E[Z]=R\,\E[S]=0$, so $\E[f_i(S)]=0$, and
\[
  \E[ZZ^\top] \;=\; RMR^\top
  \;=\; \diag(\mathcal D_1,\dots,\mathcal D_n,0,\dots,0),
\]
which implies
$\E[f_i(S)f_j(S)]=0$ for $i\neq j$ and
$\E[f_i(S)^2]=\mathcal D_i$.
Thus $\{\mathcal D_i\}_{i=1}^n$ are exactly the nonzero eigenvalues of $M$ in
these coordinates, as claimed.
\end{proof}

\subsection{Proof of \Mainprolowrank}
\begin{proof}
Under \MainAsmpSignal, each noiseless sample $S_i$ lies in
$\iota(\M)\subset\R^p$, whose linear span
$V:=\operatorname{span}(\iota(\M))$ has dimension $n$.
Hence all columns of $S_\samp=[S_1,\dots,S_N]$ lie in the same
$n$-dimensional subspace $V$, so $\rank(S_\samp)\le n$.

For the second claim, consider the case $N\ge n$.
\MainAsmpSignal \ states that the sampling measure $\rho$ on
$\iota(\M)$ has a density that is bounded and bounded away from zero on the
manifold.
In particular, when expressed in local Euclidean coordinates on $V$, the law
of $S$ is absolutely continuous with respect to Lebesgue measure on an
$n$-dimensional region.
The event that $S_1,\dots,S_n$ are linearly dependent corresponds (in such
coordinates) to the zero set of a nontrivial polynomial (the determinant of
the $n\times n$ matrix of their coordinates), hence has Lebesgue measure
zero and therefore probability zero.
Thus with probability one, $S_1,\dots,S_n$ are linearly independent, so
$\rank(S_\samp)\ge n$ whenever $N\ge n$.
Combining with $\rank(S_\samp)\le n$ gives $\rank(S_\samp)=n$ almost surely.
\end{proof}

\subsection{Proof of \MainPropmatspec}

By the orthogonal change of coordinates described below \MainAsmpSignal,
there exists $R\in\mathbb O(p)$ such that
\[
  (R\circ\iota)(x)=(f_1(x),\dots,f_n(x),0,\dots,0)\in\R^p
,\text{for all }x\in\M.
\]
For each noiseless sample $S_j=\iota(X_j)$, we have
\[
  R S_j = (f_1(X_j),\dots,f_n(X_j),0,\dots,0)^\top.
\]
Define the $n\times N$ matrix
\[
  W := \begin{pmatrix}
        f_1(X_1) & \cdots & f_1(X_N) \\
        \vdots   &        & \vdots   \\
        f_n(X_1) & \cdots & f_n(X_N)
      \end{pmatrix}
      \in\R^{n\times N},
\]
so that $W$ collects the first $n$ coordinates of the rotated samples
$R S_j$.

Let
\[
  \mathcal D_i := \E\big[f_i(X)^2\big],\quad i=1,\dots,n,
\]
and set the diagonal matrix
\[
  \mathcal D := \diag(\mathcal D_1,\dots,\mathcal D_n).
\]
As explained in \MainPropPCA, these are exactly
the nonzero eigenvalues of $M=\E[S S^\top]$:
\[
  \{\mathcal D_i\}_{i=1}^n = \mathrm{spec}(M)\setminus\{0\}.
\]

Assume:
(i) $\mathcal D_{\min}:=\min_{1\le i\le n}\mathcal D_i>0$;
(ii) the eigengap
\[
\begin{aligned}
  &\delta:=\min_{1\le i\le n}
    \min\{\mathcal D_{i-1}-\mathcal D_i,\ \mathcal D_i-\mathcal D_{i+1}\}>0,\\
  &\quad(\mathcal D_0:=+\infty,\ \mathcal D_{n+1}:=0);
\end{aligned}
\]
(iii) the coordinates $f_i(X)$ are sub-Gaussian (or bounded), so that a
matrix Bernstein bound yields, with probability at least $1-N^{-c}$,
\begin{equation}\label{eq:Eop}
  \Big\|\frac{1}{N}WW^\top-\mathcal D\Big\|_{\op}\ \le\ \varepsilon,
  \qquad \varepsilon\ :=\ C\,\sqrt{\frac{\log N}{N}}.
\end{equation}

Let the thin SVD of $W$ be
\[
  W = U\,D_\samp\,V^\top,
\]
where $U\in\R^{n\times n}$ is orthogonal,
$D_\samp=\diag(d_1,\dots,d_n)$, and
$V=[v^{(1)},\dots,v^{(n)}]\in\R^{N\times n}$ has orthonormal columns.
Since $S_{\samp}=R^\top W$ with $R$ orthogonal, $S_{\samp}$ and $W$ share
the same nonzero singular values, so $d_i$ coincides with the singular
values in \MainPropmatspec.

\medskip\noindent
\emph{Step 1: Singular values.}
By Weyl's inequality and~\eqref{eq:Eop},
\begin{equation}\label{eq:sv}
  \Big|\frac{d_i^2}{N}-\mathcal D_i\Big|\ \le\ \varepsilon,
  \qquad 1\le i\le n,
\end{equation}
and hence, using $\mathcal D_i\ge\mathcal D_{\min}>0$,
\[
  \Big|\frac{d_i}{\sqrt N}-\sqrt{\mathcal D_i}\Big|
  \le \frac{\varepsilon}{\sqrt{\mathcal D_{\min}}}.
\]
This proves the statement on singular values in
\MainPropmatspec(ii).

\medskip\noindent
\emph{Step 2: Eigenvectors of $WW^\top$.}
Write
\[
  \frac{1}{N}WW^\top = \mathcal D + E,\qquad \|E\|_{\op}\le\varepsilon.
\]
By the Davis--Kahan $\sin\Theta$ theorem applied with gap $\delta$ between
the diagonal entries of $\mathcal D$, we obtain
\begin{equation}\label{eq:dk}
  \|U-I_n\|_{\op}\ \le\ C\,\frac{\varepsilon}{\delta}.
\end{equation}
Let $u_i$ denote the $i$th column of $U$ and $e_i$ the $i$th standard
basis vector in $\R^n$.

\medskip\noindent
\emph{Step 3: Coordinate-wise approximation of right singular vectors.}
From the SVD identity $W^\top U = V D$, we have
\[
  v^{(i)} = \frac{1}{d_i} W^\top u_i,\qquad 1\le i\le n.
\]
Let $v_{ji}$ denote the $(j,i)$ entry of $V$, i.e. the $j$th coordinate of
$v^{(i)}$, and let $w_j$ be the $j$th column of $W$:
\[
  w_j = (f_1(X_j),\dots,f_n(X_j))^\top.
\]
Then
\[
\begin{aligned}
  &v_{ji}
  = \frac{1}{d_i}\,u_i^\top w_j
  = \frac{1}{d_i}\,(u_i-e_i)^\top w_j
    + \frac{1}{d_i}\,e_i^\top w_j \\
  &= \underbrace{\frac{1}{d_i}\,(u_i-e_i)^\top w_j}_{\mathrm{(A)}}
    + \underbrace{\frac{1}{d_i}\,f_i(X_j)}_{\mathrm{(B)}}.
\end{aligned}
\]
Subtracting the population target $f_i(X_j)/\sqrt{N\mathcal D_i}$ yields
\begin{equation}\label{eq:coord-decomp}
  v_{ji}-\frac{f_i(X_j)}{\sqrt{N\mathcal D_i}}
  = \underbrace{\frac{1}{d_i}\,(u_i-e_i)^\top w_j}_{\mathrm{(A)}}
    + \underbrace{\Big(\frac{1}{d_i}-\frac{1}{\sqrt{N\mathcal D_i}}\Big)\,f_i(X_j)}_{\mathrm{(B)}}.
\end{equation}

By~\eqref{eq:dk}, $\|u_i-e_i\|_2\le C\,\varepsilon/\delta$. Hence
\[
\begin{aligned}
  &\max_{1\le j\le N}|\mathrm{(A)}|
  \le \frac{1}{d_i}\,\|u_i-e_i\|_2\,\max_{1\le j\le N}\|w_j\|_2 \\
  &\le \frac{C}{\sqrt N}\,\frac{\varepsilon}{\delta}\,
       \max_{1\le j\le N}\|w_j\|_2,
\end{aligned}
\]
where in the last step we used $d_i \asymp \sqrt{N\mathcal D_i}$.

From~\eqref{eq:sv}, we also have
\[
  \Bigl|\frac{1}{d_i}-\frac{1}{\sqrt{N\mathcal D_i}}\Bigr|
  \le C\,\frac{\varepsilon}{\sqrt{N}},
\]
for some constant $C$ depending only on $\mathcal D_{\min}$, so
\[
  \max_{1\le j\le N}|\mathrm{(B)}|
  \le C\,\frac{\varepsilon}{\sqrt{N}}\,
        \max_{1\le j\le N}|f_i(X_j)|.
\]

Under the sub-Gaussian (or bounded) coordinate assumption, there exist
constants $C,c>0$ such that, with probability at least $1-N^{-c}$,
\[
\begin{aligned}
 & \max_{1\le j\le N}\|w_j\|_2 \le C\sqrt{\log N},\\
  &\max_{1\le j\le N}|f_i(X_j)| \le C\sqrt{\mathcal D_i\log N}.
\end{aligned}
\]
Combining the last four displays and using
$\varepsilon = C\sqrt{(\log N)/N}$, we obtain, with probability at least
$1-N^{-c}$,
\[
  \max_{1\le j\le N}\Big|v_{ji}-\frac{f_i(X_j)}{\sqrt{N\mathcal D_i}}\Big|
  \le C'\,\frac{\log N}{N},
\]
for some constant $C'>0$ depending only on the model parameters. This
proves the desired coordinatewise approximation in
\MainPropmatspec(iii).

\subsection{Proof of Lemma \ref{lem:B-deloc}}

\begin{proof}
Let $S_{\samp}=U_{\samp} D_\samp V_{\samp}^\top$ be the thin SVD of the sample signal
matrix, with $U_{\samp}\in\R^{p\times n}$ having orthonormal columns,
$D_\samp=\diag(d_1,\dots,d_n)$, and $V_{\samp}\in\R^{N\times n}$ having orthonormal
columns.  Write $v_j^\top$ for the $j$th row of $V_{\samp}$.

By \MainPropmatspec(iv), there exist constants
$C_V,c>0$ such that, with probability at least $1-N^{-c}$,
\begin{equation}\label{eq:V-deloc}
  \max_{1\le i\le n}\max_{1\le j\le N} |V_{\samp}(j,i)|
  \;\le\; C_V\,\frac{\sqrt{\log N}}{\sqrt N}.
\end{equation}

By definition of $B$ and the SVD of $S_{\samp}$,
\[
  B
  = \frac{1}{\sqrt N}\,S_{\samp}^\top \psi_k
  = \frac{1}{\sqrt N}\,V_{\samp} D_\samp U_{\samp}^\top \psi_k.
\]
Set $c := U_{\samp}^\top \psi_k \in \R^n$.
Since $U_{\samp}$ has orthonormal columns and $\|\psi_k\|_2=1$, we have
$\|c\|_2 \le 1$ and hence $\|c\|_1 \le \sqrt n$.

For the $j$th coordinate of $B$,
\[
  B_j
  = \frac{1}{\sqrt N}\,v_j^\top D_\samp c,
\]
so
\[
\begin{aligned}
 & |B_j|
  \;\le\; \frac{1}{\sqrt N}\,\|v_j\|_\infty\,\|D_\samp c\|_1 \\
  &\le\; \frac{1}{\sqrt N}\,\|v_j\|_\infty\,\|D_\samp\|_{\op}\,\|c\|_1\\
  &\le\; \frac{\sqrt n}{\sqrt N}\,\|v_j\|_\infty\,\|D_\samp\|_{\op}.
\end{aligned}
\]
On the high-probability event \eqref{eq:V-deloc}, this yields
\[
  \max_{1\le j\le N} |B_j|
  \;\le\; C_V\,\|D_\samp\|_{\op}\,\frac{\sqrt{n\log N}}{N}.
\]
Since $N\ge1$ implies $1/N \le 1/\sqrt N$, the above bound is in fact
stronger than the claimed inequality
\[
  \|B\|_\infty
  \;\le\; C_V\,\|D_\samp\|_{\op}\,\frac{\sqrt{n\log N}}{\sqrt N}.
\]
This proves the lemma.
\end{proof}

\section{Proof of Lemma \ref{lem:signal-resampling-general-new}}
\begin{proof}
Let $S$ denote a generic copy of the sample vectors $S_1,\ldots,S_N$,
and write
\[
  b_k := \E\big[(v_k^\top S)^2\big]
  = v_k^\top M v_k,
  \qquad
  M := \E[SS^\top].
\]
By definition,
\[
  \delta_k
  = \frac1N\sum_{i=1}^N\Big\{(v_k^\top S_i)^2 - b_k\Big\},
\]
so $\E[\delta_k]=0$.

\medskip\noindent
\emph{Raw covariance formula \eqref{eq:Gamma-general-raw-new}.}
Using independence of $S_1,\ldots,S_N$,
\[
\begin{aligned}
 & \Cov(\delta_k,\delta_j)\\
  &= \Cov\Big(
       \frac1N\sum_{i=1}^N\big((v_k^\top S_i)^2-b_k\big),\;
       \frac1N\sum_{i'=1}^N\big((v_j^\top S_{i'})^2-b_j\big)
     \Big) \\
  &= \frac1{N^2}\sum_{i=1}^N\sum_{i'=1}^N
       \Cov\big((v_k^\top S_i)^2,(v_j^\top S_{i'})^2\big).
\end{aligned}
\]
For $i\neq i'$, the terms are independent and their covariance is zero,
so only the diagonal terms survive:
\[
\begin{aligned}
 & \Cov(\delta_k,\delta_j)\\
 & = \frac1{N^2}\sum_{i=1}^N
     \Cov\big((v_k^\top S_i)^2,(v_j^\top S_i)^2\big)\\
 & = \frac1N\Cov\big((v_k^\top S)^2,(v_j^\top S)^2\big).
\end{aligned}
\]
Expanding the covariance gives
\[
\begin{aligned}
 & \Gamma^{(\mathrm{sig})}_{kj}
  = \Cov(\delta_k,\delta_j)\\
 & = \frac1N\Big\{
      \E\big[(v_k^\top S)^2(v_j^\top S)^2\big]
      - \E\big[(v_k^\top S)^2\big]\;\E\big[(v_j^\top S)^2\big]
    \Big\},
\end{aligned}
\]
which is the formula \eqref{eq:Gamma-general-raw-new}.

\medskip\noindent
\emph{Tensor form formula\eqref{eq:Gamma-general-tensor-new}.}
Let $T^{(4)} := \E[S^{\otimes4}]$ denote the fourth raw--moment tensor, with
entries $T^{(4)}_{abcd}=\E[S_a S_b S_c S_d]$.
For any $u,v\in\R^p$,
\[
  (u^\top S)^2(v^\top S)^2
  = \Big\langle u\otimes u\otimes v\otimes v,\ S^{\otimes4}\Big\rangle,
\]
so
\[
  \E\big[(v_k^\top S)^2(v_j^\top S)^2\big]
  = \Big\langle v_k\otimes v_k\otimes v_j\otimes v_j,\ T^{(4)}\Big\rangle.
\]
Since $\E[(v_k^\top S)^2]=b_k$, this yields
\[
  \Gamma^{(\mathrm{sig})}_{kj}
  = \frac1N\big\langle v_k\!\otimes v_k\!\otimes v_j\!\otimes v_j,\ T^{(4)}\big\rangle
    - \frac1N\,b_k\,b_j,
\]
which is formula\eqref{eq:Gamma-general-tensor-new}.

\medskip\noindent
Let $\Kc:=\mathrm{Cum}_4(S)$ denote the fourth central cumulant tensor of $S$.
By definition of $\Kc$ for a mean-zero random vector,
\[
  T^{(4)} \;=\; \Kc + \Sym(M\otimes M),
\]
where $\Sym(M\otimes M)$ is the symmetrized tensor with entries
\[
  \Sym(M\otimes M)_{abcd}
  := M_{ab}M_{cd} + M_{ac}M_{bd} + M_{ad}M_{bc}.
\]
Contracting with $v_k\otimes v_k\otimes v_j\otimes v_j$ gives
\[
\begin{aligned}
  &\E\big[(v_k^\top S)^2(v_j^\top S)^2\big]
  = \Big\langle v_k\otimes v_k\otimes v_j\otimes v_j,\ T^{(4)}\Big\rangle \\
  &= \Big\langle v_k\otimes v_k\otimes v_j\otimes v_j,\ \Kc\Big\rangle \\
    & + \Big\langle v_k\otimes v_k\otimes v_j\otimes v_j,\ \Sym(M\otimes M)\Big\rangle \\
  &= \Kc[v_k,v_k,v_j,v_j]
     + (v_k^\top M v_k)(v_j^\top M v_j)
     + 2\,(v_k^\top M v_j)^2.
\end{aligned}
\]
Recall that $\E[(v_k^\top S)^2]=v_k^\top M v_k = b_k$, so
\[
\begin{aligned}
  &\Cov\big((v_k^\top S)^2,(v_j^\top S)^2\big)\\
  &= \E\big[(v_k^\top S)^2(v_j^\top S)^2\big]
     - \E\big[(v_k^\top S)^2\big]\;\E\big[(v_j^\top S)^2\big] \\
  &= 2\,(v_k^\top M v_j)^2 + \Kc[v_k,v_k,v_j,v_j].
\end{aligned}
\]
Dividing by $N$ yields
\[
  \Gamma^{(\mathrm{sig})}_{kj}
  = \frac1N\Big\{\,2\,(v_k^\top M v_j)^2+\Kc[v_k,v_k,v_j,v_j]\Big\},
\]
\end{proof}

\section{Proof of Section \ref{subsec:42}}
\subsection{Proof of Lemma \ref{lem:isotropy}}

\begin{proof}
\emph{Expectation.}
Fix $x\in\R^r$ and set $v:=U_rx\in\R^p$. By the Haar distribution of $U_r$,
for each fixed $x$, the vector $v/\|x\|$ is uniformly distributed on the unit sphere $S^{p-1}$ in $\R^p$.
Therefore
\[
\begin{aligned}
&\E\!\big[(U_rx)^\top A(\xi)\,(U_rx)\,\big|\,\Sigma\big]
=\|x\|^2\,\E\!\big[v_0^\top A(\xi)\,v_0\,\big|\,\Sigma\big]\\
&=\|x\|^2\,\frac{1}{p}\Tr A(\xi)\\
&=\|x\|^2\,g(\xi),
\end{aligned}
\]
where $v_0\sim\Unif(S^{p-1})$. Since this holds for all $x$, we obtain
$\E\big[G_{U_r}(\xi)\,\big|\,\Sigma\big]=g(\xi)I_r$.

\smallskip
\emph{Deviation for a fixed direction.}
Let $x\in S^{r-1}$ and write $Y_x:=x^\top G_{U_r}(\xi)\,x=(U_rx)^\top A(\xi)\,(U_rx)$.
Let $A:=A(\xi)$ for brevity and note $A$ is diagonal in the standard basis.
Represent a uniform sphere vector as $v=z/\|z\|$, where $z\sim N(0,I_p)$ is independent of $\Sigma$.
Then
\[
Y_x \;=\; v^\top A v \;=\; \frac{z^\top A z}{\|z\|^2}.
\]
Set $Z_1:=z^\top A z-\Tr A$ and $Z_2:=\|z\|^2-p$. On the event
$\mathcal E:=\{\,\tfrac12p\le \|z\|^2\le 2p\,\}$ (which satisfies
$\P(\mathcal E^c)\le 2e^{-c_0 p}$ by $\chi^2$-concentration), we have
$\|z\|^{-2}\in[\tfrac1{2p},\tfrac2p]$ and
\[
\Big|\frac{1}{\|z\|^2}-\frac{1}{p}\Big|
=\frac{|Z_2|}{\|z\|^2\,p}
\le \frac{2|Z_2|}{p^2}\qquad\text{on }\mathcal E.
\]
Hence, still on $\mathcal E$,
\begin{align*}
\big|Y_x-g\big|
&=\Big|\frac{z^\top A z}{\|z\|^2}-\frac{\Tr A}{p}\Big|
\le \frac{|Z_1|}{\|z\|^2}+\Tr(|A|)\,\Big|\frac{1}{\|z\|^2}-\frac{1}{p}\Big|  \\
&\le \frac{2}{p}|Z_1|+\frac{2\,\Tr(|A|)}{p^2}|Z_2|.
\end{align*}
By the Hanson--Wright inequality for sub-Gaussian vectors,
\[
\P\!\Big(|Z_1|\ge t\ \Big|\ \Sigma\Big)\ \le\ 2\exp\!\Big(-c_1\min\Big\{\frac{t^2}{\|A\|_F^2},\,\frac{t}{\|A\|_{\op}}\Big\}\Big),
\]
and by $\chi^2$-concentration, $\P(|Z_2|\ge s\sqrt p)\le 2e^{-c_2 s^2}$.
Combining these bounds with the previous display and $\P(\mathcal E^c)\le 2e^{-c_0 p}$, we get, for all $u>0$,
\begin{equation}\label{eq:fixed-x-tail}
\begin{aligned}
&\P\!\Big(\big|Y_x-g\big|\ \ge\ C\Big(\frac{\|A\|_F}{p}\,u+\frac{\|A\|_{\op}}{p}\,u^2+\frac{\Tr|A|}{p^{3/2}}\,u\Big)\ \Big|\ \Sigma\Big)\\
&\ \le\ 4e^{-c u^2}.
\end{aligned}
\end{equation}
Using Cauchy--Schwarz, $\Tr|A|\le \sqrt p\,\|A\|_F$, so the third term is $\lesssim \|A\|_F/p$.

\smallskip
\emph{From a fixed direction to operator norm.}
Let $\mathcal N$ be a $1/4$-net of $S^{r-1}$ with $|\mathcal N|\le 9^r$.
For any symmetric $r\times r$ matrix $B$, $\|B\|_{\op}\le 2\sup_{x\in\mathcal N}|x^\top Bx|$.
Applying this to $B=R(\xi)$ and using \eqref{eq:fixed-x-tail} for each $x\in\mathcal N$,
a union bound with $u\asymp \sqrt{r+\log(1/\delta)}$ yields, with probability $\ge 1-\delta$,
\[
\|R(\xi)\|_{\op}
\ \le\ C\Big(\frac{\|A\|_F}{p}\sqrt r+\frac{\|A\|_{\op}}{p}\,r+\frac{\|A\|_F}{p}\Big).
\]
Since $r$ is fixed, the term $\|A\|_F/p$ is dominated by $(\|A\|_F/p)\sqrt r$.
Moreover, because $r\le p$ we have $r/p\le \sqrt{r/p}$, so the middle term admits the simpler bound
$(\|A\|_{\op}/p)\,r\ \le\ \|A\|_{\op}\sqrt{r/p}$.
This gives the stated bound \eqref{eq:Rop-bound}.

\smallskip

If $\dist(\xi,\spec(\Sigma))\ge c>0$, then $\|A(\xi)\|_{\op}\le c^{-1}$ and
$s_2(\xi)=\sum_{i=1}^p(\sigma_i-\xi)^{-2}\le p\,c^{-2}$.
Hence $\sqrt{s_2(\xi)}/p\le 1/(c\sqrt p)$ and the RHS of bound \eqref{eq:Rop-bound} \
is $O_{\P}(\sqrt{r/p})=o_{\P}(1)$ for fixed $r$.
\end{proof}

\subsection{Proof of Proposition \ref{prop:decouple}}

\begin{proof}
For each $j$, use $\det(I_r+AB)=\det(I_r+BA)$ with $A=D_r$ and $B=D_r G_j$ to get
\begin{equation}\label{eq:det-sym-asymp}
\det(I_r+D_r^2 G_j)=\det(I_r+D_r G_j D_r).
\end{equation}
Set $H_j:=D_r G_j D_r=g_j D_r^2+D_r R_j D_r$. Since $D_r^2$ solves the secular equation \MainEqRankRSecular, the left-hand side of \eqref{eq:det-sym-asymp} vanishes, hence $-1\in\spec(H_j)$. Thus there exists an index $\ell=\ell(j)$ such that
\begin{equation}\label{eq:secular-eig-asymp}
\lambda_{\ell}(H_j)=-1.
\end{equation}

By Lemma \ref{lem:isotropy}, $\|R_j\|_{\op}=o_{\mathbb P}(1)$ for each fixed $j$ (with fixed $r$ and $\dist(\xi_j,\spec(\Sigma))\ge c>0$). 
Moreover, by \MainAsmpSignal \ and \MainPropmatspec, the signal covariance is bounded, hence $\|D_r\|=\max_i d_i=O(1)$. Therefore
\[
\|D_r R_j D_r\|_{\op}\ \le\ \|D_r\|^2\,\|R_j\|_{\op} \ =\ O(1)\cdot o_{\mathbb P}(1)\ =\ o_{\mathbb P}(1).
\]
Applying Weyl's inequality to $H_j=g_j D_r^2+(D_r R_j D_r)$ gives
\begin{equation}\label{eq:weyl-asymp}
\big|\lambda_{\ell}(H_j)-g_j d_\ell^2\big|\ \le\ \|D_r R_j D_r\|_{\op} \ =\ o_{\mathbb P}(1).
\end{equation}
Combining \eqref{eq:secular-eig-asymp} and \eqref{eq:weyl-asymp}, we obtain
\[
g_j d_{\ell(j)}^2 + 1 \ =\ o_{\mathbb P}(1),
\quad\text{hence}\quad
d_{\ell(j)}^2 \ =\ -\frac{1}{g_j} + o_{\mathbb P}(1).
\]
Relabeling indices by a permutation $\pi$ yields
\[
d_{\pi(j)}^2\ =\ -\frac{1}{g_j} + o_{\mathbb P}(1),\qquad j=1,\dots,r,
\]
which proves the first claim (and the multiset convergence).

For uniqueness, assume in addition that the numbers $\{-1/g_j\}_{j=1}^r$ are pairwise separated by at least $\Delta>0$. If two distinct $j\neq j'$ shared the same $\ell$, then by the triangle inequality,
\[
\Big|-\frac{1}{g_j}+\frac{1}{g_{j'}}\Big|
\ \le\ \Big|d_{\ell}^2+\frac{1}{g_j}\Big|+\Big|d_{\ell}^2+\frac{1}{g_{j'}}\Big|
\ =\ o_{\mathbb P}(1),
\]
which contradicts the fixed separation $\Delta$ with probability tending to one. Hence the matching $\pi$ is unique with probability tending to one.
\end{proof}

\section{Proof of Section \ref{Deltamethod} \ and Section \ref{subsec:variance}}
\subsection{Proof of Proposition \ref{pro:joint-clt-uncond}}

\begin{proof}
Recall that $\boldsymbol\lambda=(\lambda_1,\dots,\lambda_r)^\top$ denotes the
vector of outlier eigenvalues of the sample covariance 
$Q:=YY^\top/N$, and that $\boldsymbol\xi
=(\xi_1,\dots,\xi_r)^\top$ are the corresponding outlier
eigenvalues of the population matrix $\widetilde\Sigma:=\Sigma+M$. We work
under \MainAsmpNoise ,\ \MainAsmpSignal \ and \MainAsmpSignalStrength.

\medskip
\noindent\textbf{Step 1: Law of total variance.}
Conditioning on the signal sample $S_{\samp}$, the randomness in
$Q$ comes solely from the noise matrix $X$. Thus
\[
\Var(\boldsymbol\lambda)
=\E\!\bigl[\Cov(\boldsymbol\lambda\mid S_{\samp})\bigr]
+\Cov\!\bigl(\E[\boldsymbol\lambda\mid S_{\samp}]\bigr).
\]
We analyse the two terms separately.

\medskip
\noindent\textbf{Step 2: Conditional (noise) term.}
For fixed $S_{\samp}$ and $\Sigma$, the ``pseudo–population'' covariance is
\[
\SigmaS \;:=\; \Sigma + \widehat M,
\qquad \widehat M := \frac{1}{N}S_{\samp}S_{\samp}^\top,
\]
whose outlier eigenvalues we denote by
\[
{\boldsymbol\pssig}
:=\bigl(\pssig_1,\dots,\pssig_r\bigr)^\top.
\]
Conditionally on $S_{\samp}$, the model is a standard spiked covariance
model with population spikes $\boldsymbol\pssig$ and noise covariance
$\Sigma$. By the spiked–covariance CLT (finite eighth moment, fixed spike separation; see Lemma \ref{lem:clt-lambda-con}), we have
\[
\sqrt N\Bigl(\boldsymbol\lambda
             -\theta\!\bigl(\boldsymbol\pssig\bigr)\Bigr)
\ \Rightarrow\ \mathcal N_r\bigl(0,\,V^{\rm cond}\bigr)
\quad\text{conditionally on }S_{\samp},
\]
where $V^{\rm cond}=V^{\rm cond}(S_{\samp})$ is the limiting conditional
covariance matrix (whose explicit decomposition is given in
Lemma \ref{lem:clt-lambda-con}). In particular,
\[
\Cov(\boldsymbol\lambda\mid S_{\samp})
=\frac{1}{N}\,V^{\rm cond} \;+\; o\!\Bigl(\frac{1}{N}\Bigr),
\]
so taking expectations in $(S_{\samp},\Sigma)$ yields
\begin{equation}\label{eq:Var-first-term}
\E\!\bigl[\Cov(\boldsymbol\lambda\mid S_{\samp})\bigr]
=\frac{1}{N}\,\E\!\bigl[V^{\rm cond}\bigr]
\;+\;o\!\Bigl(\frac{1}{N}\Bigr).
\end{equation}

\medskip
\noindent\textbf{Step 3: Signal–resampling term.}
We now study the second term, which captures the variation in
$\E[\boldsymbol\lambda\mid S_{\samp}]$ induced by resampling the
signals $S_1,\dots,S_N$.

Recall that the population signal covariance is $M=\E[SS^\top]$ and that
\[
\begin{aligned}
&\widehat M \;=\; \frac1N S_{\samp}S_{\samp}^\top
\;=\; M \;+\; \Delta_M,\\
&\Delta_M:=\widehat M - M
=\frac{1}{N}\sum_{i=1}^N\Bigl(S_i S_i^\top-\E[SS^\top]\Bigr).
\end{aligned}
\]
Define the perturbed ``population'' matrix
\[
\SigmaS
=\Sigma+\widehat M
=\widetilde\Sigma + \Delta_M.
\]
Let $(\xi_k,\popsi_k)$ be the outlier eigenpairs of
$\widetilde\Sigma=\Sigma+M$. Since the outliers are isolated with a fixed eigengap,
the first–order eigenvalue perturbation formula (Lemma \ref{lem:hadamard})
gives
\begin{equation}\label{eq:delta-sigma-sig}
\pssig_k-\xi_k
=\popsi_k^\top\Delta_M\,\popsi_k + R_k(\Delta_M)
=:\delta\xi_k^{(\mathrm{sig})},
\end{equation}
with remainder $|R_k(\Delta_M)|\le C\|\Delta_M\|_{\op}^2$ for some constant
$C>0$. Since $\Delta_M$ is an average of i.i.d.\ mean–zero matrices with finite
fourth moment, we have $\|\Delta_M\|_{\op}=O_{\mathbb P}(N^{-1/2})$, and therefore
$R_k(\Delta_M)=o_{\mathbb P}(N^{-1/2})$. Thus
\[
\delta\xi_k^{(\mathrm{sig})}
=\popsi_k^\top\Delta_M\,\popsi_k+o_{\mathbb P}(N^{-1/2}),\qquad k=1,\dots,r.
\]
Let
\[
\delta\boldsymbol\xi^{(\mathrm{sig})}
:=\bigl(\delta\xi_1^{(\mathrm{sig})},\dots,
        \delta\xi_r^{(\mathrm{sig})}\bigr)^\top.
\]

Next, since $\theta$ is $C^1$ in a neighborhood of each $\xi_k$, we may
apply a first–order Taylor expansion around the population spikes
$\boldsymbol\xi$:
\begin{align*}
\E[\boldsymbol\lambda\mid S_{\samp}]
&=\theta\!\bigl(\boldsymbol\pssig\bigr)+o_{\mathbb P}(N^{-1/2})\\
&=\theta(\boldsymbol\xi)
 +\nabla\theta(\boldsymbol\xi)\,
  \delta\boldsymbol\xi^{(\mathrm{sig})}
 +o_{\mathbb P}(N^{-1/2}),
\end{align*}
where $\nabla\theta(\boldsymbol\xi)
=\diag(\theta'_1,\dots,\theta'_r)$ with
$\theta'_k:=\theta'(\xi_k)$. Hence
\begin{equation}\label{eq:cond-mean-lambda}
\E[\boldsymbol\lambda\mid S_{\samp}]
=\theta(\boldsymbol\xi)
 +\nabla\theta(\boldsymbol\xi)\,
  \delta\boldsymbol\xi^{(\mathrm{sig})}
 +o_{\mathbb P}(N^{-1/2}),
\end{equation}
and
\begin{equation}\label{eq:Cov-second-term}
\Cov\!\bigl(\E[\boldsymbol\lambda\mid S_{\samp}]\bigr)
=\nabla\theta(\boldsymbol\xi)\,
 \Cov\!\bigl(\delta\boldsymbol\xi^{(\mathrm{sig})}\bigr)\,
 \nabla\theta(\boldsymbol\xi)^\top
 +o\!\Bigl(\frac{1}{N}\Bigr).
\end{equation}

We now express $\Cov\!\bigl(\delta\boldsymbol\xi^{(\mathrm{sig})}\bigr)$
in terms of the signal–resampling covariance derived previously.
By Lemma \ref{lem:signal-resampling-general-new} \ (applied with $v_k=\popsi_k$),
\[
\Cov\!\bigl(\popsi_k^\top\Delta_M\popsi_k,\ \popsi_j^\top\Delta_M\popsi_j\bigr)
= \Gamma^{(\mathrm{sig})}_{kj}+o\!\Bigl(\frac{1}{N}\Bigr),
\]
where
\[
\Gamma^{(\mathrm{sig})}_{kj}
=\frac{1}{N}\left\{\E\Big[(\popsi_k^\top S)^2(\popsi_j^\top S)^2\Big]
-\E\big[(\popsi_k^\top S)^2\big]\E\big[(\popsi_j^\top S)^2\big]\right\}.
\]
Collecting entries, we obtain
\begin{equation}\label{eq:Gamma-sig-final}
\Cov\!\bigl(\delta\boldsymbol\xi^{(\mathrm{sig})}\bigr)
=\Gamma^{(\mathrm{sig})}+o\!\Bigl(\frac{1}{N}\Bigr).
\end{equation}
Substituting \eqref{eq:Gamma-sig-final} into \eqref{eq:Cov-second-term}, we arrive at
\begin{equation}\label{eq:Var-second-term-final}
\Cov\!\bigl(\E[\boldsymbol\lambda\mid S_{\samp}]\bigr)
=\nabla\theta(\boldsymbol\xi)\,
 \Gamma^{(\mathrm{sig})}\,
 \nabla\theta(\boldsymbol\xi)^\top
+o\!\Bigl(\frac{1}{N}\Bigr).
\end{equation}

By the multivariate central limit theorem applied to the i.i.d.\ vectors
\[
\left(
(\popsi_1^\top S_i)^2-\mathbb E[(\popsi_1^\top S)^2],
\dots,
(\popsi_r^\top S_i)^2-\mathbb E[(\popsi_r^\top S)^2]
\right)^\top,
\]
we have
\[
\sqrt N\,
\delta\boldsymbol\xi^{(\mathrm{sig})}
\ \xrightarrow{d}\
\mathcal N_r\!\left(
0,\,
N\Gamma^{(\mathrm{sig})}
\right).
\]
Moreover, since the noise matrix is independent of $S_{\samp}$, the
conditional noise fluctuation and the signal-resampling fluctuation converge
jointly to independent Gaussian limits. Therefore,
\[
\sqrt N\bigl(
\boldsymbol\lambda-\theta(\boldsymbol\xi)
\bigr)
=
\sqrt N\bigl(
\boldsymbol\lambda-\theta(\boldsymbol\pssig)
\bigr)
+
\dot\Theta\sqrt N\bigl(
\boldsymbol\pssig-\boldsymbol\xi
\bigr)
+
o_{\mathbb P}(1),
\]
and the claimed unconditional joint CLT follows.

\medskip
\noindent\textbf{Step 4: Putting the pieces together.}
Combining \eqref{eq:Var-first-term} and \eqref{eq:Var-second-term-final}, we obtain
\[
\Var(\boldsymbol\lambda)
=\frac{1}{N}\,\E\!\bigl[V^{\rm cond}\bigr]
+\nabla\theta(\boldsymbol\xi)\,
 \Gamma^{(\mathrm{sig})}\,
 \nabla\theta(\boldsymbol\xi)^\top
+o\!\Bigl(\frac{1}{N}\Bigr),
\]
which is exactly the variance decomposition stated in
Proposition \ref{pro:joint-clt-uncond}.
Equivalently, this yields the covariance $V_*$ of
$\sqrt{N}\bigl(\boldsymbol\lambda-\theta(\boldsymbol\xi)\bigr)$ as
\[
V_*
=\E\!\bigl[V^{\rm cond}\bigr]
+ N\,\nabla\theta(\boldsymbol\xi)\,
    \Gamma^{(\mathrm{sig})}\,
    \nabla\theta(\boldsymbol\xi)^\top,
\]
since $\Var(\boldsymbol\lambda)=V_*/N+o(1/N)$.
This completes the proof.
\end{proof}

\subsection{Proof of Section \ref{Deltamethod}}
\begin{proof}
We prove the stated CLT for $\widehat d^{\,2}$ and the variance formula
\eqref{eq:vard2}, and then derive the multivariate covariance expression
\eqref{eq:SigmaPi-mult}.

\medskip
\noindent\textbf{Step 1: CLT for the sample outliers $\boldsymbol\lambda$.}
By Proposition \ref{pro:joint-clt-uncond}, under \MainAsmpNoise \ and
\MainAsmpSignal, for the vector of outlier sample eigenvalues
$\boldsymbol\lambda=(\lambda_1,\dots,\lambda_r)^\top$ we have
\begin{equation}\label{eq:clt-lambda-again}
  \sqrt N\bigl(\boldsymbol\lambda - \theta(\boldsymbol\xi)\bigr)
  \ \Rightarrow\ \mathcal N_r(0,V_*),
\end{equation}
where $\boldsymbol\xi=(\xi_1,\dots,\xi_r)^\top$ collects the
population outliers of $\widetilde\Sigma:=\Sigma+M$ and
$V_*$ is the limiting covariance matrix as in
Proposition \ref{pro:joint-clt-uncond}.

\medskip
\noindent\textbf{Step 2: Smooth mappings and their derivatives.}
Recall that for $s\notin\{\Sigma_{11},\dots,\Sigma_{pp}\}$ we define
\[
\begin{aligned}
g(s)\;=\;\frac{1}{p}\sum_{i=1}^p\frac{1}{\Sigma_{ii}-s},\\
\qquad
s_2(s)\;=\;\sum_{i=1}^p\frac{1}{(\Sigma_{ii}-s)^2},\\
\qquad
\varphi(s)\;=\;-\frac{1}{g(s)}.
\end{aligned}
\]
Then $g$ and $\varphi$ are smooth on any compact set avoiding the spectrum of
$\Sigma$. Differentiating $g$ we obtain
\[
\begin{aligned}
&g'(s)
=\frac{{\rm d}}{{\rm d}s}
\Bigl[\frac{1}{p}\sum_{i=1}^p(\Sigma_{ii}-s)^{-1}\Bigr]\\
&=\frac{1}{p}\sum_{i=1}^p(\Sigma_{ii}-s)^{-2}\\
&=\frac{1}{p}\,s_2(s).
\end{aligned}
\]
By the chain rule,
\[
\varphi'(s)
=\frac{{\rm d}}{{\rm d}s}\Bigl(-\frac{1}{g(s)}\Bigr)
=\frac{g'(s)}{g(s)^2}
=\frac{s_2(s)}{p\,g(s)^2}.
\]

Next, define the componentwise mappings
\[
\Theta:\R^r\to\R^r,\quad
\Theta(\boldsymbol\xi)
:=\theta(\boldsymbol\xi)
:=\bigl(\theta(\xi_1),\dots,\theta(\xi_r)\bigr)^\top,
\]
\[
\Phi:\R^r\to\R^r,\quad
\Phi(\boldsymbol\xi)
:=\bigl(\varphi(\xi_1),\dots,\varphi(\xi_r)\bigr)^\top
= d^{2},
\]
and $H:=\Phi\circ\Theta^{-1}$. In particular,
\[
H(\boldsymbol\lambda)
=\Phi\bigl(\Theta^{-1}(\boldsymbol\lambda)\bigr)
=\bigl(\varphi(\theta^{-1}(\lambda_1)),\dots,\varphi(\theta^{-1}(\lambda_r))\bigr)^\top.
\]
By construction, our estimator of the signal strengths is exactly
\[
\widehat d^{\,2}
=\bigl(\widehat d_1^{\,2},\dots,\widehat d_r^{\,2}\bigr)^\top
=H(\boldsymbol\lambda),\qquad
d^{2}=\Phi(\boldsymbol\xi)=H\bigl(\theta(\boldsymbol\xi)\bigr).
\]

By assumption, for each $j$ we have $\theta'(\xi_j)\neq 0$, so $\theta$ is
locally invertible near $\xi_j$ and $\theta^{-1}$ is differentiable near
$\theta(\xi_j)$. Since $\theta$ acts componentwise, its Jacobian at
$\boldsymbol\xi$ is the diagonal matrix
\[
\mathrm D\Theta(\boldsymbol\xi)
=\diag\bigl(\theta'(\xi_1),\dots,\theta'(\xi_r)\bigr),
\]
and hence by the inverse function theorem,
\[
\mathrm D\Theta^{-1}\bigl(\theta(\boldsymbol\xi)\bigr)
=\bigl(\mathrm D\Theta(\boldsymbol\xi)\bigr)^{-1}
=\diag\Bigl(\frac{1}{\theta'(\xi_1)},\dots,\frac{1}{\theta'(\xi_r)}\Bigr).
\]
Similarly, since $\Phi$ is also componentwise,
\[
\begin{aligned}
\mathrm D\Phi(\boldsymbol\xi)
=\diag\bigl(\varphi'(\xi_1),\dots,\varphi'(\xi_r)\bigr)\\
=\diag\Bigl(\frac{s_2(\xi_1)}{p\,g(\xi_1)^2},\dots,
            \frac{s_2(\xi_r)}{p\,g(\xi_r)^2}\Bigr).
\end{aligned}
\]

Therefore, the Jacobian of $H=\Phi\circ\Theta^{-1}$ at the point
$\theta(\boldsymbol\xi)$ is
\[
\mathrm D H\bigl(\theta(\boldsymbol\xi)\bigr)
=\mathrm D\Phi(\boldsymbol\xi)\;
 \mathrm D\Theta^{-1}\bigl(\theta(\boldsymbol\xi)\bigr)
=\diag(\Gamma_1,\dots,\Gamma_r)
=: \Gamma,
\]
where for each $j$,
\begin{align*}
\Gamma_j
&=\Bigl[\mathrm D H\bigl(\theta(\boldsymbol\xi)\bigr)\Bigr]_{jj}
=\varphi'(\xi_j)\cdot\frac{1}{\theta'(\xi_j)}\\
&=\frac{s_2(\xi_j)}{p\,g(\xi_j)^2}\cdot\frac{1}{\theta'(\xi_j)},
\end{align*}
which is exactly the expression stated in the proposition.

\medskip
\noindent\textbf{Step 3: Application of the multivariate delta method.}
We now treat $H$ as a smooth mapping in a neighborhood of
$\theta(\boldsymbol\xi)$. Combining \eqref{eq:clt-lambda-again} with the
multivariate delta method yields
\[
\sqrt{N}\Bigl(H(\boldsymbol\lambda)
              -H\bigl(\theta(\boldsymbol\xi)\bigr)\Bigr)
=\sqrt{N}\bigl(\widehat d^{\,2}-d^{2}\bigr)
\ \Rightarrow\ \mathcal N_r\bigl(0,\;\Gamma\,V_*\,\Gamma\bigr).
\]
This proves the claimed CLT
\[
\sqrt{N}\bigl(\widehat d^{\,2}-d^{2}\bigr)\ \Rightarrow\
\mathcal N_r\!\bigl(0,\;\Gamma\,V_*\,\Gamma\bigr),
\]
and shows that the limiting covariance matrix is exactly
$\Gamma\,V_*\,\Gamma$.

\medskip
\noindent\textbf{Step 4: Componentwise variance formula.}
Taking the $j$th diagonal element in the covariance $\Gamma V_* \Gamma$ gives
\[
\Var\!\bigl(\widehat d_j^{\,2}\bigr)
=\Gamma_j^2\,\Var(\lambda_j)\;+\;o\!\Bigl(\frac{1}{N}\Bigr),
\]
where we used that $V_*$ is the covariance of
$\sqrt{N}(\boldsymbol\lambda-\theta(\boldsymbol\xi))$.
Substituting
\[
\Gamma_j
=\frac{s_2(\xi_j)}{p\,g(\xi_j)^2}\cdot\frac{1}{\theta'(\xi_j)}
\]
yields
\[
\Var\!\bigl(\widehat d_j^{\,2}\bigr)
=\Bigl(\frac{s_2(\xi_j)}{p\,g(\xi_j)^2}\Bigr)^{\!2}
\frac{\Var(\lambda_j)}{(\theta'(\xi_j))^{2}}\;+\;o\!\Bigl(\frac{1}{N}\Bigr),
\]
which is exactly the variance expression \eqref{eq:vard2}.

\medskip
\noindent\textbf{Step 5: Matrix form for $\Sigma_{d^2}$.}
Define the asymptotic covariance
\[
\Sigma_{d^2}
:=
\Gamma V_*\Gamma.
\]
Equivalently, the finite-sample covariance satisfies
\[
\Cov(\widehat d^{\,2})
=
\frac{1}{N}\Sigma_{d^2}
+
o\!\left(\frac{1}{N}\right).
\]
\medskip
\noindent\textbf{Step 6: CLT and covariance for the normalized proportions $\Pi$.}
Let $a=\mathbf 1^\top d^2>0$ and define the mapping
\[
F:\R_+^r\to\Delta^{r-1},\qquad
F(d^2)=\Pi:=\frac{d^2}{\mathbf 1^\top d^2}.
\]
Its Jacobian at $d^2$ is the matrix $J^\Pi=(J^\Pi_{kj})_{k,j\le r}$ with
\[
J^\Pi_{kj}
=\frac{\partial \Pi_k}{\partial d_j^2}
=\frac{\delta_{kj}\,a-d_k^2}{a^2},
\qquad a=\mathbf 1^\top d^2,
\]
as in formula\eqref{eq:varJ}. Since we already have the CLT
\[
\sqrt N\bigl(\widehat d^{\,2}-d^2\bigr)
\Rightarrow
\mathcal N_r(0,\Sigma_{d^2}),
\qquad
\Sigma_{d^2}=\Gamma V_*\Gamma.
\]
an application of the multivariate delta method to the smooth map $F$ yields
\[
\sqrt N\bigl(\widehat\Pi-\Pi\bigr)
=\sqrt N\bigl(F(\widehat d^{\,2})-F(d^{2})\bigr)
\ \Rightarrow\
\mathcal N_r\bigl(0,\Sigma_\Pi\bigr)
\]
on the tangent space
$T_\Pi\Delta^{r-1}=\{u\in\R^r:\mathbf 1^\top u=0\}$, where
\begin{equation}
\label{eq:SigmaPi-mult-again}
\Sigma_\Pi
=
J^\Pi\Sigma_{d^2}(J^\Pi)^\top
=
J^\Pi\Gamma V_*\Gamma(J^\Pi)^\top.
\end{equation}
Consequently,
\[
\Cov(\widehat\Pi)
=
\frac{1}{N}\Sigma_\Pi
+
o\!\left(\frac{1}{N}\right).
\]
In particular, $\Sigma_\Pi\mathbf 1=\mathbf 0$ and $\rank(\Sigma_\Pi)\le r-1$.
This is exactly the covariance matrix stated in formula~\eqref{eq:SigmaPi-mult}.
\end{proof}
\begin{proof}
We first justify the componentwise confidence intervals for a single
profile and for the difference between two independent profiles.

\medskip
\noindent\emph{1. Confidence intervals for single-sample components.}
By the profile CLT~\eqref{eq:Pi-CLT}, for each fixed $t$ such that
$[\Sigma_\Pi]_{tt}>0$,
\[
\sqrt{N}\,
(\widehat\Pi_t-\Pi_t)
\xrightarrow{d}
\mathcal N\!\left(
0,[\Sigma_\Pi]_{tt}
\right).
\]
Equivalently,
\[
\frac{
\widehat\Pi_t-\Pi_t
}{
\sqrt{[\Sigma_\Pi]_{tt}/N}
}
\xrightarrow{d}
\mathcal N(0,1).
\]

Let $\widetilde V_\Pi$ be a covariance estimator satisfying the
operator-norm consistency condition stated in the main text:
\[
\left\|
N\widetilde V_\Pi-\Sigma_\Pi
\right\|_{\mathrm{op}}
\xrightarrow{\mathbb P}0.
\]
In particular,
\[
N[\widetilde V_\Pi]_{tt}
\xrightarrow{\mathbb P}
[\Sigma_\Pi]_{tt},
\]
and hence
\[
\frac{
[\widetilde V_\Pi]_{tt}
}{
[\Sigma_\Pi]_{tt}/N
}
\xrightarrow{\mathbb P}1.
\]
Slutsky's theorem therefore gives
\[
\frac{
\widehat\Pi_t-\Pi_t
}{
\sqrt{[\widetilde V_\Pi]_{tt}}
}
\xrightarrow{d}
\mathcal N(0,1).
\]
Consequently,
\[
\Pr\!\left(
\Pi_t\in\mathrm{CI}_t^{(\Pi)}
\right)
\longrightarrow
1-\alpha.
\]

\medskip
\noindent\emph{2. Confidence intervals for profile differences.}
Consider two independent datasets and define
\[
\Delta\widehat\Pi
:=
\widehat\Pi_{r,1}-\widehat\Pi_{r,2},
\qquad
\Delta\Pi
:=
\Pi_r(\rho_1)-\Pi_r(\rho_2).
\]
Let
\[
V_\Delta
:=
\frac{\Sigma_{\Pi,1}}{N_1}
+
\frac{\Sigma_{\Pi,2}}{N_2}
\]
denote the asymptotic covariance of
$\Delta\widehat\Pi-\Delta\Pi$ at its natural finite-sample scale.
The two independent profile CLTs imply that, for each fixed $t$ such
that $[V_\Delta]_{tt}>0$,
\[
\frac{
(\Delta\widehat\Pi)_t-(\Delta\Pi)_t
}{
\sqrt{[V_\Delta]_{tt}}
}
\xrightarrow{d}
\mathcal N(0,1).
\]

Let
\[
\widetilde V_\Delta
:=
\widetilde V_\Pi^{(1)}
+
\widetilde V_\Pi^{(2)},
\]
where the individual covariance estimators satisfy
\[
\left\|
N_i\widetilde V_\Pi^{(i)}-\Sigma_{\Pi,i}
\right\|_{\mathrm{op}}
\xrightarrow{\mathbb P}0,
\qquad i=1,2.
\]
It follows that
\[
\left\|
\widetilde V_\Delta-V_\Delta
\right\|_{\mathrm{op}}
=
o_{\mathbb P}
\left(
N_1^{-1}+N_2^{-1}
\right).
\]
Therefore, for each fixed $t$ with a positive limiting variance,
\[
\frac{
[\widetilde V_\Delta]_{tt}
}{
[V_\Delta]_{tt}
}
\xrightarrow{\mathbb P}1.
\]
Another application of Slutsky's theorem yields
\[
\frac{
(\Delta\widehat\Pi)_t-(\Delta\Pi)_t
}{
\sqrt{[\widetilde V_\Delta]_{tt}}
}
\xrightarrow{d}
\mathcal N(0,1).
\]
Consequently,
\[
\Pr\!\left(
(\Delta\Pi)_t
\in
\mathrm{CI}_t^{(\Delta)}
\right)
\longrightarrow
1-\alpha.
\]
\end{proof}

\begin{proof}
Let
\[
\Delta\Pi
:=
\Pi_r(\rho_1)-\Pi_r(\rho_2),
\qquad
\Delta\widehat\Pi
:=
\widehat\Pi_{r,1}-\widehat\Pi_{r,2},
\]
and write
\[
d
:=
\|\Delta\Pi\|_2,
\qquad
\widehat d
:=
\|\Delta\widehat\Pi\|_2.
\]
By assumption, $d$ is bounded away from zero. Hence the map
$h:\mathbb R^r\to\mathbb R$ defined by $h(x)=\|x\|_2$ is
differentiable at $\Delta\Pi$, with gradient
\[
\nabla h(\Delta\Pi)
=
\frac{\Delta\Pi}{d}.
\]

Define
\[
V_\Delta
:=
\frac{\Sigma_{\Pi,1}}{N_1}
+
\frac{\Sigma_{\Pi,2}}{N_2}.
\]
The two-sample profile CLT and the multivariate delta method give
\[
\frac{
\widehat d-d
}{
\sqrt{v_d}
}
\xrightarrow{d}
\mathcal N(0,1),
\]
where
\[
v_d
:=
\nabla h(\Delta\Pi)^\top
V_\Delta
\nabla h(\Delta\Pi)
=
\frac{
\Delta\Pi^\top
V_\Delta
\Delta\Pi
}{
d^2
}.
\]

Let
\[
\widetilde V_\Delta
:=
\widetilde V_\Pi^{(1)}
+
\widetilde V_\Pi^{(2)}
\]
and define
\[
\widetilde v_d
:=
\frac{
\Delta\widehat\Pi^\top
\widetilde V_\Delta
\Delta\widehat\Pi
}{
\widehat d^{\,2}
}.
\]
By the consistency of $\Delta\widehat\Pi$, we have
\[
\Delta\widehat\Pi
\xrightarrow{\mathbb P}
\Delta\Pi,
\qquad
\widehat d
\xrightarrow{\mathbb P}
d.
\]
Moreover, the covariance-consistency conditions imply
\[
\left\|
\widetilde V_\Delta-V_\Delta
\right\|_{\mathrm{op}}
=
o_{\mathbb P}
\left(
N_1^{-1}+N_2^{-1}
\right).
\]
Since the limiting variance of the nMSD estimator is positive, these
relations yield
\[
\frac{\widetilde v_d}{v_d}
\xrightarrow{\mathbb P}1.
\]
Slutsky's theorem therefore gives
\[
\frac{
\widehat d-d
}{
\sqrt{\widetilde v_d}
}
\xrightarrow{d}
\mathcal N(0,1).
\]
By the definition of $\mathrm{CI}_{\textup{nMSD}}$, it follows that
\[
\Pr\!\left(
d_r(\rho_1,\rho_2)
\in
\mathrm{CI}_{\textup{nMSD}}
\right)
\longrightarrow
1-\alpha.
\]
\end{proof}

\section{Proof of \MainSubsecThirtyThree}
\subsection{Proof of \MainThmWaldNull}

\begin{proof}
From the joint profile central limit theorems and the independence of
the two datasets, under $H_0$,
\[
\sqrt{N_{\mathrm{eff}}}\,
\Delta_\Pi
\xrightarrow{d}
\mathcal N_r(0,\Sigma_\Delta),
\]
where
\[
\Sigma_\Delta
=
\alpha\Sigma_{\Pi,1}
+
(1-\alpha)\Sigma_{\Pi,2}.
\]
Since both estimated profiles lie in the simplex,
$\mathbf 1^\top\Delta_\Pi=0$.
Moreover, normalization gives
$\Sigma_\Delta\mathbf 1=\mathbf 0$, and the assumed positive definiteness of
$\Sigma_\Delta$ on the simplex tangent space implies
\[
\operatorname{null}(\Sigma_\Delta)
=
\operatorname{span}\{\mathbf 1\},
\qquad
\operatorname{rank}(\Sigma_\Delta)
=
r-1.
\]

Therefore,
\[
\begin{aligned}
T_\Pi^\circ
&=
\bigl(
\sqrt{N_{\mathrm{eff}}}\Delta_\Pi
\bigr)^\top
\Sigma_\Delta^+
\bigl(
\sqrt{N_{\mathrm{eff}}}\Delta_\Pi
\bigr)
\xrightarrow{d}
\chi^2_{r-1}.
\end{aligned}
\]
\end{proof}
\subsection{Proof of \MainThmWaldPower}

\begin{proof}
Under the fixed alternative,
\[
\Delta_\Pi
=
\Delta
+
O_{\mathbb P}
\bigl(
N_{\mathrm{eff}}^{-1/2}
\bigr).
\]
Hence,
\[
\begin{aligned}
\frac{T_\Pi^\circ}{N_{\mathrm{eff}}}
&=
\Delta_\Pi^\top
\Sigma_\Delta^+
\Delta_\Pi=
\Delta^\top
\Sigma_\Delta^+
\Delta
+
o_{\mathbb P}(1)
\xrightarrow{\mathbb P}
c_0.
\end{aligned}
\]
The remaining statements follow from the upper-tail bound for the
$\chi^2_{r-1}$ reference distribution.
\end{proof}
\subsection{Proof of \MainCorSep}
\begin{proof}[Proof of \MainCorSep]
Recall that $\Delta=\Pi_r(\rho_1)-\Pi_r(\rho_2)\in\R^r$ and
\[
  c_0 \;=\; \Delta^\top  \Sigma_\Delta^{+}\Delta,
\]
where $\Sigma_\Delta$ is a symmetric positive–semidefinite matrix
(the asymptotic covariance in the two–sample CLT for $\sqrt{N_{\mathrm{eff}}}\,\Delta_\Pi$)
and $\Sigma_\Delta^{+}$ is its Moore–Penrose pseudoinverse.

\paragraph*{Step 1: lower bound on $c_0$}
Let the spectral decomposition of $ \Sigma_\Delta$ be
\[
   \Sigma_\Delta \;=\; U\Lambda U^\top,
\]
where $U$ is orthogonal and
$\Lambda=\mathrm{diag}(\lambda_1,\ldots,\lambda_q,0,\ldots,0)$ with
$\lambda_1\ge\cdots\ge\lambda_q>0$ the nonzero eigenvalues.
Then
\[
  \Sigma_\Delta^{+} \;=\; U\Lambda^{+}U^\top,
  \qquad
  \Lambda^{+}=\mathrm{diag}\Bigl(\frac1{\lambda_1},\ldots,\frac1{\lambda_q},
  0,\ldots,0\Bigr).
\]

Because each $\Pi_r(\rho_k)$ lies in the simplex,
$\mathbf 1^\top\Delta=0$, so $\Delta\in\mathbf 1^\perp$.
Normalization gives $\Sigma_\Delta\mathbf 1=\mathbf 0$, and the assumed
positive definiteness of $\Sigma_\Delta$ on $\mathbf 1^\perp$ implies
\[
\operatorname{null}(\Sigma_\Delta)
=
\operatorname{span}\{\mathbf 1\}.
\]
Hence $\Delta\in\operatorname{range}(\Sigma_\Delta)$, and therefore its
coordinates $z:=U^\top\Delta$ satisfy $z_i=0$ whenever $\lambda_i=0$.
Then
\[
  c_0
  = \Delta^\top \Sigma_\Delta^{+}\Delta
  = \sum_{i:\lambda_i>0} \frac{z_i^2}{\lambda_i}.
\]
Let $\lambda_{\max}(\Sigma_\Delta)$ denote the largest nonzero eigenvalue,
that is, $\lambda_{\max}( \Sigma_\Delta)=\lambda_1$.
For each $i$ with $\lambda_i>0$ we have
$\lambda_i\le\lambda_{\max}( \Sigma_\Delta)$, hence
\[
  \frac1{\lambda_i}
  \;\ge\; \frac{1}{\lambda_{\max}( \Sigma_\Delta)}.
\]
Therefore
\[
  c_0
  = \sum_{i:\lambda_i>0}\frac{z_i^2}{\lambda_i}
  \;\ge\;
  \frac{1}{\lambda_{\max}( \Sigma_\Delta)}\sum_{i:\lambda_i>0} z_i^2
  = \frac{\|\Delta\|_2^2}{\lambda_{\max}( \Sigma_\Delta)}.
\]
This proves the first inequality.

By Theorem~\MainThmWaldPower,
\[
\frac{T_\Pi^\circ}{N_{\mathrm{eff}}}
\xrightarrow{\mathbb P}
c_0.
\]
Combining this convergence with the preceding lower bound, for every
$\varepsilon\in(0,1)$,
\[
\Pr\left\{
T_\Pi^\circ
\ge
(1-\varepsilon)N_{\mathrm{eff}}
\frac{\|\Delta\|_2^2}
{\lambda_{\max}(\Sigma_\Delta)}
\right\}
\longrightarrow 1.
\]
If $\|\Delta\|_2\ge C>0$, this gives
\[
T_\Pi^\circ
\gtrsim
N_{\mathrm{eff}}
\frac{C^2}
{\lambda_{\max}(\Sigma_\Delta)}
\]
with high probability, as claimed.
\end{proof}

\section{Proof of \MainSubsecKernelMSS}
\begin{proposition}[Operator–norm concentration for kernel covariance]
\label{prop:kernel-op-conc}
Work in the RKHS setup of \MainSubsecKernelMSS \ for a fixed
dataset $k\in\{1,2\}$.
Let $M_{K,k}$ be the centered population kernel covariance operator
and $\widehat M_{K,k}$ its empirical counterpart,
\[
\begin{aligned}
  &M_{K,k}
  := \E_{X\sim\rho_k}\Bigl[(\phi(X)-m_k)\otimes(\phi(X)-m_k)\Bigr],\\
  &\widehat M_{K,k}
  := \frac{1}{N_k}\sum_{i=1}^{N_k}
      \bigl(\phi(x_{k,i})- m_k\bigr)\otimes
      \bigl(\phi(x_{k,i})- m_k\bigr),
\end{aligned}
\]
with $m_k:=\E_{X\sim\rho_k}[\phi(X)]$.
Assume $\sup_{x\in\X}K(x,x)\le\kappa^2<\infty$.
Then there exist constants $c_1,c_2>0$ (depending only on $\kappa$) such that
for all $\delta\in(0,1)$,
\[
  \Pr\!\left(
    \bigl\|\widehat M_{K,k}-M_{K,k}\bigr\|_{\op}
    >
    \frac{c_1\kappa^2\log(c_2/\delta)}{\sqrt{N_k}}
  \right)
  \;\le\;\delta .
\]
In particular,
\[
  \bigl\|\widehat M_{K,k}-M_{K,k}\bigr\|_{\op}
  = O_{\mathbb P}\!\bigl(N_k^{-1/2}\bigr).
\]
\end{proposition}

\begin{proof}
We recall a standard concentration inequality for Hilbert–space valued
random variables; see, for example, Smale and Zhou (2009)\citep{SmaleZhou2009Geometry}:

\medskip\noindent
\emph{Hilbert–space concentration.}
Let $(\mathcal H,\langle\cdot,\cdot\rangle)$ be a Hilbert space and let
$\zeta_1,\dots,\zeta_{N_k}$ be i.i.d.\ random elements in $\mathcal H$ such that
$\|\zeta_i\|\le B$ almost surely and
$\E[\zeta_i]=\mu$.
Then there exist universal constants $a_1,a_2>0$ such that for all
$\delta\in(0,1)$,
\[
  \Pr\!\left(
    \left\|
      \frac{1}{N_k}\sum_{i=1}^{N_k}(\zeta_i-\mu)
    \right\|>
    \frac{a_1B\log(a_2/\delta)}{\sqrt{N_k}}
  \right)
  \le\delta .
\tag{$\star$}
\label{eq:HS-conc}
\]

We apply \eqref{eq:HS-conc} in the Hilbert space
$HS(\mathcal H_K)$ of Hilbert–Schmidt operators on $\mathcal H_K$ with inner
product $\langle A,B\rangle_{HS}=\operatorname{tr}(B^*A)$
and norm $\|\cdot\|_{HS}$.
Define
\[
  \zeta_i
  :=
  \bigl(\phi(x_{k,i})- m_k\bigr)\otimes
  \bigl(\phi(x_{k,i})- m_k\bigr),
  \qquad i=1,\dots,N_k.
\]
By boundedness of the kernel we have
\[
  \|\phi(x)\|_{\mathcal H_K}^2
  = K(x,x)\le\kappa^2
  \quad\Rightarrow\quad
  \|\phi(x)-m_k\|_{\mathcal H_K}\le 2\kappa,
\]
so almost surely we obtain
\[
  \|\zeta_i\|_{HS}
  = \|\phi(x_{k,i})-m_k\|_{\mathcal H_K}^2
  \le (2\kappa)^2 = 4\kappa^2
\]
Thus we may take $B=4\kappa^2$ in \eqref{eq:HS-conc}.

Next note that
\[
\begin{aligned}
&  \frac{1}{N_k}\sum_{i=1}^{N_k}\zeta_i
  = \widehat M_{K,k},\\
&  \mu
  = \E_{X\sim\rho_k}
      \bigl[(\phi(X)-m_k)\otimes(\phi(X)-m_k)\bigr]
  = M_{K,k}.
\end{aligned}
\]
Therefore,
\[
  \left\|\frac{1}{N_k}\sum_{i=1}^{N_k}(\zeta_i-\mu)\right\|_{HS}
  = \|\widehat M_{K,k}-M_{K,k}\|_{HS}.
\]
Applying \eqref{eq:HS-conc} and using
$\|A\|_{\op}\le\|A\|_{HS}$ for any Hilbert–Schmidt operator $A$, we obtain
\[
  \Pr\!\left(
    \bigl\|\widehat M_{K,k}-M_{K,k}\bigr\|_{\op}
    >
    \frac{a_1(4\kappa^2)\log(a_2/\delta)}{\sqrt{N_k}}
  \right)
  \le\delta .
\]
Renaming constants $c_1:=4a_1$ and $c_2:=a_2$ gives the stated inequality,
and the $O_{\mathbb P}(N_k^{-1/2})$ rate follows immediately.
\end{proof}

\begin{proposition}[Convergence of kernel principal variances]
\label{prop:kernel-eig-conv}
In the setting of \MainSubsecKernelMSS, fix $k\in\{1,2\}$.
Let $\{\mathcal D_{k,j}\}_{j\ge1}$ be the eigenvalues of $M_{K,k}$ arranged in
nonincreasing order and let
$\hat{\mathcal D}_{k,1}\ge\cdots\ge\hat{\mathcal D}_{k,r}>0$ be the leading $r$
eigenvalues of the normalized centered Gram matrix
\[
  G_k := \frac{1}{N_k}\,\widetilde K_k .
\]
Assume the eigengap condition
\[
\begin{aligned}
 & \mathcal D_{k,1}\ge\cdots\ge\mathcal D_{k,r}>\mathcal D_{k,r+1},\\
 & \gamma_k := \min_{1\le j\le r}
    \bigl(\mathcal D_{k,j}-\mathcal D_{k,j+1}\bigr) > 0 .
\end{aligned}
\]
Then, under the bounded-kernel assumption $\sup_x K(x,x)\le\kappa^2$,
for each fixed $1\le j\le r$,
\[
  \hat{\mathcal D}_{k,j}
  \;\xrightarrow{\ \mathbb P\ }\;
  \mathcal D_{k,j},
  \quad\text{as }N_k\to\infty.
\]
\end{proposition}
\begin{proof}
We work with the centered feature map associated with dataset $k$.
Recall from \MainSubsecKernelMSS \ that
\[
\begin{aligned}
 & M_{K,k}
  := \E_{X\sim\rho_k}\Bigl[(\phi(X)-m_k)\otimes(\phi(X)-m_k)\Bigr],\\
  &m_k := \E_{X\sim\rho_k}[\phi(X)],
\end{aligned}
\]
so the covariance is taken with respect to the centered feature
$\phi(X)-m_k$.

\medskip\noindent
\textbf{Step 1: Centered feature map and Gram matrix.}
Define the centered feature map
\[
  \tilde\phi_k(x) \;:=\; \phi(x)-m_k,\qquad x\in\X,
\]
and note that
\[
  \langle \tilde\phi_k(x_{k,i}),\tilde\phi_k(x_{k,j})\rangle_{\mathcal H_K}
  \;=\; [\widetilde K_k]_{ij},
\]
so the centered Gram matrix can be written as
\[
  \widetilde K_k
  = \bigl(\langle \tilde\phi_k(x_{k,i}),\tilde\phi_k(x_{k,j})\rangle_{\mathcal H_K}\bigr)_{1\le i,j\le N_k}.
\]
By definition,
\[
  G_k \;=\; \frac{1}{N_k}\,\widetilde K_k.
\]

\medskip\noindent
\textbf{Step 2: Sampling operator representation.}
Define the sampling operator
$S_{c,k}:\mathcal H_K\to\R^{N_k}$ by
\[
  (S_{c,k} f)_i
  := \langle f,\tilde\phi_k(x_{k,i})\rangle_{\mathcal H_K},
  \qquad i=1,\dots,N_k.
\]
Its adjoint $S_{c,k}^*:\R^{N_k}\to\mathcal H_K$ satisfies
\[
  S_{c,k}^* c
  = \sum_{i=1}^{N_k} c_i\,\tilde\phi_k(x_{k,i}),
  \qquad c=(c_1,\dots,c_{N_k})^\top\in\R^{N_k}.
\]

A direct calculation shows
\[
  S_{c,k}^* S_{c,k}
  = \sum_{i=1}^{N_k} \tilde\phi_k(x_{k,i})\otimes\tilde\phi_k(x_{k,i})
  \;=\; N_k\,\widehat M_{K,k},
\]
where
\[
  \widehat M_{K,k}
  := \frac{1}{N_k}\sum_{i=1}^{N_k}
     \tilde\phi_k(x_{k,i})\otimes\tilde\phi_k(x_{k,i})
\]
is the empirical covariance operator of the centered features.

Similarly,
\[
  S_{c,k} S_{c,k}^*
  = \Bigl(\langle \tilde\phi_k(x_{k,i}),\tilde\phi_k(x_{k,j})\rangle_{\mathcal H_K}\Bigr)_{i,j}
  = \widetilde K_k,
\]
so
\[
  G_k
  = \frac{1}{N_k}\,S_{c,k}S_{c,k}^*.
\]

\medskip\noindent
\textbf{Step 3: Equality of nonzero eigenvalues.}
Let $\mu_1\ge\cdots\ge\mu_{N_k}\ge0$ be the eigenvalues of
$S_{c,k}S_{c,k}^*$ (in nonincreasing order), and let
$\tilde\mu_1\ge\tilde\mu_2\ge\cdots\ge0$ be the eigenvalues of
$S_{c,k}^*S_{c,k}$.
A basic linear-algebra fact for bounded operators $A$ between Hilbert
spaces (applied here with $A=S_{c,k}$) is that $A^*A$ and $AA^*$ have
the same nonzero eigenvalues, with the same multiplicities.
Thus the nonzero eigenvalues of $S_{c,k}S_{c,k}^*$ and $S_{c,k}^*S_{c,k}$
coincide:
\[
  \{\mu_j:\mu_j>0\}
  = \{\tilde\mu_j:\tilde\mu_j>0\}.
\]

Now, the eigenvalues of $\widehat M_{K,k}$ are
\[
  \lambda_j(\widehat M_{K,k})
  = \frac{1}{N_k}\,\tilde\mu_j,
\]
while the eigenvalues of $G_k=(1/N_k)S_{c,k}S_{c,k}^*$ are
\[
  \lambda_j(G_k)
  = \frac{1}{N_k}\,\mu_j.
\]
By the previous paragraph, the nonzero eigenvalues of $G_k$ and
$\widehat M_{K,k}$ are equal (counting multiplicities). Therefore, when
we arrange them in nonincreasing order,
\[
  \hat{\mathcal D}_{k,j}
  = \lambda_j(G_k)
  = \lambda_j(\widehat M_{K,k}),
  \qquad j\ge1.
\]

\medskip\noindent
\textbf{Step 4: Operator--norm convergence and Weyl's inequality.}
By Proposition~\ref{prop:kernel-op-conc} (applied to the bounded kernel
$K$ and centered feature map $\tilde\phi_k$), we have
\[
  \bigl\|\widehat M_{K,k}-M_{K,k}\bigr\|_{\op}
  = O_{\mathbb P}\!\bigl(N_k^{-1/2}\bigr)
  \;\xrightarrow{\ \mathbb P\ }\; 0.
\]

Since $M_{K,k}$ and $\widehat M_{K,k}$ are self-adjoint, Weyl's
eigenvalue inequality gives, for each $j\ge1$,
\[
  \bigl|\lambda_j(\widehat M_{K,k})-\lambda_j(M_{K,k})\bigr|
  \le \bigl\|\widehat M_{K,k}-M_{K,k}\bigr\|_{\op}
  \xrightarrow{\ \mathbb P\ } 0.
\]
By definition, $\lambda_j(M_{K,k})=\mathcal D_{k,j}$, so for each fixed
$1\le j\le r$,
\[
  \bigl|\lambda_j(\widehat M_{K,k})-\mathcal D_{k,j}\bigr|
  \xrightarrow{\ \mathbb P\ } 0.
\]

Finally, using $\hat{\mathcal D}_{k,j}=\lambda_j(\widehat M_{K,k})$ from
Step~3, we conclude that for each fixed $1\le j\le r$,
\[
  \hat{\mathcal D}_{k,j}
  \xrightarrow{\ \mathbb P\ } \mathcal D_{k,j},
\]
which is exactly the desired convergence.
\end{proof}

\begin{theorem}[Consistency of the kernelized manifold spectral dissimilarity]
\label{thm:kernel-mss-consistency}
Under the RKHS setup of \MainSubsecKernelMSS, assume:
\begin{enumerate}
\item The kernel is bounded: $\sup_{x\in\X}K(x,x)\le\kappa^2<\infty$.
\item For each dataset $k\in\{1,2\}$, the eigenvalues
$\{\mathcal D_{k,j}\}_{j\ge1}$ of $M_{K,k}$ satisfy the spiked/eigengap
condition
\[
\begin{aligned}
  \mathcal D_{k,1}\ge\cdots\ge\mathcal D_{k,r}>\mathcal D_{k,r+1},\\
  \min_{1\le j\le r}
    \bigl(\mathcal D_{k,j}-\mathcal D_{k,j+1}\bigr) > 0.
\end{aligned}
\]
\end{enumerate}
Let
\[
\begin{aligned}
 & \Pi^{(K)}_r(\rho_k)
  := \frac{\bigl(\mathcal D_{k,1},\ldots,\mathcal D_{k,r}\bigr)}
           {\sum_{j=1}^r \mathcal D_{k,j}}
  \in\Delta^{\,r-1},\\
 & d^{(K)}_r(\rho_1,\rho_2)
  := \bigl\|\Pi^{(K)}_r(\rho_1)-\Pi^{(K)}_r(\rho_2)\bigr\|_2 ,
\end{aligned}
\]
and define the empirical quantities
\[
\begin{aligned}
  \widehat\Pi^{(K)}_{r,k}
  := \frac{\bigl(\hat{\mathcal D}_{k,1},\ldots,\hat{\mathcal D}_{k,r}\bigr)}
           {\sum_{j=1}^r \hat{\mathcal D}_{k,j}},\\
  \widehat d^{(K)}_r
  := \bigl\|\widehat\Pi^{(K)}_{r,1}-\widehat\Pi^{(K)}_{r,2}\bigr\|_2 .
\end{aligned}
\]
Then, as $N_1,N_2\to\infty$,
\[
  \widehat\Pi^{(K)}_{r,k}
  \;\xrightarrow{\ \mathbb P\ }\; \Pi^{(K)}_r(\rho_k),
  \qquad k=1,2,
\]
and consequently
\[
  \widehat d^{(K)}_r
  \;\xrightarrow{\ \mathbb P\ }\; d^{(K)}_r(\rho_1,\rho_2).
\]
\end{theorem}

\begin{proof}
Fix $k\in\{1,2\}$.
By Proposition~\ref{prop:kernel-eig-conv}, for each fixed $1\le j\le r$,
\[
  \hat{\mathcal D}_{k,j}\ \xrightarrow{\ \mathbb P\ }\ \mathcal D_{k,j}.
\]
Hence the vector of leading eigenvalues converges in probability:
\[
  \bigl(\hat{\mathcal D}_{k,1},\ldots,\hat{\mathcal D}_{k,r}\bigr)
  \xrightarrow{\ \mathbb P\ }
  \bigl(\mathcal D_{k,1},\ldots,\mathcal D_{k,r}\bigr).
\]
The sum in the denominator is a continuous map, so
\[
  \sum_{j=1}^r \hat{\mathcal D}_{k,j}
  \xrightarrow{\ \mathbb P\ } \sum_{j=1}^r \mathcal D_{k,j} =: s_k > 0.
\]
The normalization map $v\mapsto v/(\mathbf1^\top v)$ is continuous on
$\{v:\mathbf1^\top v>0\}$, and $s_k>0$ by assumption, so the continuous
mapping theorem yields
\[
  \widehat\Pi^{(K)}_{r,k}
  = \frac{\bigl(\hat{\mathcal D}_{k,1},\ldots,\hat{\mathcal D}_{k,r}\bigr)}
         {\sum_{j=1}^r \hat{\mathcal D}_{k,j}}
  \xrightarrow{\ \mathbb P\ }
  \frac{\bigl(\mathcal D_{k,1},\ldots,\mathcal D_{k,r}\bigr)}
         {\sum_{j=1}^r \mathcal D_{k,j}}
  = \Pi^{(K)}_r(\rho_k).
\]

Finally, the Euclidean norm on $\R^r$ is continuous, so another application
of the continuous mapping theorem gives
\[
  \widehat d^{(K)}_r
  = \bigl\|\widehat\Pi^{(K)}_{r,1}-\widehat\Pi^{(K)}_{r,2}\bigr\|_2
  \xrightarrow{\ \mathbb P\ }
  \bigl\|\Pi^{(K)}_r(\rho_1)-\Pi^{(K)}_r(\rho_2)\bigr\|_2
  = d^{(K)}_r(\rho_1,\rho_2),
\]
as claimed.
\end{proof}
\subsection{Pathway enrichment for genes expressed at lower levels in each MSD-defined group}

As a complement to the pathway enrichment analysis of genes expressed at
higher levels in each MSD-defined hepatocyte group, we also performed
over-representation analysis for genes expressed at lower levels in each
group relative to the remaining two groups. Reactome and Gene Ontology
Biological Process (GO BP) terms were analyzed separately using the same
strict MAST marker criterion as in the main analysis.

\begin{figure}[ht]
  \centering
  \includegraphics[width=0.95\linewidth]{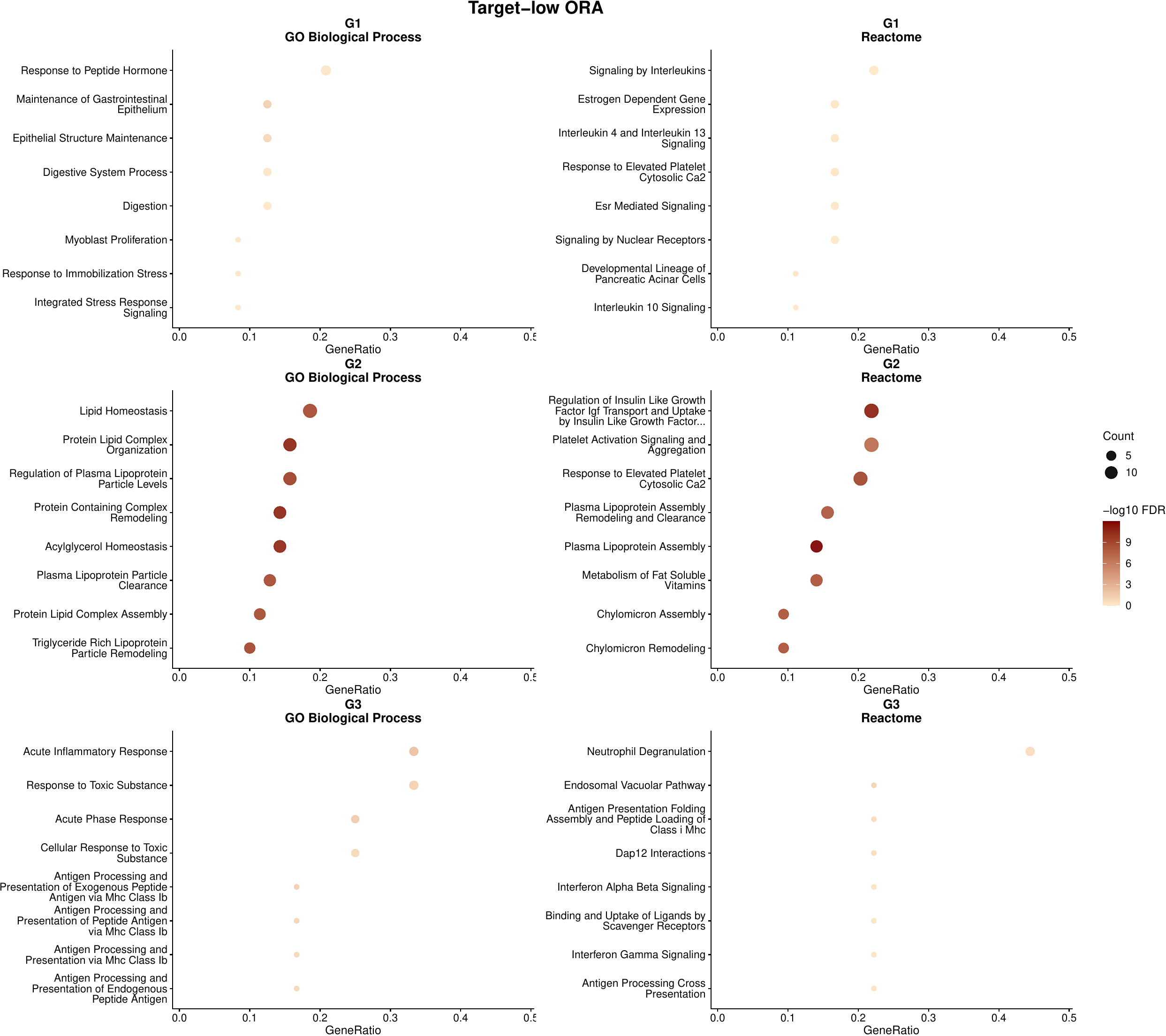}
  \caption{Reactome and GO Biological Process pathway enrichment for genes
  expressed at lower levels in each MSD-defined hepatocyte group. }
  \label{fig:hcc_pathway_low}
\end{figure}

\bibliographystyle{abbrvnat}
\bibliography{reference}